%% file: ASTARS.tex
\documentclass[11pt]{amsart}
\usepackage{amssymb, amsmath, graphicx, caption, enumerate, url}
\graphicspath{ {images/} }
\usepackage{amsthm,amsaddr}
\usepackage{xargs}
\usepackage{scalerel}
\usepackage[ ]{algorithm2e}
\usepackage[dvipsnames]{xcolor}
\usepackage{multicol}
\usepackage{lineno}
\usepackage{amsrefs}

\makeatletter
\renewcommand{\@algocf@capt@plain}{above}
\makeatother

\newcommand{\D}{\mathcal{D}}
\newcommand{\A}{\mathcal{A}}
\newcommand{\N}{\mathbb{N}}
\newcommand{\I}{\mathcal{I}}

\newcommand{\R}{\mathbb{R}}

\newcommand{\ds}{\displaystyle}

\usepackage[colorinlistoftodos,prependcaption,textsize=tiny]{todonotes}
\newcommandx{\question}[2][1=]{\todo[linecolor=red,backgroundcolor=red!25,bordercolor=red,#1]{#2}}
\newcommandx{\change}[2][1=]{\todo[linecolor=blue,backgroundcolor=blue!25,bordercolor=blue,#1]{#2}}
\newcommandx{\add}[2][1=]{\todo[linecolor=OliveGreen,backgroundcolor=OliveGreen!25,bordercolor=OliveGreen,#1]{#2}}
\newcommandx{\improve}[2][1=]{\todo[linecolor=Plum,backgroundcolor=Plum!25,bordercolor=Plum,#1]{#2}}
\newcommandx{\thiswillnotshow}[2][1=]{\todo[disable,#1]{#2}}
\newcommandx{\remove}[2][1=]{\todo[linecolor=yelllow,backgroundcolor=yellow!10,bordercolor=red,#1]{#2}}

\newcommand\reallywidehat[1]{\arraycolsep=0pt\relax%
\begin{array}{c}
\stretchto{
  \scaleto{
    \scalerel*[\widthof{\ensuremath{#1}}]{\kern-.5pt\bigwedge\kern-.5pt}
    {\rule[-\textheight/2]{1ex}{\textheight}} 
  }{\textheight} %
}{0.5ex}\\           
#1\\                 
\rule{-1ex}{0ex}
\end{array}
}

\calclayout

\title[Accelerating Derivative-Free Optimization]{Accelerating Derivative-Free Optimization with Dimension Reduction and Hyperparameter Learning}
\author[J.R. Hall and V. Carey]{Jordan R. Hall and Varis Carey}
\date{\today}

\begin{document}

\begin{abstract}

We consider convex, black-box functions $f$ with additive or multiplicative noise with a high-dimensional parameter space $\Lambda$ and a data space $\mathcal{D}$ of lower dimension, where $\nabla f$ exists, but may be inaccessible.  
We investigate Derivative-Free Optimization (DFO) in this setting and propose a novel method, Active STARS (ASTARS), based on STARS \cite{CW} and dimension reduction in $\Lambda$ via Active Subspace methods \cite{Constantine2015}. 
STARS hyperparmeters are inversely proportional to the known dimension $P$ of $\Lambda$, resulting in heavy smoothing and small step sizes for large $P$.
When possible, ASTARS learns a lower-dimensional Active Subspace, $\mathcal{A} \subset \Lambda$, defining a set of directions in $\Lambda$ causing the majority of the variance in $f$. 
ASTARS iterates are updated with steps only taken in $\mathcal{A}$, reducing the value of $f$ more efficiently than STARS, which updates iterates in the full variables, $\Lambda$. 
In addition to computational savings made by stepping only in $\A$ when it exists, computational costs may be reduced further by estimating hyperparameters and $\A$ using STARS iterates, reducing the total evaluations of $f$ and eliminating the requirement that the user specify hyperparameters, which may be unknown in our setting.
We call this method Fully Automated ASTARS (FAASTARS).
We show that STARS and ASTARS will both converge -- with a certain complexity -- even with inexact, estimated hyperparemters.
We also find that FAASTARS converges with the use of estimated $\A$ and hyperparameters.
We explore the effectiveness of ASTARS and FAASTARS in numerical examples which compare ASTARS and FAASTARS to STARS.

\end{abstract}

\maketitle

\setcounter{tocdepth}{1} 
\tableofcontents
\newpage

\section{Introduction}

In this paper, we present an efficient optimization algorithm for expensive, noisy objective functions which avoids the use of full gradient information.
Efficiency is achieved by performing parameter space dimension reduction whenever possible, which reduces the burden of many computational tasks, including optimization.
Additionally, we present a fully-automated version of our algorithm which estimates its hyperparameters and performs approximate dimension reduction by reusing iterates, saving costly function evaluations.

Optimizing a deterministic function in the presence of its noise is a problem of \textit{Optimization Under Uncertainty} (\textit{OUU}) and quite often arises as a necessary step in important mathematical and statistical applications.
Many objective functions appearing in \textit{Uncertainty Quantification (UQ)} problems and applications involve post-processing evaluations of noisy functions (i.e., physical models). 
For instance, \textit{Stochastic Inverse Problems} (\textit{SIPs}) may be solved by posing equivalent deterministic, convex optimization problems (under the heavy assumptions that the true model $f$ is a linear and distributions are Gaussian). We carefully address this topic and propose the use of our methods for optimization in this setting in a forthcoming paper.

We define a parameter space $\Lambda$ with $P:= \dim \left( \Lambda \right)$,  a data space $\mathcal{D}$, and a  \textit{parameter-to-observable map} or \textit{model} $f: \Lambda \rightarrow \mathcal{D}$, which we assume is convex.
We assume $ \dim \left( \mathcal{D} \right) =: D < P$. 
(We focus on the case $\mathcal{D} \subseteq \mathbb{R}$.)
Points in $\mathcal{D}$ may be known values of $f(\lambda), \lambda\in \Lambda$; we may write $d=f(\lambda)$ to denote the particular datum corresponding to the evaluation of a point $\lambda \in \Lambda$ under $f$. 

Many physical systems of interest possess turbulent or chaotic behavior.  The physical state of a time-dependent system $u(t,\lambda)$ and the corresponding parameter-to-observable map may be modeled as a stochastic process, $\hat{f}(u(t,\lambda))$, a deterministic function with additive or multiplicative noise.  
We model a noisy signal with $\hat{f}(\lambda;\xi)=f(\lambda)+\epsilon(\xi)$ (additive noise) or $\hat{f}(\lambda;\xi)=f(\lambda)(1+\epsilon(\xi))$ (multiplicative noise). 
Now $f$ represents some true signal and $\hat{f}$ represents a noisy, polluted signal. 
In both cases $\epsilon(\cdot)$ denotes a random variable specified by realizations or draws $\epsilon(\xi)$ and we assume $\epsilon(\cdot)$ has bounded variance $\sigma^2<\infty$.
We postulate that $\hat{f}$ is an expensive, black-box function so that $f$ lacks closed form entirely and may take minutes to hours for even a single evaluation.

The efficient and accurate extraction of gradients of $\hat{f}$ in parameter space is a challenging undertaking, as popular techniques based on
linearization, including adjoint methods, are inaccurate \cites{lea2000, Qiqi2014}.  
The finite-difference approximations to $\nabla f$ 
involve  $P=\text{dim}(\Lambda)$ 
additional, usually nonlinear model solves for physical system states $u(t, \lambda + \delta)$, and may be greatly polluted by the noise in $\hat{f}=f+\epsilon$ or $\hat{f}=f(1+\epsilon)$.

As a consequence of these difficulties, optimization in this setting often needs to be performed by algorithms which
do not require gradient approximations. 
Otherwise, we expect to surpass a reasonable computational budget since even one gradient approximation involves $\mathcal{O}(P)$  evaluations of $\hat{f}$ using, for instance, finite differencing.
With $\nabla f$ considered unavailable and unfeasible to approximate we consider \textit{Derivative-Free Optimization (DFO)} techniques which avoid evaluations of $\nabla f$, as well as expensive and possibly inaccurate approximations to full pointwise gradients.
In the following section we provide a discussion of DFO, DFO hyperparameter learning, and parameter space dimension reduction via Active Subspace (AS) methods.

\subsection{Derivative-Free Optimization (DFO)}

As in \cites{CW, Nesterov, ARDFDS}, we consider the additive noise OUU problem

\begin{eqnarray} \label{eq:1}
\min_{\lambda \in \Lambda} \quad \mathbb{E}\left[f(\lambda)+\epsilon(\xi)\right],
\end{eqnarray} 

\noindent and the multiplicative noise OUU problem

\begin{eqnarray} \label{eq:2}
\min_{\lambda \in \Lambda} \quad \mathbb{E}\left[f(\lambda)(1+\epsilon(\xi))\right],
\end{eqnarray} 

\noindent where we assume:

\begin{enumerate}[(i.)]

\item $f: \Lambda=\R^P \to \R=\D$ is continuously differentiable, convex, and $\nabla f$ has a real-valued Lipschitz constant $L_1>0$;

\item $\epsilon(\xi)$ is a random variable with probability density $\pi_{\epsilon}(\epsilon(\xi))$;

\item for all $\lambda$ and $\xi$ the noise $\epsilon(\xi)$ is independent and identically distributed, has bounded variance $\sigma^2$, and is unbiased; i.e., $\mathbb{E}_\xi (\epsilon(\xi))=0$;

\item for multiplicative OUU, the additional assumptions that: the \textit{signal-to-noise ratio} is bounded -- i.e., $\mathbb{E}((1+\epsilon(\xi))^{-1})<b$, $b>0$ -- and the support of $\epsilon(\xi)$ is bounded by $\pm a$, $a<1.$

\end{enumerate}

We consider DFO algorithms suited for additive and multiplicative noise. DFO algorithms require nothing more than the ability to evaluate the noisy model and randomly draw from a normal distribution; $\nabla f$ is not required.

Given an initial iterate $\lambda^{(0)}$, many DFO algorithms find subsequent iterates $\lambda^{(k)}$ by random walks or \textit{ball walks} in $\Lambda$ specified by draws from a $P$-dimensional Gaussian distribution.
Iterates are often controlled by prescribed hyperparameters (e.g., step size) which depend on potentially unknown properties of the model.
For example, in the DFO algorithm considered in this paper \cite{CW}, the hyperparameters depend on the variance in the noise, $\sigma^2$, and on the first-degree or gradient Lipschitz constant of $f$, denoted $L_1$, which we assume to exist; i.e., the real number $L_1> 0: ||\nabla f(\lambda^1)-\nabla f(\lambda^2)|| \leq L_1 ||\lambda^1 -\lambda^2||$ $\forall \lambda^1, \lambda^2 \in \Lambda$, where $||\cdot||$ denotes a norm on $\Lambda$ (which is a \textit{normed linear space}).  We discuss strategies for estimating or \textit{learning} DFO hyperparameters in Section \ref{ss:hyper_learn}.

\subsubsection{STARS}  \label{ss:methods_stars}

The authors in \cite{CW} present the \textit{STep-size Approximation in Randomized Search (STARS)}, a DFO algorithm used to numerically solve the additive and multiplicative OUU problems \eqref{eq:1} and \eqref{eq:2} under  assumptions (i.)-(iv.) STARS uses small perturbations of iterates in $\Lambda$ by the addition of a random vector with components drawn from a normal distribution, computes the noisy function value at the randomly perturbed point, and updates iterates using a Gaussian-smoothed finite-difference for approximate gradient information in a gradient-descent-type scheme. STARS only  requires the ability to evaluate $\hat{f}$ and take random draws from a normal distribution. All in all, the algorithm can be implemented in about 10 lines of code in any standard computing language. (We used Python 3.7 for the results presented in this paper.) 



STARS requires prescribing the values of two key hyperparameters -- denoted $\mu^*_k$ and $h$ -- which are the algorithm's smoothing factor and step size, respectively. ($k=1,\ldots,M$ denotes STARS iterations 1 through $M$, the maximum number of iterations.) The step length $h$ will remain constant for all iterations regardless of the type of noise in $\hat{f}$.
In the case of Additive OUU, $\mu^*_k$ will also be constant, i.e., $\mu^*_k=\mu^*$, a fixed value for all iterations $k=1, \ldots, M$. However in the case of Multiplicative OUU, 
the smoothing factor $\mu^*_k$ will need to be adjusted at every iteration $k$.  
For the additive noise OUU problem \eqref{eq:1}, the values for $\mu^*$ and $h$ are 

\begin{eqnarray} \label{eq:9}
\mu^*:=\left( \frac{8\sigma^2 P}{L_1^2(P+6)^3}\right)^{1/4} \quad  h:=(4L_1(P+4))^{-1},
\end{eqnarray} 

\noindent which are proven as optimal values for STARS' convergence in \cite{CW}.

In the multiplicative noise OUU problem \eqref{eq:2}, the step length $h$ remains the same, exactly as in \eqref{eq:9} above, held constant for each iteration. However the smoothing factor must be updated for each iteration $k=1,\ldots, M.$ As shown in \cite{CW}, the optimal smoothing factor for an iterate $k$ is given by

\begin{eqnarray} \label{eq:10}
\mu^*_k:=\left( \frac{16 \sigma^2 \hat{f}(\lambda^{(k)})^2 P}{L_1^2(1+3\sigma^2)(P+6)^3}\right)^{1/4}.
\end{eqnarray}

We present STARS suited for additive or multiplicative OUU as in \eqref{eq:1} and \eqref{eq:2} in the  pseudocode below.

\vspace{.25cm}

\begin{algorithm}[H] \label{alg:STARS-algo}

\SetAlgoLined

	\KwIn{\texttt{maxit}$=:M$; $\lambda^{(0)}$; $f_0:=\hat{f}(\lambda^{(0)})$; $h$; $k=1$}
 
    \While{$k\leq M$}{
    
  Form smoothing factor $\mu^*_k$
    
  Draw $u^{(k)}$, where $u^{(k)}_p \sim N(0,1)$ for $p=1,\ldots,P$\;
  
  Evaluate $g_k:=\hat{f}(\lambda^{(k-1)}+\mu^*_k\cdot u^{(k)})$\;
  
  Set $ \ds d^{(k)}:=\frac{g_k-f_{k-1}}{\mu^*_k}\cdot u^{(k)}$\;
  
  Set $\lambda^{(k)}=\lambda^{(k-1)}-h\cdot d^{(k)}$\;
  
  Evaluate $f_k:=\hat{f}(\lambda^{(k)})$\; 
  
  Set $k=k+1$\;}

    \KwOut{($\lambda^{(M)}$, $f_M$)}

	\caption{STARS \cite{CW}}

\end{algorithm}

\vspace{.25cm}

STARS typically converges to a minimum when the hyperparameters $\mu_k^*$ and $h$ are within an order of magnitude of their true values. 
The closer the user-defined $\mu_k^*$ and $h$ values are to the truth, the faster STARS converges. 
If $\mu_k^*$ and $h$ are underestimated, STARS will take very small and heavily smoothed steps, converging slowly; however, if the values are overestimated, the algorithm may cause divergence, in the sense that function evaluations will grow with each iterate.

It is then of high interest to tune the values $\mu_k^*$ and $h$ so that function evaluations are not wasted. 
In the proceeding subsection, we discuss learning the values $\mu_k^*$ and $h$ depend on -- namely $\sigma^2$ and $L_1$.

\subsubsection{Learning STARS Hyperparameters}  \label{ss:hyper_learn}
STARS \cite{CW}, like many DFO algorithms \cite{Nesterov}, exhibits optimal convergence if and only if its hyperparameters -- namely the smoothing factor and step size -- are properly tuned. 
Tuning STARS hyperparameters is a matter of learning $\sigma^2$ and $L_1$. 
We examine and propose algorithmic methods of estimating STARS hyperparameters so that one need not specify the hyperparameters at all, fully automating the process of solving  problems \eqref{eq:1} and \eqref{eq:2}.

To estimate $\sigma^2$, we rely on the ECNoise algorithm \cite{MW} which even for $P$ large requires few evaluations of $\hat{f}$ -- often 6 to 10 evaluations will suffice. Briefly, ECNoise uses a set of nearby samples $s^i:=(\lambda^i,\hat{f}(\lambda^i))$ of points $\lambda^i$ along a line in $\Lambda$. Forming a classical difference table of iterative residuals, the authors show that estimators $\hat{\sigma^2}$ to $\sigma^2$ may be formed using scaled averaging and selection criteria discussed more in the next section. Learning $\sigma^2$ is performed prior to STARS, and herein is viewed as a computational expense one must pay up front to ensure convergence of the algorithm.

Learning $L_1$ in our setting is a challenge mainly due to our postulated lack of access to $\nabla f$. 
Given $S$ pointwise estimates to the gradient, which we denote with $\hat{\nabla}f(\lambda^i)$, $i=1,\ldots, S$, and assuming $\lambda^i \neq \lambda^ j$ for $i\neq j$ we could consider an estimator such as

\begin{equation} \label{eq:3}
\reallywidehat{L_1}:= \max_{i \neq j} \frac{\left|\left|\reallywidehat{\nabla}f(\lambda^{i})-\reallywidehat{\nabla}f(\lambda^{j})\right|\right|-2\epsilon^*}{||\lambda^{i}-\lambda^{j}||},
\end{equation} 

\noindent $\epsilon^*=\sup_\xi |\epsilon(\xi)|, i,j=1,\ldots, S,$ given in \cite{Calliess}. There are several problems with such an approach in this setting. Mainly, forming each $\hat{\nabla} f$ is expensive, requiring at least $P+1$ evaluations of $\hat{f}$ for each approximation. Even if one uses only 3 samples $s^i$, $i=1,2,3$, forming $\hat{L_1}$ requires $3(P+1)$ $f$ evaluations, which will often exceed the number of function evaluations needed for STARS to converge, assuming its hyperparameters are reasonably tuned.

Another challenge involves specifying $\epsilon^*$ in \eqref{eq:3}, which is subtracted from the estimator's numerator to control for noise in $\hat{f}$. To save computational expense, we propose forming an initial estimate to $L_1$ by re-using the samples from ECNoise to approximate \textit{directional derivatives} in a finite-difference fashion, avoiding the expensive approach above (as well as approximating the full $\nabla f$). 
Then, once STARS is initiated, $\hat{L_1}$ may be updated using information from approximate directional derivatives formed from iterates (and their corresponding function values). We shall see that the iterates (and intermediate iterates) formed during STARS lend themselves naturally to finite differencing to estimate $L_1$ -- if a larger (more pessimistic) value of $L_1$ is discovered, we replace or \textit{update} $\hat{L_1}$ with this new value.

We also propose fitting a surrogate using $\hat{f}$ values often collected from STARS iterates. 
One may use the closed-form surrogate to estimate $L_1$.
We observe the surrogate-based method for estimating $L_1$ is typically more accurate than finite differences, which -- as we discussed above -- are sensitive to noise.

Updates to $\hat{L_1}$ may be performed between iterations of STARS.
As additional $\hat{f}$ evaluations are obtained, one may update $\hat{L_1}$ anytime an estimate $L_1^\text{update}$ for which $\hat{L_1}<L_1^\text{update}$.  
The updated value $\hat{L_1}=L_1^\text{update}$ may be formed using an un-centered second-order finite difference scheme, which is similar to how $L_1^\text{init}$ is formed, or by refitting a surrogate with the newly-gathered samples of $\hat{f}$. (One must use care to adjust the finite difference formulas to account for iterates and intermediate function evaluations, which are generally un-centered.)

\subsection{Dimension Reduction via the Active Subspace (AS) Method}
We consider functions which map a high-dimensional space $\Lambda$ to a data space $\D$ of smaller dimension; i.e., $\dim \D=D << P=\dim \Lambda$. 
Many functions of interest
actually represent post-processed quantities from the solution of complex physical models. 
It is not often the case that every parameter has an equal impact on function values; usually some parameters matter more than others. If it is possible to mimic the response of $f$ by processing fewer parameters, we can expect computational savings.

We consider AS methods described by Paul Constantine in \cite{Constantine2015} and an equivalent method by T.M. Russi in \cite{Russi}. These techniques seek to explain outputs $f(\Lambda)$ in an AS denoted $\A \subset \Lambda$ for which $\dim (\A) <P$. Here we discuss the theoretical formulation of $\A$. The details of finding $\A$ algorithmically is discussed in the proceeding section.

We note that AS requires, at the very least, approximations to $\nabla f$. For the discussion in this section, we continue with the understanding that $\nabla f$ is approximated in some fashion, the details of which will be discussed in the proceeding Methods section. We assume that $\nabla f(\lambda)$ is square integrable in $\Lambda$ with respect to a probability density $\pi_\Lambda(\lambda)$ that is positive everywhere in $\Lambda$ and 0 otherwise.
 
In AS techniques -- and many other dimension reduction techniques \cite{Russi} -- we transform inputs $\lambda$ to a bounded domain with some fixed variance, typically so that $\lambda\in [-1,1]^P$ for all $\lambda$. Then, as in \cite{ConstantineMC}, we write the \textit{sensitivity matrix}

\begin{equation} \label{eq:4}
W=\int_\Lambda \nabla f(\lambda) \nabla f(\lambda)^\top  \pi_\Lambda(\lambda) d\lambda,
\end{equation} 

\noindent which is a $P\times P$ symmetric positive semi-definite matrix defining a certain covariance of $\nabla f$ in $\Lambda$. This interpretation of \eqref{eq:4} suggests the computation of the eigendecomposition of $W$,

\begin{equation} \label{eq:5}
W=VQV^\top,
\end{equation} 

\noindent where $V$ is $P \times P$ unitary with columns given by the eigenvectors $v^i$, $i=1,\ldots,P$ of $W$ and $Q$ is a diagonal matrix containing the ordered eigenvalues of $W$, $\{q_i\}_{i=1}^P$. To find the AS, we seek a drop-off in magnitude between some pair of eigenvalues, $q_{j}$ and $q_{j+1}$, $1\leq j \leq j+1 \leq P$, where $q_{j}>>q_{j+1}$. 
In this paper, we typically use $95\%$ of the eigenvalues by weight, so that $q_1+\cdots+q_j\geq \tau (q_1+\cdots+q_P)$, where $\tau=0.95$.
We shall see that sometimes the \textit{eigenvalue threshold} $\tau \in (0,1)$ must be changed -- depending on the problem -- to obtain a quality AS.
The \textit{active subspace of $f$} is the span of $v^1,\ldots,v^{j}$, which we denote 
\begin{equation} \label{eq:6}
\A(f)=\text{span}\{v^1,\ldots,v^j\}.
\end{equation}

Likewise, we define the \textit{inactive subspace} of $f$ with

\begin{equation} \label{eq:fix-1}
\I(f)=\text{span}\{v^{j+1},\ldots,v^P\}.
\end{equation}

The fact that $v^1,\ldots,v^{j}$ correspond to large eigenvalues is precisely why they account for the most amount of variance in function values. In fact, one can view an AS as a reasonable choice of principal components after a \textit{full} Principal Component Analysis (PCA) is performed in the gradient space $W$; for more details on this viewpoint, we refer the reader to Section 6.4 in \cite{Russi}.

For a point $\lambda \in \Lambda$, we define

\begin{equation} \label{eq:7}
  \mathcal{P}_\A(\lambda)=\sum_{i=1}^{j}\left( ({v^i})^T \lambda\right)v^i \in \A, 
\end{equation}

\noindent which is a projection of the point $\lambda$ into the AS of $f$. We call this projection an \textit{active variable}, which is a point in the AS $\A$. 
We define a submatrix of $V$ with $V_{\A}:=V_{1:P,1:j}$, the first $j$ columns of $V$ from the eigendecomposition of $W$ in \eqref{eq:2}. 
Then \eqref{eq:7} can be rewritten as $\mathcal{P}_\A(\lambda)=V_\A V_\A^\top \lambda.$
We have arrived at the property that 

\begin{equation} \label{eq:8}
f\left(\mathcal{P}_\A(\lambda)\right) \approx f(\lambda).
\end{equation}

The above property  gives the ability to save computational expense in many scenarios in UQ, including optimization, approximation,
and solving inverse problems \cite{Constantine2015}.
We analogously define a projection into the inactive variables with

\begin{equation} \label{eq:fix-2}
\mathcal{P}_\I(\lambda)=\sum_{i=j+1}^{P}\left( ({v^i})^T \lambda\right)v^i \in \I.
\end{equation}

\noindent We define another submatrix of $V$ with $V_{\I}:=V_{1:P,j+1:P}$, the last $P-j$ columns of $V$ from the eigendecomposition of $W$ in \eqref{eq:5}. 
Then \eqref{eq:fix-2} can be rewritten as $\mathcal{P}_\I(\lambda)=V_\I V_\I^\top \lambda.$

In cases in which $j=1$ or $2$, we can use visualizations to check the extent to which the AS accounts for functions values $\hat{f}(\lambda)$ 
by checking for resolution in a \emph{sufficient summary plot} \cite{Constantine2015}, where one plots active variables against function values. Interpolating these values results with a curve ($j=1$) or surface ($j \geq 2$) forms what is called a \textit{response surface.} In these plots, we hope to see a pattern between the active variables versus their function values. For example, if $f$ is quadratic in its active variables, then we expect to see a quadratic response surface. When $\tilde{j}>3,$ one may check the correlation coefficient between the response surface and true function values since visualization techniques become unfeasible.

In the event that $\nabla f$ is unavailable, the methods above are unusable and approximation methods are required.
We turn our attention to briefly discuss estimating an AS.

\subsubsection{AS Learning} \label{ss:methods_as_learn}

We present two methods we utilized for learning an AS from $\hat{f}$ evaluations. 
In both cases, we often heavily violate a standard assumption that $\Lambda$ samples are independent and random -- the samples we use are deterministic samples from evaluations of $\hat{f}$ from ECNoise (usually 7-10 samples) and STARS (where we often only take enough steps in full variables to train a linear or quadratic surrogate).
\footnote{Recall that points in ECNoise are drawn along a line, but the points in STARS will be at least somewhat random, dictated by a random vector.}
In practice, finding an AS of $f$ without $\nabla f$ will require forming an approximation to $W$ in \eqref{eq:4} in some fashion \cite{ConstantineMC} and necessarily involves full $\nabla f$ approximations.
We present two methods for AS learning we considered -- an approach involving a Monte Carlo approximation to $W$, and an approach involving the use of a closed-form surrogate to obtain $\nabla f$ approximations.

The Monte Carlo approach is simple to implement, found in \cite{Russi}, and equivalent to the method in \cite{ConstantineMC}. For a draw $\lambda^i \in \Lambda$, we obtain an approximation to $\nabla f$ and store the row vector $\nabla f(\lambda^i)^\top$ in a matrix $\tilde{W}$. The eigendecomposition of $\tilde{W}=\tilde{V} Q \tilde{V}^\top$ gives approximations to the eigenvectors and eigenvalues of $W$ as in \eqref{eq:5}.

One initializes the method by performing $S$ random draws of $\lambda^i \in \Lambda$. We then compute $\hat{f}(\lambda^i)$ for all $i=1,\ldots,S$ samples, which we note will require $S$ evaluations of $f$; in a realistic setting, this would require $S$ model solves. We define $D_S:=\{s^i\}_{i=1}^S$, a set of $S$ pairs of samples $\lambda^i$ and their function values. 
Next, we need $\nabla f$ -- which we assume that we do not have in closed form -- evaluated at $\lambda^i$ for all $i=1,\ldots,S$. Hence, we generally need a gradient approximation method \cites{Constantine2015, Smith}.  
Here we form a local linear, global linear, or global quadratic surrogate to $f$ using $D_S$ along the lines of \cite{Constantine2015}. We also consider RBFs as an additional surrogate method.

The gradient of the closed-form surrogate is used to approximate $\nabla f$. \footnote{Despite the fact that we now have formed a potentially global approximation to $\nabla_\Lambda f$ via surrogate, following \cite{CW}, we will still prefer the directional derivative approximations to take steps in STARS since we take steps in exactly the direction we are approximating $\hat{D}_v^i$ for some decent direction $v$.} 
Using  this approximation, we denote each estimation to $\nabla f(\lambda^i)$ with $\hat{\nabla}f(\lambda^i)$ and we define the $P \times S$ matrix $\tilde{W}$ (which is presented below as $\tilde{W}^\top$)

\begin{equation} \label{eq:15}
\tilde{W}^\top:=\begin{bmatrix}
\reallywidehat{\nabla}f(\lambda^1)
\cdot \cdot \cdot
\reallywidehat{\nabla}f(\lambda^S)\\
\end{bmatrix}.
\end{equation}  

Applying an eigendecomposition to $\tilde{W}^\top \tilde{W}$, which is $P \times P$, we obtain $\tilde{W}^\top \tilde{W}=\tilde{V}\tilde{Q}\tilde{V}^\top$ and search for a drop off in the magnitude of the numerical eigenvalues $\{\tilde{q}_i\}_{i=1}^P$ (using an eigenvalue threshold $\tau \in (0,1)$). Assuming such a drop off occurs for an index $\tilde{j}:1\leq\tilde{j}\leq P$, we let

\begin{equation} \label{eq:16}
\tilde{\A}\left(\hat{f}; D_S \right):=\text{span}\{\tilde{v}^1,\ldots,\tilde{v}^{\tilde{j}}\}
\end{equation}

\noindent denote the AS of $\hat{f}$ with respect to the samples $D_S$. (We use this notation including $D_S$ to emphasize the dependence of the approximated AS, $\tilde{\A}$, on the samples taken in $\Lambda$.)

Instead of performing the Monte Carlo method above, one may use a closed-form surrogate to obtain an approximate (but also closed-form) gradient function, which may be used in place of the exact gradient in the formulations of the last section.
We often prefer the surrogate-based approach for its simplicity, as well as its estimation quality and performance (which is often better than the Monte Carlo method).

\begin{center}
  $\ast$~$\ast$~$\ast$
\end{center}

In the proceeding Methods section (\ref{s:Methods}), we present the \textit{Active STARS (ASTARS)} and \textit{Fully-Automated ASTARS (FAASTARS)} algorithms.
ASTARS leverage AS dimension reduction in $\Lambda$ to perform STARS steps in a lower-dimensional space when possible.
FAASTARS is fully-automated in the sense that the user need not specify hyperparameters nor provide the true AS of $f$.
In the Results section (\ref{s:Results}), we analyze and compare the performance of STARS, ASTARS ,and FAASTARS in a series of examples. 
Finally, in the Conclusion and Discussion section (\ref{s:Conc}), we review the extent of ASTARS and FAASTARS efficiency. Limitations are discussed and future research questions are posed. 

\section{Methods}  \label{s:Methods}

\subsection{Active STARS (ASTARS)} \label{ss:methods_astars}
Given $\hat{f}$ and the exact AS $\A$ of stochastic-free signal $f$, we are interested in investigating the effectiveness of optimizing $\hat{f}$ in its active variables. There are several approaches one may consider, and some of those approaches and their corresponding results are discussed in the remainder of the paper; in this section, we focus on performing STARS in the true, known AS of the true signal, $f$, denoted $\A$.

Active STARS, or \textit{ASTARS}, is a modification of STARS in which iterates only take random walks in directions lying in $\A$. In detail, at iteration $k$, STARS uses random walks in directions given by drawing a random vector $u^{(k)}$ of dimension $P$ in which every component $u_i^{(k)},i=1,\ldots,P$ of $u^{(k)}$ is  drawn from a specified normal distribution. 
Instead, given the first $j$ eigenvectors $v^1,\ldots, v^j$ spanning $\A$, ASTARS takes $j$ draws from a specified normal distribution, which we denote $r_i\sim N(0,\omega_i^2)$, $i=1,\ldots,j$, defining the random vector $u_\A^{(k)}$ for the $k$-th random \textit{active direction} as 

\begin{eqnarray}\label{eq:17}
u_\A^{(k)}=\sum_{i=1}^j r_i v^i, \quad  r_i\sim N(0,\omega_i^2), \quad i=1,\ldots,j.
\end{eqnarray}

The direction $u_\A^{(k)}$ is a randomly-weighted linear combination of the active directions of $f$ and is the direction used in place of $u^{(k)}$ in STARS. In the case that there is not a large drop-off in the spectrum of $\tilde{W}$, then all $P$ directions are active, and ASTARS reduces to performing STARS (in all variables).

In ASTARS (Algorithm 2 below), we equally weight the $j$ active directions using unit variance, in the sense that $\omega_i^2=1$ for $i=1,\ldots,j$. The use of unit variance in the random coefficients matches the theoretical assumptions in \cite{CW}, and is an assumption in the proofs in our technical report.
However, one may suspect other choices of $\omega_i$ may improve ASTARS performance in some cases.

Other weighting schemes considered include taking $\omega_i=\sqrt{\tilde{q_1}}/\sqrt{\tilde{q_i}}$, $i=1,\ldots,j$, where we recall $\tilde{q_i}$ denotes the $i$-th numerical eigenvalue obtained from the eigendecomposition of $\tilde{W}^\top \tilde{W}$ which are indexed so that $\tilde{q_1}\geq \cdots \geq \tilde{q_j}$.
In our numerical experiments, this and other alternate weighting schemes exhibited promising ASTARS convergence, and further research is needed, which we plan to present in a follow-up paper investigating this and other extensions of the ASTARS framework.

ASTARS requires modifying the initialization and changing the second step of STARS (Algorithm 1) by replacing $u^{(k)}$ with $u^{(k)}_\A$ as we discussed above. 
ASTARS also uses modified STARS hyperparameters. For the additive noise OUU problem, we define the \textit{active hyperparameters} $\mu_\A^*$ and $h_\A$

\begin{eqnarray} \label{eq:18}
\mu_\A^*:=\left( \frac{8\sigma^2 j}{L_1^2(j+6)^3}\right)^{1/4} \, h_\A:=(4L_1(j+4))^{-1},
\end{eqnarray} 

\noindent the active smoothing factor and step length, respectively. For the multiplicative noise OUU problem, the step length $h_\A$ remains the same, exactly as in \eqref{eq:18} above but the optimal \textit{active smoothing factor} for the $k$-th iterate of ASTARS, $k=1,\ldots, M$, is given by

\begin{eqnarray} \label{eq:19}
(\mu_\A^*)_k:=\left( \frac{16 \sigma^2 \hat{f}(\lambda^{(k)})^2 j}{L_1^2(1+3\sigma^2)(j+6)^3}\right)^{1/4}.
\end{eqnarray} 

We may use $\mu_\A^*$ to generally denote and discuss the active smoothing factor hereon with the understanding that in the case of multiplicative noise one actually uses $(\mu_\A^*)_k$.

In Algorithm 2 below, we present ASTARS.

\begin{algorithm}[h] 

\SetAlgoLined

	\KwIn{\texttt{maxit}$=:M$; $\lambda^{(0)}$; $f_0:=\hat{f}(\lambda^{(0)})$; $h_\A$; $\tilde{V}_\A:=\tilde{V}_{1:P,1:j}$; $k=1$}
 
    \While{$k\leq M$}{
    
  Form smoothing factor $(\mu^*_\A)_k$\;
    
  Draw $r^{(k)}$, where $r^{(k)}_p \sim N(0,1)$ for $p=1,\ldots,j$ and set $u^{(k)}_\A:= \tilde{V}_\A r^{(k)}$\;
  
  Evaluate $g_k:=\hat{f}(\lambda^{(k-1)}+ (\mu^*_\A)_k u_\A^{(k)})$\;
  
  Set $ d^{(k)}:=\frac{g_k-f_{k-1}}{(\mu^*_\A)_k}u_\A^{(k)}$\;
  
  Set $\lambda^{(k)}=\lambda^{(k-1)}-h_\A\cdot d^{(k)}$\;
  
  Evaluate $f_k:=\hat{f}(\lambda^{(k)})$\; 
  
  Set $k=k+1$\;}

    \KwOut{($\lambda^{(M)}$, $f_M$)}

	\caption{ASTARS}
	
	\label{alg:ASTARS-algo}

\end{algorithm}


\noindent \textbf{ASTARS Corollary.}
\textit{Let the vectors $u_\A^{(k)}$ denote those drawn using Algorithm \ref{alg:ASTARS-algo} (zero mean, unit variance in each component); let $f\in \mathcal{C}^{1,1}(\Lambda)$ and assume $f$ is convex; and assume that the i.i.d. noise draws $\epsilon(\xi)$ are additive, zero mean, with bounded variance $\sigma^2$ for all $\xi$.
Fixing the step size $h_\A$ in \eqref{eq:18}, the active smoothing parameter $\mu_\A^*$ in \eqref{eq:18} minimizes the error between the gradient oracle in Algorithm \ref{alg:ASTARS-algo},} given by
$$\frac{\hat{f}(\lambda^{(k-1)}+\mu^*_\A u_\A^{(k)}) - \hat{f}(\lambda^{(k-1)})}{\mu^*_\A} u_\A^{(k)},$$ 
\noindent \textit{and the true directional derivative of $f$ in the direction $u_\A^{(k)}$ in the $j$-dimensional $\A$.
}

\noindent \textbf{Remark.} \textit{The preceding corollary implies that using the fixed step length $h_\A$, ASTARS uses an optimal choice of smoothing parameter $\mu^*_\A$, in the sense that $\mu^*_\A$ minimizes the error between our approximate directional derivative formed  in Algorithm \ref{alg:ASTARS-algo}.
Since ASTARS takes steps in the $j$-dimensional space $\A$, the hyperparameters from STARS must be redefined to remain optimal. The hyperparameters in \eqref{eq:9} and \eqref{eq:010} are proven to be optimal hyperparameters for the convergence of STARS in the $P$-dimensional space $\Lambda$ by the authors in Theorem 4.3 of \cite{CW}. 
Now, replacing $P=\dim \Lambda$ in \eqref{eq:9} and \eqref{eq:10} with $j=\dim \A$, we obtain \eqref{eq:18} and \eqref{eq:19}.
We will present a proof of the corollary above as ASTARS Corollary 2 in our technical report.
Note also that the authors \cite{CW} have an analogous result for the case of multiplicative noise, and so we note that we, too, have an optimal active smoothing constant given in \eqref{eq:19} in the case of multiplicative noise.}

In our technical report, we prove that the complexity of ASTARS is

\begin{eqnarray} \label{eq:complex-summary}
M \sim \mathcal{O}\left( \frac{L_1 j R^2_\A}{\epsilon_{\text{tol}}} \right).
\end{eqnarray}

\noindent Here, $R^2_\A$ is a bound on the squared norm of the difference between any initial iterate and the true minimum of $f$, both projected in the inactive subspace $\I$.
$\epsilon_{\text{tol}}>0$ is a final accuracy which is bounded below by terms involving the variance in the noise, as well as by terms involving the inactive subspace of $f$; for details, refer to the technical report.

\subsection{Fully-Automated ASTARS (FAASTARS)} \label{ss:methods_faastars}
We now introduce a fully-automated version of ASTARS (Algorithm 2), in the sense that the user need not specify anything beyond an initial iterate $\lambda^{(0)}$, its evaluation $\hat{f}(\lambda^{(0)})$, the black-box objective function $\hat{f}$, and a maximum number of iterations, $M$. We call this algorithm \textit{Fully-Automated ASTARS (FAASTARS)}.

We first note that $\hat{f}$ need not be in closed form; we only require that $\hat{f}$ is a callable function, which we recall may actually represent a post-processed quantity from, for instance, a noisy solution of a PDE evaluated at some point in parameter (e.g., phase) space. 
FAASTARS estimates $\sigma^2$ and $L_1$ from a handful of upfront samples taken from performing ECNoise (and recycled for $L_1$ learning). As well, $\A$ is estimated from regular STARS iterates (supplemented with the original ECNoise samples as well) during a \textit{STARS burn-in} phase.

In the following, we outline FAASTARS by breaking it down into three phases: (1.) a hyperparameter learning phase; (2.) a STARS burn-in phase; (3.) an approximate ASTARS phase.

In the first phase, the estimator $\hat{\sigma}^2$ will be formed immediately by using ECNoise on 7 to 10 sampled points in $\Lambda$. As discussed above, the samples $s^i$ created by ECNoise lend themselves to forming $L_1^\text{init}$ in a finite-difference scheme approximating directional derivatives and $\hat{L_1}=L_1^\text{init}$ is used.

We can also form a surrogate $F$ from the ECNoise points to estimate $L_1$.
By obtaining the value of the closed-form hessian of the surrogate $F$, $\nabla^2 F$, we obtain a lower bound $L_1$, as we will see. The surrogate can be improved by incorporating STARS iterates into the set of samples used to form the surrogate, which we will present in the next section, about adaptive sampling.

We will first need approximated hyperparameters, given by

\begin{equation} \label{eq:hattedhypers}
\hat{\mu}^*:=\left( \frac{8\hat{\sigma}^2 P}{\hat{L_1}^2(P+6)^3}\right)^{1/4} \quad \quad \hat{h}=(4\hat{L_1}(P+4))^{-1}.
\end{equation}

Note that the approximated active smoothing factor will be optimal in $\Lambda$ -- given the available information -- in a somewhat analogous fashion as our result in ASTARS Corollary 2 (see technical report), but with a specified loss of optimality as the estimates $\hat{L_1}$ and $\hat{\sigma}$ stray from their true values.
(See STARS Theorem 4.3 (Modified) in \ref{ss:stars-conv}.)
In particular, we find that STARS will either diverge or make no progress when the values are underestimated or overestimated, respectively.
That is, underestimation of either value may lead to STARS' divergence, but for distinct reasons in each case.

When the variance in the noise $\sigma^2$ is underestimated, $\mu^*$ is also underestimated, and thus we may not take a large enough step to successfully perturb $\hat{f}$ enough to see a change in function value larger than the noise level itself, leading to inaccurate derivative information and potentially descent steps of poor quality.
(Note $\sigma^2$ does not appear in the step size.)
When the gradient Lipschitz constant $L_1$ is underestimated, both $\mu^*$ and $h$ will be too large, meaning we may take too large of a step for the finite difference approximation to be accurate \textit{and} we may take too large of a descent step (in a bad direction), causing a quick rise in the function valeus we see.
Indeed, underestimating $L_1$ is to be avoided at all costs.

We may also form a linear surrogate $F$ from ECNoise samples (if we have enough data) to initiate our approximation to $L_1$ by computing the matrix norm of the closed form Hessian of $F$, $||\nabla^2 F||$, for which we have $||\nabla^2 F|| \approx ||\nabla^2 f||  \leq L_1 $. 
The first phase of FAASTARS is given in the algorithm below. 

\vspace{.25cm}

\begin{algorithm}[H]  \label{alg:FAASTARS-algo-1}

\SetAlgoLined

	\KwIn{$\lambda^{(0)}$; $k=0$}
	
	\While{$k=0$}{
	
	Run ECNoise using $\lambda^{(0)}$ as base point and obtain $\hat{\sigma}^2$\;
	
	Initialize storage array $D_S$ for samples of $f$ and store $\{s^i\}$ formed by ECNoise\;
	
	Use $\{s^i\}$ to form second-order FD approximation (or form linear surrogate $F$ if $S>P+1$ to compute $||\nabla^2 F||$) for $L_1^\text{init}$\;

	}
	
	\KwOut{$\hat{\sigma}^2$; $L_1^\text{init}$; $D_S$}
	
	\caption{FAASTARS for Additive or Multiplicative OUU, Phase 1: Hyperparameter Learning}
	
\end{algorithm}

\vspace{.25cm}

Next, in the second phase, we perform standard STARS (Algorithm \ref{alg:STARS-algo}) until enough samples are obtained to perform AS analysis via a surrogate to form needed $\nabla f$ evaluations. We let $M_\A$ denote the number of iterations needed to find the AS from samples, as we see in the first phase of FAASTARS above. We note that $M_\A$ will depend on the type of surrogate formed (e.g., using local linear versus global quadratic versus RBFs). FAASTARS will not begin its ASTARS routine until enough samples have been taken to form a surrogate, based on the chosen surrogate method -- RBFs are the default surrogate method when none is provided using known formulas. For example, if one wishes to use a globally quadratic surrogate, $(P+1)(P+2)/2$ samples of $f$ are required \cite{Smith}. 
Recalling that every STARS step requires \textit{two} evaluations of $f$, $M_\A=(P+1)(P+2)/4$ will suffice for quadratic surrogates.
$M_\A$ is not a value the user has to provide; FAASTARS will automatically terminate its STARS burn-in period as soon as $k$ is large enough for the chosen or default surrogate to be formed.

After all steps of approximate STARS in Phase 2 are taken, $\tilde{\A}$ is found using the Monte Carlo method described in Chapter 1 using the samples/iterates $\{s^i\}$ collected from both ECNoise and $M_\A$ steps of standard STARS. 
We present the second phase of FAASTARS in the pseudocode below.

\vspace{.25cm}

\begin{algorithm}[H]  \label{alg:FAASTARS-algo-2}

	\KwIn{$\lambda^{(0)}$; $f_0:=\hat{f}(\lambda^{(0)})$; $k=1$; Surrogate method (optional, default is RBFs); $\hat{\sigma}^2$; $L_1^\text{init}$; $D_S$; eigenvalue threshold $0<\tau<1$ (optional, default is $\tau=0.95$)}

	Determine $M_\A$ based on chosen surrogate method or RBFs if none is provided\;
	
	Form step length $\hat{h}$ using $\hat{\sigma}^2$ and $L_1^\text{init}$
	
	\While{$ 1\leq k \leq M_\A$}{
    
  Form smoothing factor $\hat{\mu}^*_k$ using $\hat{\sigma}^2$ and $L_1^\text{init}$\;

  Find $\lambda^{(k)}$ and evaluate $f_k:=\hat{f}(\lambda^{(k)})$ via STARS (Algorithm \ref{alg:STARS-algo})\; 
  
    Store $(\lambda^{(k-1)}+\hat{\mu}^*_k u^{(k)},g_k)$ and $(\lambda^{(k)},f_k)$ as samples $\{s^i\}$ in $D_S$\; 
  
  \textbf{Optional:} Form $(L_1^\text{update})_k$ via FD with $D_S$ (or use $||\nabla^2 F||$ for surrogate $F$ with $D_S$)\;
  
  \textbf{Optional:} If $(L_1^\text{update})_k>L_1^\text{init}$, set $L_1^\text{init}=(L_1^\text{update})_k$ and re-compute $h$\;

  Set $k=k+1$\;}
	
	Use $D_S$ to form surrogate via selected method\; 
	
	Form $\tilde{W}$, and apply SVD to obtain $\tilde{W}=\tilde{V}\tilde{Q}\tilde{V}^\top$\;
	
	Find a drop-off index $\tilde{j}$ for which $\tilde{q}_1+\cdots+\tilde{q}_{\tilde{j}}\geq \tau (\tilde{q}_1+\cdots+\tilde{q}_P)$
  
    Set $\tilde{V}_{\tilde{\A}}:=\tilde{V}_{1:P,1:\tilde{j}}$ and form $h_\A$ using $\tilde{j}$ for $\dim \tilde{\A}$; $\hat{\sigma}^2$ and $L_1^\text{init}$\;
    
    \KwOut{$\lambda^{(M_\A)}$, $f_{M_\A}:=\hat{f}(\lambda^{(M_\A)})$, $\tilde{V}_{\tilde{\A}}$}
    
    	\caption{FAASTARS for Additive or Multiplicative OUU, Phase 2: Approximate STARS burn-in}

\end{algorithm}

\vspace{.25cm}

Note that at the end of each standard STARS iteration in the burn-in phase, we have the functionality to form candidates $L_1^\text{update}$ for $\hat{L_1}$ via finite difference (FD) approximations to the directional derivatives in the direction of the corresponding descent step.
We can also use a surrogate as in Phase 1.
Regardless, we reject the candidate update anytime $L_1^\text{update} \leq L_1^\text{init}$, since we are always searching for the most pessimistic bound to $L_1$ to avoid divergence.

With $\tilde{\A}$ in hand, we must first update the hyperparameters so that they are computed with the value $\tilde{j}=\dim \tilde{\A}$ (and not $j = \dim \A$, since it is generally unknown in this setting).
We define the approximated active hyperparamters,

\begin{equation} \label{eq:hattedactives}
\hat{\mu}^*_{\tilde{\A}}:=\left( \frac{8\hat{\sigma}^2 \tilde{j}}{\hat{L_1}^2(\tilde{j}+6)^3}\right)^{1/4} \quad \quad \hat{h}_{\tilde{\A}}:=(4\hat{L_1}(\tilde{j}+4))^{-1}.
\end{equation}

In the third phase, we pick up where standard STARS left off, and we perform ASTARS in the approximated AS for the remaining iterations until the maximum number of iterations, $M$, is met. 
We have the functionality to also update $\hat{L_1}$ in a similar fashion as the burn-in phase during the ASTARS phase.
That is, we may continue to use finite differences, or use the surrogate $F$ that we form for AS approximations to update our approximation to $L_1$.
(In practice, this surrogate could be formed from the initial ECNoise points and also be recalculated at every step as we gain more samples of $\hat{f}$.)
We may take $L_1^\text{update} = ||\nabla^2 F||$ and similarly to above reject the candidate update anytime $L_1^\text{update} \leq L_1^\text{init}$; otherwise, we have a new initial $L_1$, where $L_1^{\text{init}}=||\nabla^2 F||$.

 \vspace{.25cm}

\begin{algorithm}[H]   \label{alg:FAASTARS-algo-3}

	\KwIn{\texttt{maxit}$=:M$; $\lambda^{(M_\A)}$; $f_{M_\A}:=\hat{f}(\lambda^{(M_\A)})$; $k=M_\A+1$; $\hat{\sigma}^2$; $L_1^\text{init}$; $\tilde{V}_{\tilde{\A}}$}
 
    \While{$M_\A < k\leq M$}{
    
  Form step size $\hat{h}_{\tilde{\A}}$ and smoothing factor $\hat{\mu}^*_{\tilde{\A}}$ using $\tilde{j}$ for $\dim \tilde{\A}$; $\hat{\sigma}^2$ and $L_1^\text{init}$\;
    
  Draw $r^{(k)}$, where $r^{(k)}_p \sim N(0,1)$ for $p=1,\ldots,\tilde{j}$ and set $u^{(k)}_{\tilde{\A}}:= \tilde{V}_{\tilde{\A}} r^{(k)}$\;
  
  Evaluate $g_k:=\hat{f}(\lambda^{(k-1)}+ \hat{\mu}^*_{\tilde{\A}}\cdot  u^{(k)}_{\tilde{\A}})$\;
  
  Set $ d^{(k)}:=\frac{g_k-f_{k-1}}{\hat{\mu}^*_{\tilde{\A}}} u^{(k)}_{\tilde{\A}}$\;
  
  Set $\lambda^{(k)}=\lambda^{(k-1)}-\hat{h}_{\tilde{\A}}\cdot d^{(k)}$\;
  
  Evaluate $f_k:=\hat{f}(\lambda^{(k)})$\;
  
    \textbf{Optional:} Update $L_1^{\text{init}}$ with FD or surrogates as in Phase 2 (requires updating $D_S$)\;
      
  Set $k=k+1$\;}

    \KwOut{($\lambda^{(M)}$, $f_M$)}

    	\caption{FAASTARS for Additive or Multiplicative OUU, Phase 3: Approximate ASTARS}

\end{algorithm}

\vspace{.25cm}

We now have all the necessary components for performing FAASTARS. In summary, FAASTARS has three major phases: (1.) an initial, relatively inexpensive learning phase where we acquire the estimates $\hat{\sigma}^2$ and $L_1^\text{init}$; (2.) a STARS burn-in phase (in full variables) where we acquire enough samples to compute an AS using the Monte Carlo methods above; and (3.) an ASTARS phase, where we use the learned AS, $\tilde{\A}$. Note that in both of the latter phases, we will update $\hat{L_1}$, if and only if we obtain a more pessimistic (larger) estimate.

In our technical report, we show that the approximately active hyperparameters \eqref{eq:hattedactives} are chosen to minimize the error in finite difference approximations to directional derivatives. 
We also prove that the complexity of FAASTARS is similar to that of ASTARS, but with a specified inflation of certain bounds, which arise from approximating both hyperparameters and $\A$.

\subsubsection{AS Retraining}

We introduce a method we call \textit{AS retraining} which can be optionally applied in place of Phase 3 FAASTARS.
AS retraining begins with an initial approximation $\tilde{\A}$ for $\A$, which is obtained once the minimum number of samples $M_\A$ are gathered during FAASTARS Phase 2.
We then take $k_T \geq 1$ approximate ASTARS steps using $\tilde{\A}$ (as in FAASTARS Phase 3) and recompute $\tilde{\A}$ with our new set of samples of $f$.
We continue in this fashion until the maximum iteration count $M$ is reached. 

Samples obtained during FAASTARS Phase 3 contain information which is more local to a given Phase 3 iterate.
Hence, adaptive sampling incorporates local information to approximate $\A$, more aligned with how $f$ changes in the region of recent samples.
It is almost always better to include more samples whenever we can to improve $\tilde{A}$, borne out in numerical results. 
Indeed, we use AS retraining to produce all FAASTARS results in this work.
We present the adaptive sampling algorithm below.

\vspace{.25cm}

\begin{algorithm}[H]

\SetAlgoLined

	\KwIn{\texttt{maxit}$=:M$; $\lambda^{(M_\A)}$; $f_{M_\A}:=\hat{f}(\lambda^{(M_\A)})$; re-training phase length $k_T$; $D_S$; $\hat{\sigma}^2$; $L_1^\text{init}$; $\tilde{V}_{\tilde{\A}}$ from FAASTARS Phase 2; $l = 1$; $k=M_\A+1$;}

    \While{$k\leq M_\A+lk_T$}{
    
  	Form approximate active step size $\hat{h}_{\tilde{\A}}$ and approximate active smoothing factor $(\hat{\mu}^*_{\tilde{\A}})_k$\;
    
  	Take $k_T$ steps of (approximate) ASTARS with $\tilde{\A}$ (as in FAASTARS Phase 3) and store all samples of $f$ in $D_S$\;
  
  	Retrain $\tilde{\A}$ and recompute $\tilde{V}_{\tilde{\A}}$ with $D_S$ (containing $k_T$ new samples)\; 
  
  	Set $l=l+1$ and exit loop only if $k_S+lk_T \geq  M$\;}

    \KwOut{($\lambda^{(M)}$, $f_M$)}
    
    \caption{AS Retraining}

\end{algorithm}

\vspace{.25cm}

At times, a challenge with adaptive sampling is poor quality in the samples obtained from the FAASTARS Phase 3 steps, after the original formation of $\tilde{A}$.
When the step size is small and $f$ changes very little -- sometimes due to inaccurate descent directions -- samples are generally uninformative about the behavior of $f$.
This challenge is also present, more generally, in all AS approximations involving samples of $\hat{f}$ taken in a partially deterministic DFO algorithm. 
Again, these samples are often clustered together in $\Lambda$ due to the small steps we take both in finite differencing and in a descent step.
Some of these challenges are addressed in a forthcoming extensions paper.
Regardless, we find incorporating as many samples as possible for computations involving $\A$ usually improves numerical results.

\section{Results} \label{s:Results}

In this section, we present results of using STARS, ASTARS, and FAASTARS in several examples with noisy objective functions with high-dimensional parameter spaces. 

\noindent \textbf{Example 1. (Toy)} \textit{Let $\hat{f}: \Lambda=\R^{P} \to \R =\D$. Fixing a weight vector, $w \in \Lambda$, where $w \neq 0$, we define}

\begin{eqnarray} \label{eq:103}
\hat{f}(\lambda; \xi):=\left(w^\top \lambda\right)^2 + \epsilon(\xi),  
\end{eqnarray}

\noindent \textit{$\epsilon(\cdot) \sim N(0,\sigma^2)$, where $\sigma^2=1 \times 10^{-4}.$
We note the noise-free signal $f$ is convex. 
The minimum of $\hat{f}$ is given by $0 \in \Lambda$ with minimizer $f^*=0$. 
Here, the noise-free signal $f$ has a one-dimensional AS (i.e., $j=1$) in the direction $w$. 
Also note $L_1=2$. 
We considered $P=20$ and note $D=1$. 
We took $w_i=1$ in all directions and an initial iterate $\lambda^{(0)}$ with components drawn from a zero-mean Gaussian with unit variance, scaled by 10. 
First, we performed 1000 trials each of STARS using exact hyperparameters, ASTARS using exact active hyperparemters as well as the exact AS, and
using exact hyperparameters but learned AS, we performed FAASTARS.
A maximum iteration count of $2P^2=800$ was used.
AS's were re-trained every $2P=40$ steps for FAASTARS using an eigenvalue threshold of 99 percent.
No noise regularization was required.
We then performed 500 trials each of STARS and FAASTARS using estimated hyperparameters and a maximum iteration count of $500$.
In this case, AS's were relearned every $P=20$ steps for FAASTARS using an eigenvalue threshold $\tau = 0.95$ and a noise regularization at the level of $\sigma^2$.
}

\begin{figure}[h]
	\centering
	\includegraphics[scale=.155]{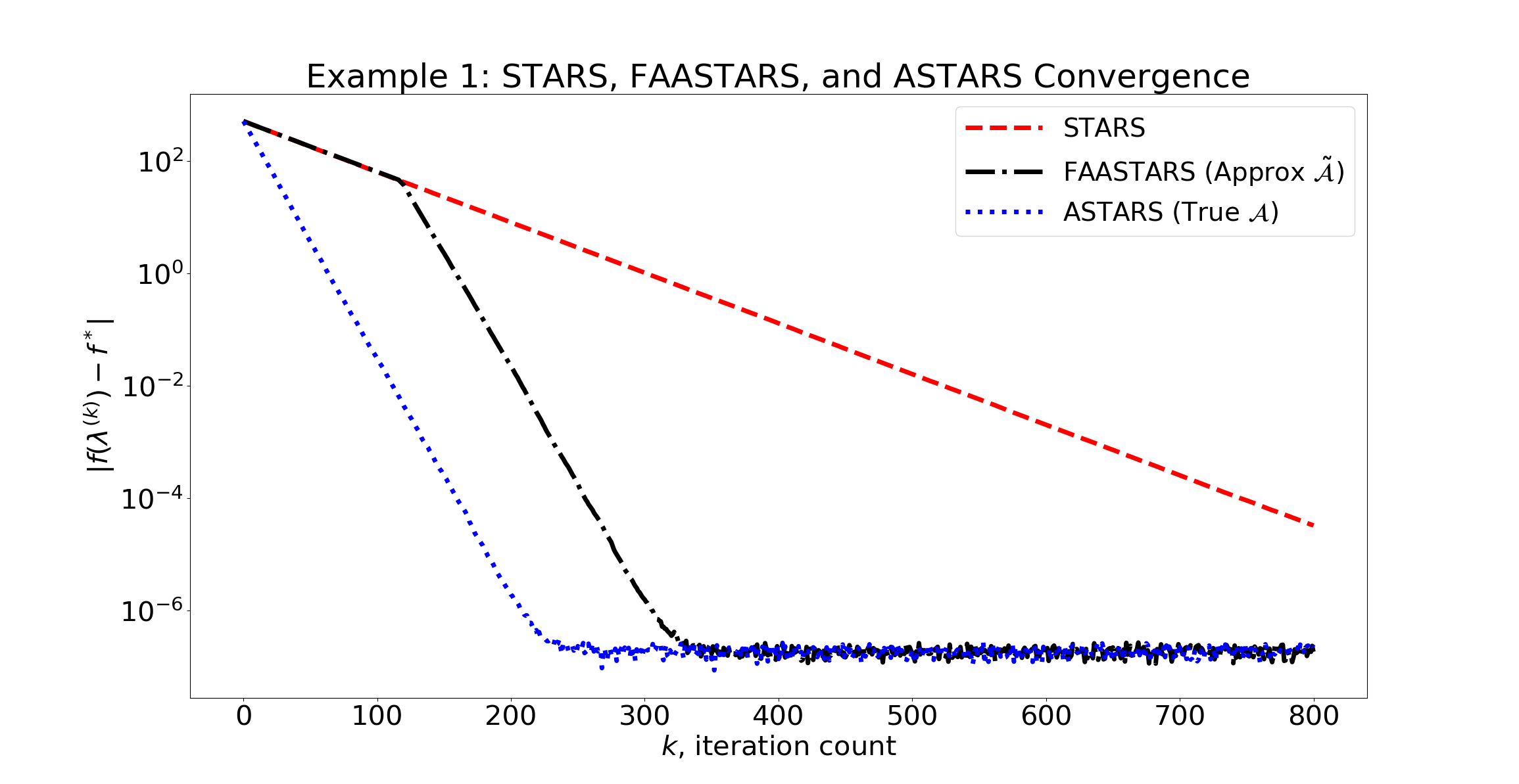} 
    \includegraphics[scale=.155]{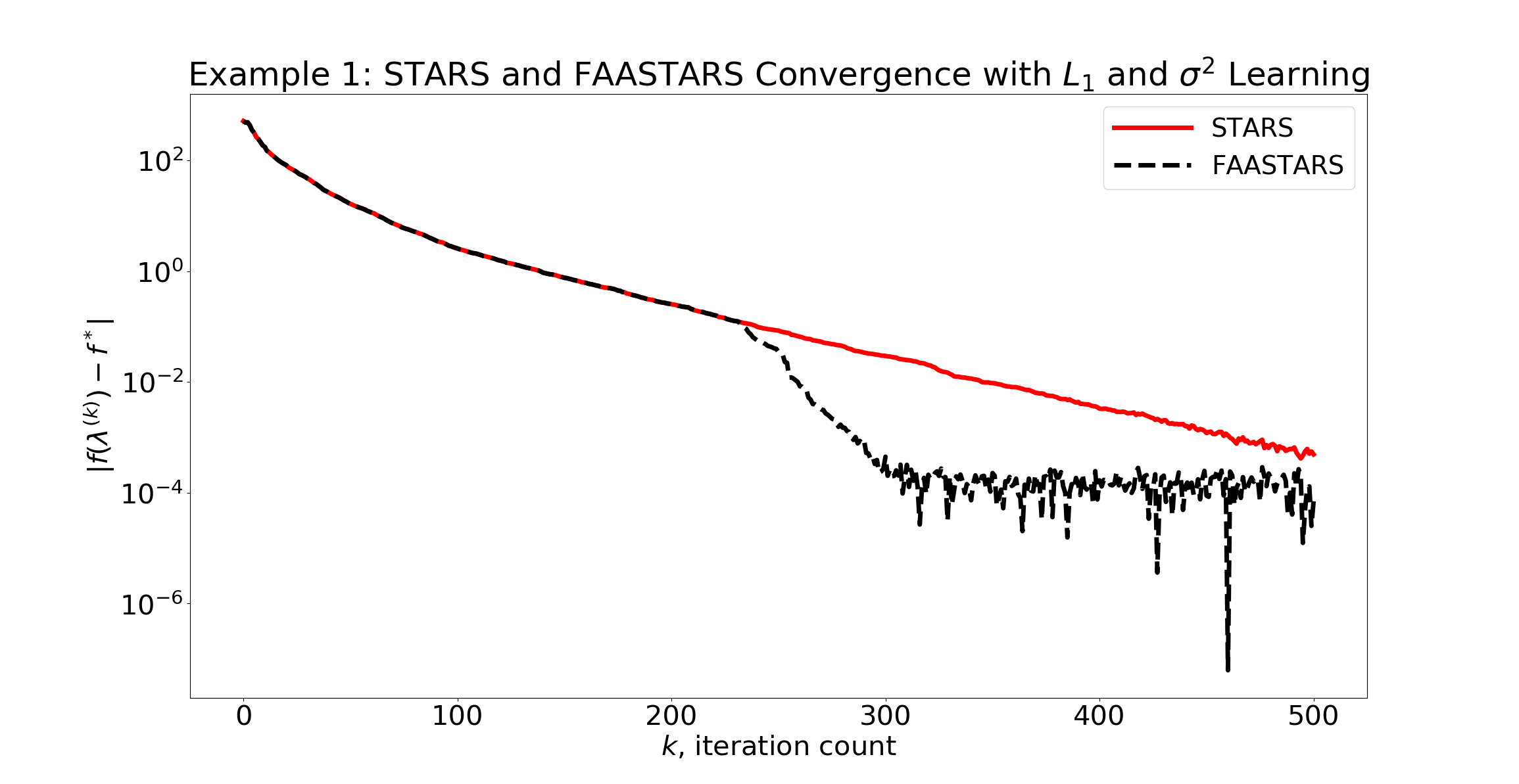}
	\caption{We show the convergence of STARS, FAASTARS, and ASTARS with and without hyperparameter learning for Example 1.}
	\label{fig:ex1}	
\end{figure}

\noindent \textbf{Example 2. (Active sphere)} \textit{Let $\hat{f}: \Lambda=\R^{P} \to \R =\D$,}

\begin{eqnarray} \label{eq:104}
\hat{f}(\lambda; \xi):=\sum_{i=1}^j \lambda_i^2 + \epsilon(\xi),  
\end{eqnarray}

\noindent \textit{$\epsilon(\cdot) \sim N(0,\sigma^2)$, where $\sigma^2=1 \times 10^{-3}.$
We note the noise-free signal $f$ is convex. 
The minimum of $\hat{f}$ is given by $0 \in \A$ with arbitrary components in $\I$ with minimizer $f^*=0$. 
Here, $\hat{f}$ has a $j$-dimensional AS spanned by the first $j$ standard basis vectors in $\Lambda$. 
Also note $L_1=2$. 
We considered $P=20$, $j=10$, and note $D=1$. 
We took initial iterate $\lambda^{(0)}$ with components drawn from a zero-mean Gaussian with unit variance, scaled by a factor of 10. 
First, we performed 100 trails each of STARS using exact hyperparameters, ASTARS using exact active hyperparemters as well as the exact AS, and
using exact hyperparameters but learned AS, we performed FAASTARS.
A maximum iteration count of $2P^2=800$ was used.
AS's were re-trained every $P=20$ steps using an eigenvalue threshold of 99.9 percent.
Noise regularization of $\sigma^2$ improved results.
We then performed 250 trials of STARS and FAASTARS using estimated hyperparameters and a maximum iteration count of $4P^2=1600$.
In this case, AS's were relearned every $P=20$ steps for FAASTARS using an eigenvalue threshold $\tau = 0.9995$ and a noise regularization at the level of $\sigma^2$.}

\begin{figure}[h]
	\centering
	\includegraphics[scale=.155]{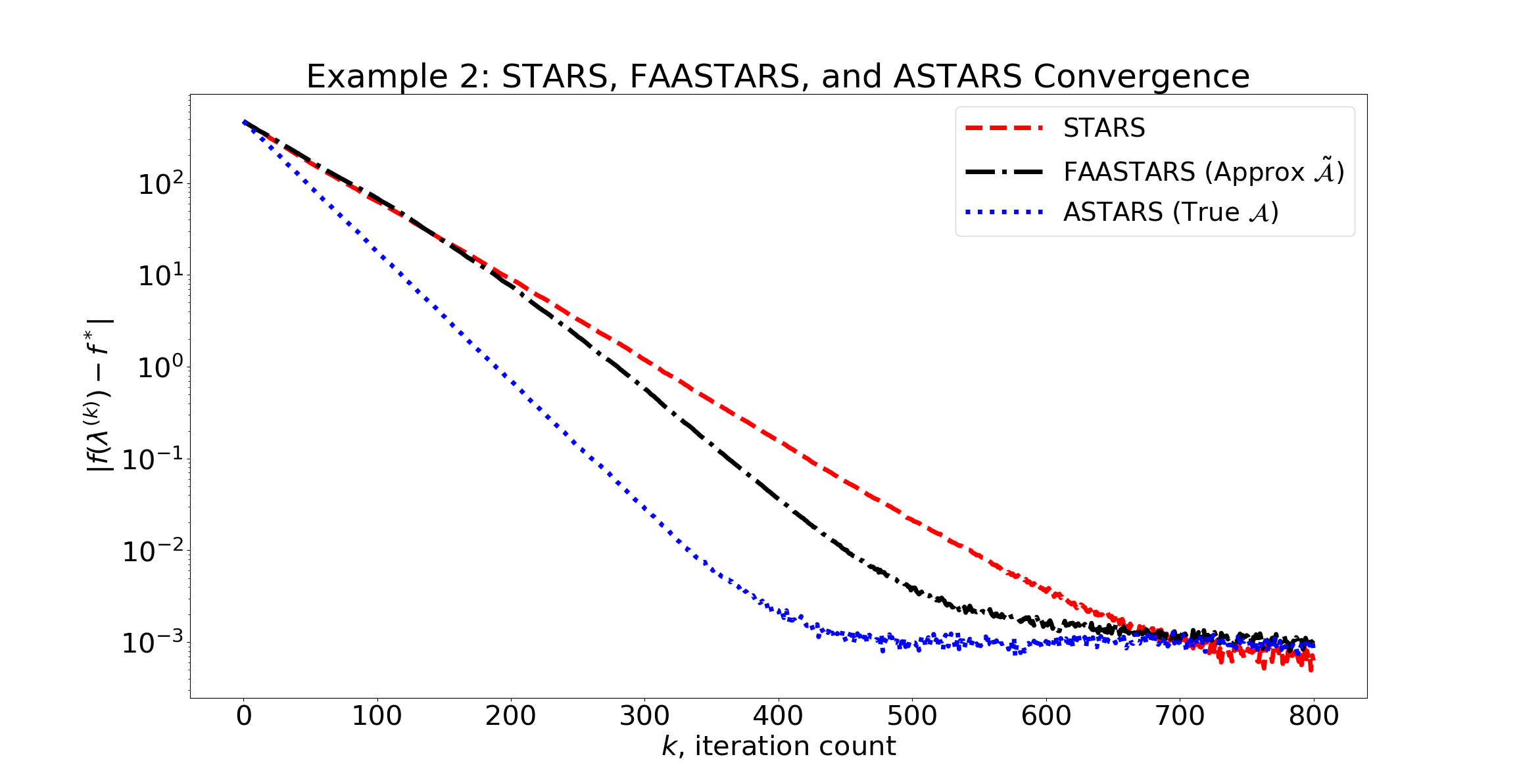}
	\includegraphics[scale=.155]{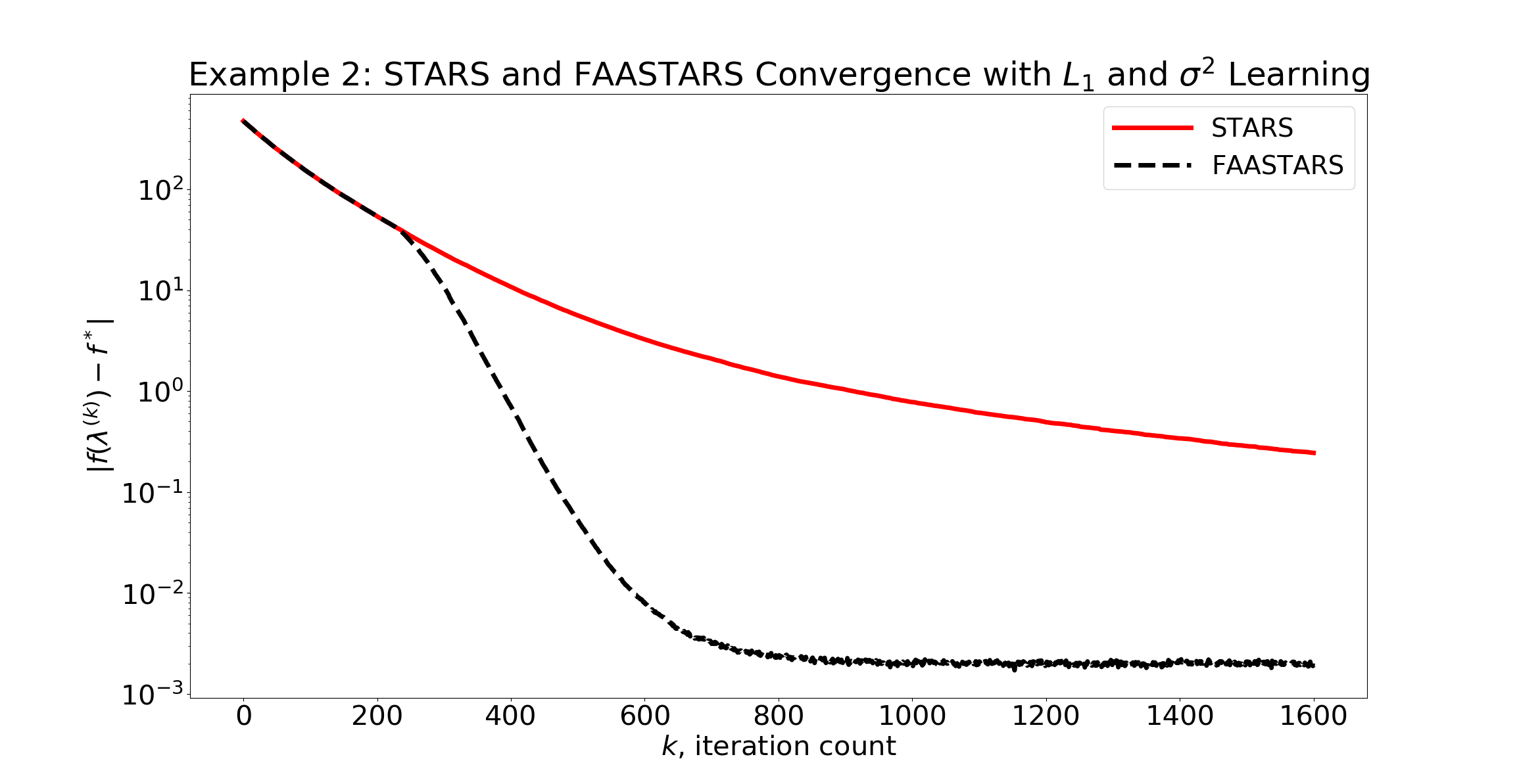}
	\caption{We show the convergence of STARS, FAASTARS, and ASTARS with and without hyperparameter learning for Example 2.}
	\label{fig:ex2}	
\end{figure}

\newpage

\noindent \textbf{Example 3. (Nesterov's function -- active version)} \textit{Let} $\hat{f}: \Lambda=\R^P \to \R =\D,$

\begin{eqnarray} \label{eq:105}
\hat{f}(\lambda; \xi)=\frac{1}{2}\left(\lambda_1^2 + \sum_{i=1}^{j-1} (\lambda_i-\lambda_{i+1})^2+\lambda_j^2\right)-\lambda_1+\epsilon(\xi), 
\end{eqnarray}

\noindent \textit{$\epsilon(\cdot) \sim N(0,\sigma^2)$, where $\sigma^2=1 \times 10^{-4}.$
$\hat{f}$ possesses additive noise with variance $\sigma^2.$ 
This function is a test function used in \cite{Nesterov} we have modified so that there is a distinct AS.
We considered $P=50$ and $j=5$, note $D=1$, and the non-stochastic $f$ is convex. 
Here, we have a $j$-dimensional $\A$ spanned by the standard basis vectors $e^i$, $i=1,\ldots,j$.
Note that the minimum of $f$, $\lambda^*$ is given by $\lambda^*_i = 1-i/(j+1)$ for $i=1,\ldots,j$ and $\lambda^*_i$ is arbitrary for $i=j+1,\ldots, P$ and has minimizer $f^*=-1/2(1-1/(j+1))$ \cite{Nesterov}. 
Also, $L_1=4$ \cite{Nesterov}.
We performed 50 trials each of STARS using exact hyperparameters, ASTARS using exact active hyperparemters as well as the exact AS, and
using exact hyperparameters but learned AS, we performed FAASTARS.
A maximum iteration count of $3P^2=7500$ was used.
AS's were relearned every $2P=100$ steps using an eigenvalue threshold $\tau = 0.999$.
Noise regularization at the level of $\sigma^2$ improved results.
We do not present hyperparameter learning in this example -- further extensions to our methods are needed and discussed in our next paper.}

\begin{figure}[h]
	\centering
	\includegraphics[scale=.175]{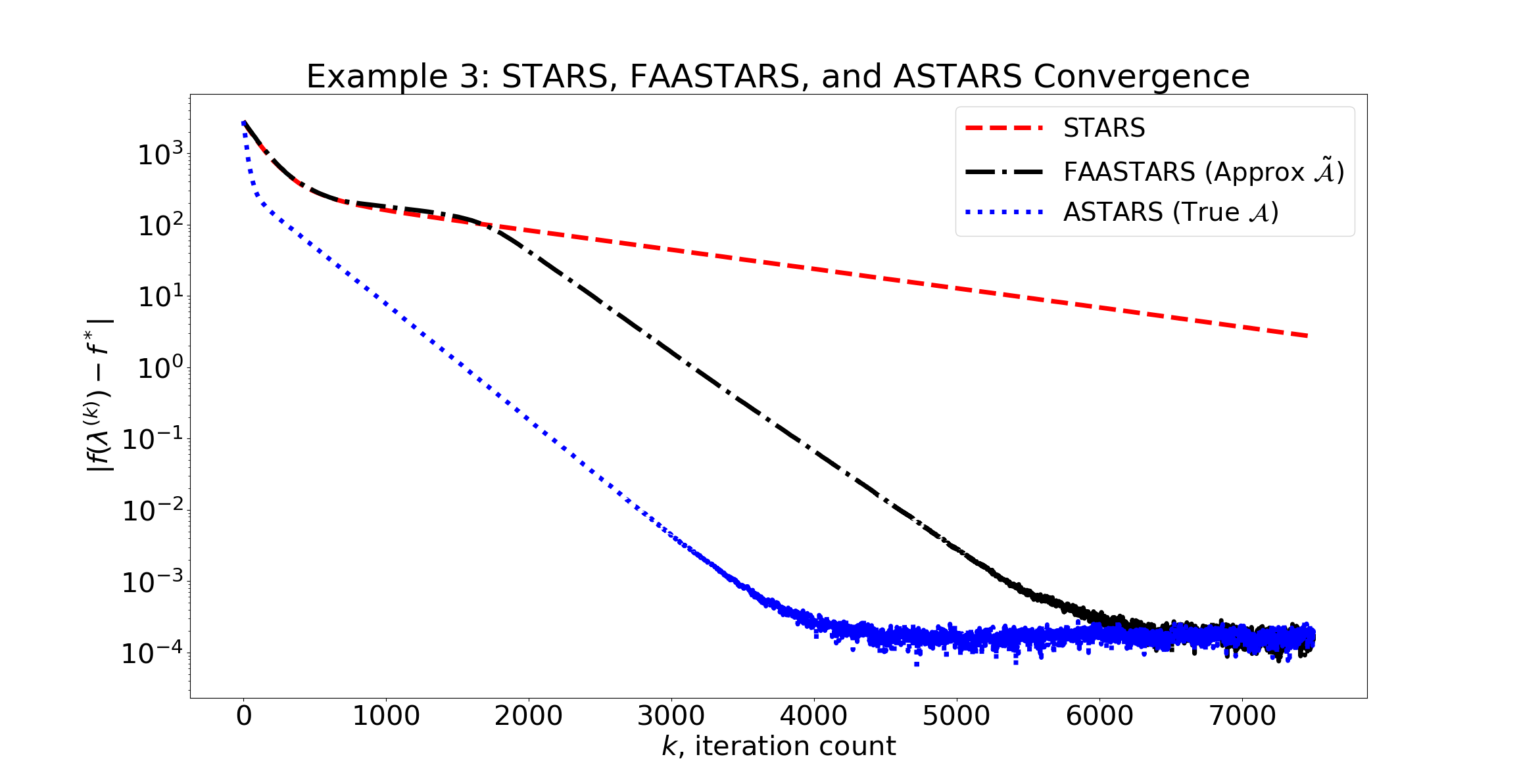}
	\caption{We show the convergence of STARS, FAASTARS, and ASTARS (without hyperparameter learning) for Example 3.}
	\label{fig:ex3}	
\end{figure}

\noindent \textbf{Example 4. (Full sphere)} \textit{Let} $\hat{f}: \Lambda=\R^P \to \R =\D,$

\begin{eqnarray} \label{eq:106}
\hat{f}(\lambda; \xi):=\sum_{i=1}^P \omega_i \lambda_i^2 + \epsilon(\xi), 
\end{eqnarray}

\noindent \textit{$\epsilon(\cdot) \sim N(0,\sigma^2)$, where $\sigma^2=1 \times 10^{-5}.$ We took $P=10$ and $\omega_i=1$, $i=1,\ldots,P$ so $\hat{f}$ the \textit{sphere function} in $\R^P$ with additive noise.
Note $L_1=2$ here.
We examined the extent to which the presented theoretical bounds for $L_1$ estimates hold in practice.
Define $K_1>0$ as the scale factor for which $L_1^2 = K_1 \hat{L_1}^2$.
In STARS Theorem 4.5 (in our technical report), we find the requirement $0<K_1<4$.
We write $\hat{L_1}=c \, L_1$, where $c=1/\sqrt{K_1}$ and we have $\frac{1}{2}<c<\infty$.
We used STARS to minimize \eqref{eq:106} where we used the correct $\sigma^2$ (not estimated) and fixed $\hat{L_1}=c \, L_1$ for $c=0.1, \, 0.2, \, 1,$ and $4$.
We performed 100 trials for each value $c$ and a maximum iteration count of $2P^3=2000$.
Note that we allow $c<1/2$; in some cases we find $c$ may be less than $1/2$ and STARS with estimated hyperparameters will still converge. In general, overestimation of $L_1$ ($c>1$) slow convergence and as $c \to \infty$, perturbations are not taken in $\Lambda$, and the value of $\hat{f}$ cannot be meaningfully changed at all. Underestimation of $L_1$ ($c<1$) may in fact improve convergence in some cases (as we see in this example), but as $c \to 0^+$, STARS will eventually diverge.}

\begin{figure}[h]
	\centering
	\includegraphics[scale=0.175]{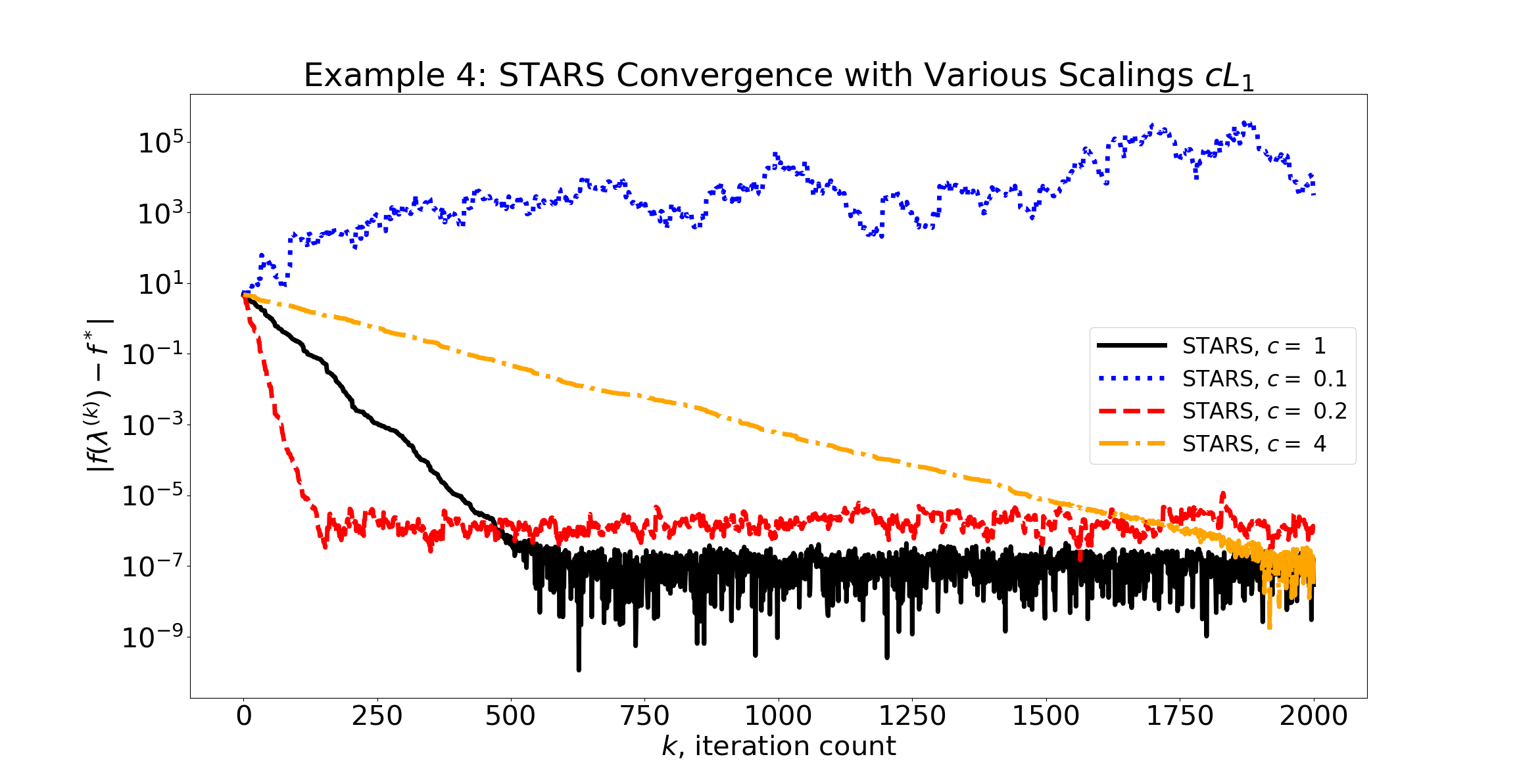}
	\caption{We show the convergence of STARS for various $c \, L_1$ for Example 4.}
	\label{fig:ex4}	
\end{figure}

\noindent \textbf{Example 5. (Nesterov-inspired, Nesterov 2)} \textit{Let} $\hat{f}: \Lambda=\R^P \to \R =\D,$

\begin{eqnarray} \label{eq:107}
\hat{f}(\lambda; \xi)=\sum_{i=1}^{P} 2^{(-1)^{(i-1)} (i-1)}\lambda_i^2+\epsilon(\xi), 
\end{eqnarray}

\noindent \textit{$\epsilon(\cdot) \sim N(0,\sigma^2)$, where $\sigma^2=1 \times 10^{-3}.$ 
We note $f$ is convex. 
We considered $P=10$ and note $D=1$. 
Note that the minimum of $\hat{f}$ is given by $0 \in \Lambda$. 
Here, as $i$ increases, terms in $\hat{f}$ become either more important or less important, depending on whether $i$ is even or odd.  
We take an initial iterate similar to prior examples. Also note  $L_1=2^{P+1}$ as long as $P$ is odd; if $P$ is even then $L_1=2^P$. Thus, with $P=10$, we have $L_1=1024$.
Here, the determination of which variables are active depends completely on one's choice of threshold. 
We show convergence of  STARS compared with FAASTARS using fixed active dimensions $\tilde{j}=2, \, 4,$ and $8$ in turn.
We performed 25 trials for each method and used a maximum iteration count of $5P^3=5000$.
}
\newpage

\begin{figure}[h]
	\centering
	\includegraphics[scale=0.175]{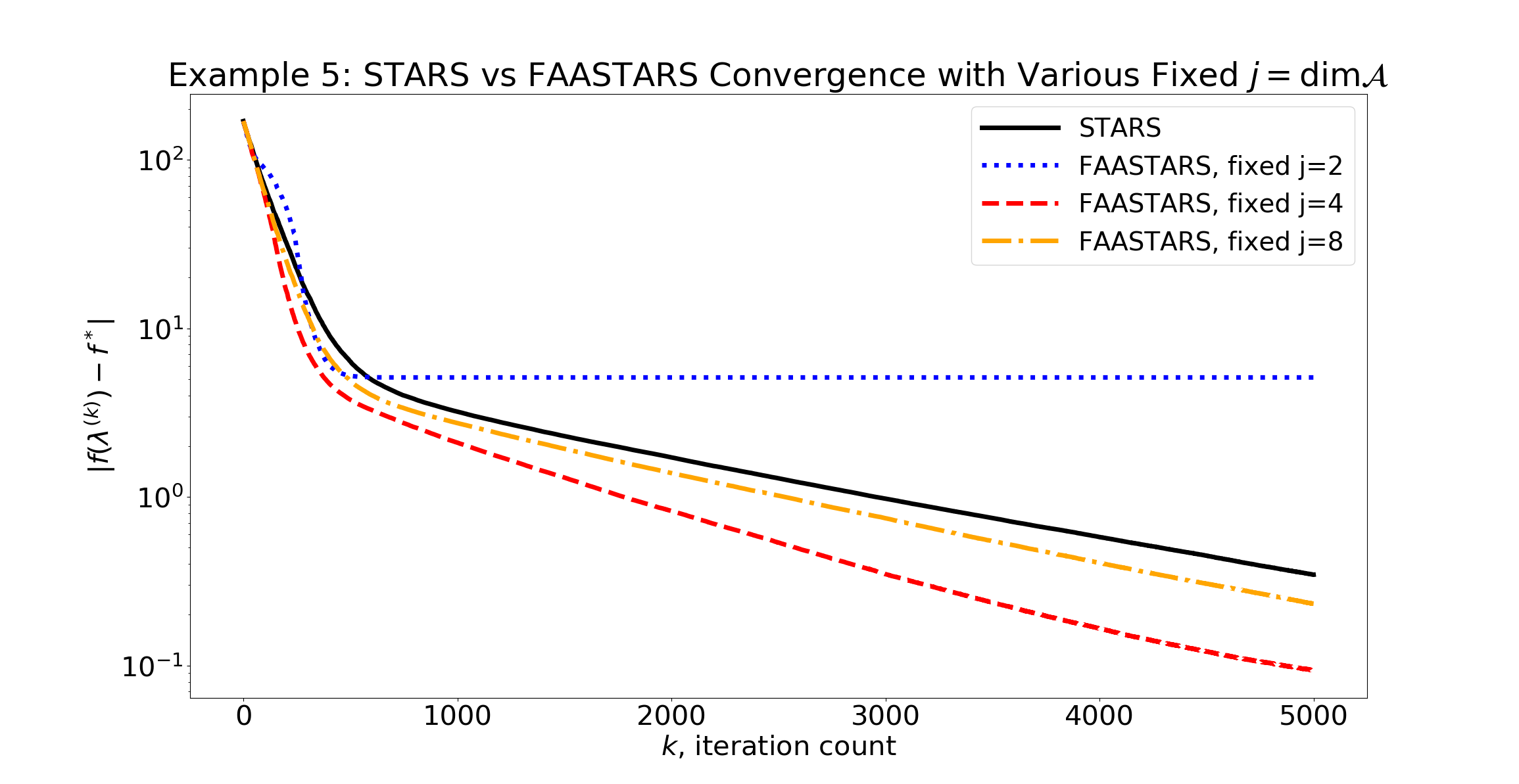}
	\caption{We show the convergence of STARS versus FAASTARS for various fixed $\tilde{j}$ for Example 5.}
	\label{fig:ex5}	
\end{figure}

In our examples, we see the computational benefit in stepping in the AS instead of in the full variables. Recall that the hyperparameters in STARS are dimension-dependent, so anytime an AS resolves $\hat{f}$ well -- which occurred in the above examples by obvious design -- we expect ASTARS in the exact active variables to converge more quickly than STARS in full variables.  

\subsection{Software Dissemination}
A Python 3.7 package called ASTARS was used to produce results in this section, and is open-source and publicly-available online \cite{ASTARS}.
The ASTARS package has the functionality to perform STARS, ASTARS, and FAASTARS.
STARS is not otherwise publicly-available, to our knowledge.
The algorithms used here which are open-source and publicly-available online include the AS software of Constantine et al \cite{AS2} -- updated in Python 3.6 (required for ASTARS package) by Varis Carey \cite{AS3} -- and ECNoise by Wild and Mor\'e \cite{ECN}.

\section{Conclusion and Discussion} \label{s:Conc}

We presented combinations and modifications made to well-documented algorithms including STARS \cite{CW}, Monte Carlo AS learning \cite{ConstantineMC}, noise variance learning \cite{MW}, and Lipschitz constant learning \cite{Calliess} to produce the fully-automated ASTARS algorithm. In addition, we presented several model  problems that were used for testing ASTARS and FAASTARS.

There is not any guarantee that a general stochastic-free mapping $f: \Lambda \to \D$ permits dimension reduction via AS. AS methods may fail to improve STARS for many reasons -- sometimes occurring in combination with each other -- including (nearly) equal importance of modeled parameters, a $\nabla f$ that is poorly approximated by surrogates, or too much noise. Regardless, in the case that no AS is clearly defined, recall that ASTARS is equivalent to STARS in full variables. 

If an AS exists, we observed that performing ASTARS and FAASTARS provided computational savings in our numerical examples when compared to STARS, since its convergence outpaces STARS convergence on average; see our technical report for detailed theoretical statements. 
We note that at times, any of the presented algorithms may take a very \textit{lucky} step in $\Lambda$, causing quick convergence. After all, the considered DFO algorithms depend on random (but smoothed) steps. 

We also observed that it is possible to learn the AS of the functions considered using the samples obtained from deterministic ECNoise and STARS iterates, which greatly reduces the usual computational expense of an AS computation. 
This result is not obvious, since AS discovery typically relies on (many) random, iid samples in $\Lambda$.

Sometimes, when FAASTARS has minimized $\hat{f}$ in $\tilde{\A}$ as much as possible, there may be remaining variables in the inactive subspace that are not minimized, as in Example 5 above.
Indeed, we saw that for various fixed AS dimensions $\tilde{j}$, ASTARS may behave almost identically to STARS, or much worse than STARS, depending on $\tilde{j}$.
Even without a fixed AS dimension, we still find that if FAASTARS determines $\tilde{j}< j$ (or $\tilde{j}>j$), iterates may not provide enough information for FAASTARS to update $\tilde{j}$ closer to $j$, incorrectly fixing $\tilde{j}$, and causing behavior identical to the flat-lining we produced in Example 5.

At its core, this problem is directly related to the choice of eigenvalue threshold $\tau$, which determines how many directions to include in $\tilde{\A}$.
When $\tau$ is fixed and samples are not informative enough to increase (or decrease) $\tilde{j}$, $\tilde{\A}$ cannot be substantially updated based on local information.

In an upcoming follow-up paper, we address the flat-lining behavior witnessed in some numerical examples by introducing a method of \textit{adaptive thresholding}, which will change $\tau$ if and when flat-lining occurs.
Other extensions of the ASTARS method, such as alternative weighting schemes and other approaches to address flat-lining, may also improve convergence and behavior, and will be considered as well.

In their unpublished manuscript, the authors in \cite{ARDFDS} present an accelerated version of DFO algorithms of Nesterov's classic approach called $RS_\mu$ in \cite{Nesterov} by leveraging techniques of \textit{mirror descent}.
(Note, $RS_\mu$ and STARS exhibit identical complexities \cites{Nesterov, CW}.)
As authors in \cite{ARDFDS} observe, complexity in $M$ -- the maximum number of iterations -- is equivalent to the complexity of oracle calls (i.e., the number of times we must approximate certain directional derivatives). These complexities are then proportional to the number of $\hat{f}$ evaluations, as well. Since we postulate that $\hat{f}$ calls are expensive, we always seek methods that evaluate our map as few times as possible. We will note, though, that we also know that there is a trade-off between fewer evaluations (samples) and the quality of the numerically estimated AS.

For a more concise comparison between the complexity results here and in \cite{CW, Nesterov, ARDFDS}, we let $\epsilon = \epsilon_{\text{tol}}$ and recall $R^2=||\lambda^{(0)}-\lambda^*||_2^2$ denotes the distance from the true minimizer to our initial iterate.
Recall, for all STARS-based methods, we achieve the main complexity statement in \eqref{eq:066} so long as 

\begin{equation} \label{eq:var-req}
\sigma \leq \frac{\sqrt{2K_2}(2-\sqrt{K_1})\epsilon_{\text{tol}}}{8(P+4)C_4},
\end{equation}

\noindent where  $K_1>0$ and $K_2 >0$ are the scale factors for which $L_1^2 = K_1 \hat{L_1}^2$ and $\sigma^2 = K_2 \hat{\sigma}^2$.

If we satisfy \eqref{eq:var-req}, then we obtain the results shown in Table \ref{table:dfo-complexity} below. Other methods will have $\sigma^2$ appearing in their complexity results, since no bounds are assumed directly on the variance in the noise.

\begin{table}[h]
\caption{Complexity of DFO Algorithms}
\begin{tabular}{c c } 
 \hline
\\ \textbf{ DFO Algorithm} & \textbf{Complexity}, $\mathcal{O}(\cdot)$  \\ [1ex] 
 \hline \\
 $RS_\mu$, \cite{Nesterov} and STARS, \cite{CW}  & $\frac{L_1 P R^2}{\epsilon}$ \\  [1ex]
 STARS (estimated hyperparameters) & $\frac{L_1 P R^2}{\sqrt{K_1}(2-\sqrt{K_1})\epsilon} $  \\  [1ex]
 ASTARS (estimated hyperparameters) & $ \frac{L_1 j R^2}{\sqrt{K_1}(2-\sqrt{K_1})\epsilon} $  \\ [1ex]
  FAASTARS (estimated hyperparameters) & $ \frac{L_1 \tilde{j} R^2}{\sqrt{K_1}(2-\sqrt{K_1})\epsilon} $  \\ [1ex]
 ARDFDS, \cite{ARDFDS} & $\max \left\{ P^{\frac{1}{2}+\frac{1}{q}}\sqrt{\frac{L_1 R^2}{\epsilon}}, \frac{P^{\frac{2}{q}}\sigma^2 R^2}{\epsilon^2} \right\}$ \\  

\end{tabular}
\label{table:dfo-complexity}
\end{table}

The methods in \cite{ARDFDS} achieve accelerated convergence by utilizing a different random direction in forming their gradient oracle, drawn randomly from the unit hypersphere in $\Lambda$ in a mirror descent scheme.
For future work, we propose investigating whether dimension reduction could accelerate the methods in \cite{ARDFDS}.

As a final note of caution, anytime the variance of the noise $\sigma^2$ approaches the magnitude of $f$ values, we expect failure of all methods presented. The assumptions we made in Chapter 1 -- assumptions ubiquitous in the DFO literature -- forbid noise of this order, and with good reason: anytime the noise is on the order of function evaluations, it becomes difficult (and eventually impossible) to distinguish the true signal from the effects of noise. Filtering or smoothing methods must be used in this scenario, which is outside the scope and focus of this dissertation. 

\newpage


\bibliographystyle{siam}
\bibliography{ref/references,ref/references-code}

\newpage

\begin{center}
\large \textbf{Technical Report}
\end{center}

\vspace{.25cm}

\setcounter{equation}{0}

We present theoretical results regarding the convergence of the methods provided in the preceding section.
In the first part of this section, we provide key results needed for proofs in the remaining parts.
In the second part of this section, we prove a series of modified STARS results culminating in a statement about the convergence of the algorithm with approximate hyperparameters. 
In the third part, we prove a series of results showing the convergence of ASTARS with exact hyperparameters and an exact AS.
Finally, we prove a series of FAASTARS results culminating in a statement about the convergence of FAASTARS with approximate hyperparamters and with an approximate AS.

Broadly, our contribution is showing: (1) STARS will still converge if $L_1$ and $\sigma$ are unknown and replaced with estimates $\hat{L_1}$ and $\hat{\sigma}$ in the formation of STARS hyperparameters; (2) ASTARS will converge with exact hyperparameters and an exact AS; and (3) FAASTARS will analogously converge also with uncertain hyperparameters, and even with an approximated $\tilde{\A}$ in place of the true $\A$.

\input{updated_proofs.tex}

\end{document}

%% file: updated_proofs.tex
\subsection{Preliminaries}

We provide the equations and results needed from \cite{Nesterov}, which are also summarized in \cite{CW}.
We also modify certain key results needed for ASTARS theoretical arguments.

We focus only on the case of additive noise in $\hat{f}$ and we assume $f$ is convex and differentiable in $\Lambda$.
Recall that $\Lambda = \R^P$ with $\lambda$'s denoting vectors, $\lambda \in \Lambda$; at times $u,x,y,z \in \Lambda$ will denote vectors, too. Also, we note that for the FAASTARS , we leave out the optional steps in FAASTARS related to the updating of $L_1^\text{init}$.
We will fix our approximation to $L_1$ at the beginning of FAASTARS using the samples formed via ECNoise \cite{MW}; i.e., $\hat{L_1}=L_1^\text{init}$ throughout.
We also note that we will not consider adaptive sampling methods for learning $\tilde{\A}$, so $\tilde{\A}$ will be fixed after FAASTARS Phase 2 (STARS burn-in phase).

We assume the true signal $f$ is convex and differentiable -- this is Assumption 4.1 in \cite{CW}. The signal we access is $\hat{f}$, which has additive noise so $\hat{f}(\lambda;\xi):=f(\lambda)+\epsilon(\xi)$, where $\mathbb{E}_\xi(\epsilon(\xi))=0$ $\forall \xi$ and $0<$Var$_{\xi}(\epsilon(\xi))=\sigma^2<\infty$ $\forall \xi$. These assumptions make up Assumption 4.2 in \cite{CW}. First, let

\begin{equation} \label{eq:01}
f_\mu(\lambda) := \mathbb{E}_u [f(\lambda + \mu u)], \, u \in \R^P, \, u_i \sim N(0,1), \, i= 1, \ldots, P,
\end{equation}

\noindent which is the expectation of the Gaussian-smoothed form of $f$.

Also, note for a direction $u=u_\A \in \A \subset \R^P$, the last $P-j$ components of $u_\A$ are zero and $u_i \sim N(0,1)$ for $i=1,\ldots , j$.
Hence, \eqref{eq:01} becomes

\begin{equation} \label{eq:active-01}
f_\mu^\A(\lambda) = \mathbb{E}_{u_\A} [f(\lambda + \mu u_\A)], \, u_\A \in \A \subset \R^P, \, u_i \sim N(0,1), \, i= 1, \ldots, j,
\end{equation}

\noindent and similarly for $\tilde{\A}$ with $\dim \tilde{\A}=\tilde{j}$.

The existence of $L_1$ implies

\begin{equation} \label{eq:02}
\left|f(y)-f(x)-\langle \nabla f(x), y - x \rangle \right| \leq \frac{L_1}{2}||x -y||^2 \quad \forall x,y \in \Lambda,
\end{equation}

\noindent Now let $\lambda^*$ denote a global minimizer of $f$. Then

\begin{equation} \label{eq:03}
||\nabla f(\lambda)||^2 \leq 2L_1 (f(\lambda)-f(\lambda^*)) \quad \forall \lambda \in \Lambda,
\end{equation}

\noindent proven in \cite{Zhou}.
A differentiable function $f$ is convex iff

\begin{equation} \label{eq:04}
f(y) \geq f(x) + \langle \nabla f(x) , y-x \rangle \quad \forall x,y \in \Lambda,
\end{equation}

\noindent Note that this implies that the left-hand side of \eqref{eq:02} is nonnegative, so long as $f$ is convex. An interpretation of \eqref{eq:04} is that a convex function is always underestimated by its linear approximation.

We now present the needed results on Gaussian smoothing from \cite{Nesterov}, also presented in \cite{CW}. First, some notation; for $\mu > 0$ and $u$ a Gaussian random vector (as in \eqref{eq:01} or \eqref{eq:active-01}),

\begin{equation} \label{eq:05}
g_\mu (x) := \frac{f(x+\mu u)-f(x)}{\mu}u,
\end{equation}

\noindent which is a first-order approximation to the directional derivative of $f$ in the direction of $u.$
If $u = u_\A$, then $g_\mu^\A$ will be the same object as \eqref{eq:05}, but with $u=u_\A$. 
Then $g_\mu^\A$ estimates the directional derivative of $f$ for directions $u_\A$ strictly in $\A$.

Now for $p \geq 0$, let

\begin{equation} \label{eq:06}
M_p(u) := \mathbb{E}_u \left( ||u||^p \right),
\end{equation}

\noindent the $p$-th moment of the norm of the random vector $u$, $||u||$.
Here and throughout this section, $||\cdot||$ denotes the 2-norm in $\R^P$. 
Let $u$ continue to denote a random vector as in \eqref{eq:01}.
For $M_p(u)$ defined in \eqref{eq:06} above,

\begin{equation} \label{eq:07}
M_p(u) \leq (P)^{p/2}, \quad p \in [0,2] \quad \text{and} \quad M_p(u) \leq (P+p)^{p/2}, \quad p>2,
\end{equation}

\noindent where we recall $P = \dim \Lambda$.

Now for $u=u_\A$ (as in \eqref{eq:active-01}) and $j<P$, we have a sharper but analogous result for these moments, which we explain by considering the case of $p=2$.
From \cite{Nesterov}, we have

\begin{equation} \label{eq:moment-act-1}
M_2(u) = \frac{1}{c} \int_{\R^P} ||u||^2 \, e^{-\frac{1}{2}||u||^2}\, \, du=B^{-1},
\end{equation}

\noindent where $B$ is the matrix specifying the norm $||u||^2 = \langle Bu, u \rangle$ for a given inner product $\langle \cdot, \cdot \rangle$ in $\R^P$ and $c$ is a normalization factor.
Throughout, we will let $B=I_P$ and use the Euclidean inner product, so that $||\cdot||$ will continue to denote the standard (Euclidean) 2-norm in $\Lambda$.

Note that since $u_\A \in \A$, we only integrate over the $j$ parameters corresponding to $\A$ since the components of $u_\A$ in $\I$ are fixed and zero in expectation.
Also, since $u_\A$ is zero in its last $P-j$ components, we can compute the norm of $u_\A$ truncated after its $j$-th entry in $\R^j$ instead of in $\R^P$.
In particular, define $\underline{u_\A}=(u_\A)_{1:j}$ and let $||\cdot||_{\underline{\A}}$ denote the 2-norm in $\R^j$
We have $||u_\A|| = ||u_{\underline{\A}}||_{\underline{\A}}$, and therefore

\begin{equation} \label{eq:moment-act-2}
M_2(u_\A) = \frac{1}{c} \int_{\R^j} ||\underline{u_{\A}}||^2_{\underline{\A}} \, e^{-\frac{1}{2}||\underline{u_{\A}}||^2_{\underline{\A}}} \, \,d \underline{u_{\A}} = I_{j}.
\end{equation}

As in \cite{Nesterov}, taking the inner product of the left-hand and right-hand sides with $I_j$ shows $M_2(u_\A)=j$.
Using similar arguments, one can prove the generalized bounds

\begin{equation} \label{eq:moment-act-3}
M_p(u_\A) \leq (j)^{p/2}, \quad p \in [0,2] \quad \text{and} \quad M_p(u_\A) \leq (j+p)^{p/2}, \quad p>2,
\end{equation}

\noindent which involves making similar changes to those outlined above.
In particular, one must substitute $u$ with $\underline{u_{\A}}$ and $||\cdot||$ with $||\cdot||_{\underline{\A}}$ and rewrite the integrals corresponding to the expectation over all $\R^P$ as integrals over $\R^j$, as we presented for $M_2$.

We now present Gaussian smoothing results in the case for $f$ convex.
If $f$ is convex, then

\begin{equation} \label{eq:08}
f_\mu (x) \geq f(x) \quad \forall x \in \Lambda,
\end{equation}

\noindent which can be verified by writing $\mathbb{E}_u(f(x+\mu u)) \geq \mathbb{E}_u(f(x)+\langle \nabla f(x), \mu u \rangle)$, where the inequality arises from applying the definition of convexity (\eqref{eq:04} with $y=x+\mu u$ and $x=x$).
Then since $f(x)$ is constant with respect to $u$ and $\mathbb{E}_u (\langle \nabla f(x), \mu u \rangle)=0$ (since each component of $u$ is zero-mean), we obtain $\mathbb{E}_u(f(x)+\langle \nabla f(x), \mu u \rangle)=f(x)$, verifying \eqref{eq:08}.
Also note that \eqref{eq:08} holds with $f_\mu=f_\mu^\A$ (as in \eqref{eq:active-01}).

If $f$ is convex and $f\in \mathcal{C}^{1,1}(\Lambda)$, which is the space of functions $f: \Lambda=\R^P \to \R$ that are continuously differentiable and possess a $L_1$ Lipschitz constant, then

\begin{equation} \label{eq:09}
f_\mu (x) -f(x) \leq \frac{L_1 \mu^2}{2} P \quad \forall x\in \Lambda,
\end{equation}

\noindent where we recall \eqref{eq:08} implies the left-hand side of \eqref{eq:09} above is nonnegative. 

To verify \eqref{eq:09}, we use the fact that $f$ is convex (i.e., \eqref{eq:04} with $y=x+\mu u$ and $x=x$) to write
$$0 \leq f(x+\mu u)-f(x)-\langle \nabla f(x),\mu u \rangle.$$

\noindent Applying \eqref{eq:02},
$$ f(x+\mu u)-f(x)-\langle \nabla f(x),\mu u \rangle \leq \frac{L_1 \mu^2}{2}||u||^2.$$
Applying the expectation in $u$ to both sides and using \eqref{eq:07} for the second moment of $||u||$, we obtain \eqref{eq:09}.
If $u=u_\A$ and we have $f_\mu^\A$ (as in \eqref{eq:active-01}), then using \eqref{eq:moment-act-3},

\begin{equation} \label{eq:active-09}
f_\mu^\A (x) -f(x) \leq \frac{L_1 \mu^2}{2} j \quad \forall x\in \Lambda.
\end{equation}

We will also work with the objects $\nabla f_\mu (x)$ and $\nabla f_\mu^\A (x)$.
To quote \cite{CW}, we note that with "$\nabla f_\mu (x)$ we denote the gradient (with respect to $x$) of the Gaussian approximation" $f_\mu(x)$ defined in \eqref{eq:01} and likewise for $\nabla f_\mu^\A (x)$, for $f_\mu^\A$ as in \eqref{eq:active-01}.
In \cite{Nesterov}, the authors obtain the form in \eqref{eq:grad-integral-form} below for $\nabla f_\mu$ by performing a substitution in the integral corresponding to the expectation in $u$ arising in the definition of $f_\mu$ and differentiate both sides in $x$, yielding

\begin{equation} \label{eq:grad-integral-form}
\nabla f_\mu(x) = \frac{1}{c} \int_{\R^P} \frac{f(x+\mu u)-f(x)}{\mu} u e^{-\frac{1}{2}||u||^2} \, du= \frac{1}{c} \int_{\R^P} g_\mu(x) e^{-\frac{1}{2}||u||^2} \, du,
\end{equation}

\noindent which shows by definition

\begin{equation} \label{eq:010}
 \nabla f_\mu (x) = \mathbb{E}_u \left( g_\mu (x) \right)  \quad \forall x \in \Lambda
\end{equation}

Analogously, we have

\begin{equation} \label{eq:active-010}
 \nabla f_\mu^\A (x) = \mathbb{E}_u \left( g_\mu^\A (x) \right)  \quad \forall x \in \Lambda
\end{equation}

\noindent by taking $u=u_\A$ in our arguments above.

If $f$ is differentiable at $x$ and $f \in\mathcal{C}^{1,1}$, 

\begin{equation} \label{eq:011}
\mathbb{E}_u \left( ||g_\mu(x)||^2 \right) \leq 2(P+4)||\nabla f(x)||^2+\frac{\mu^2}{2}L_1^2(P+6)^3 \quad \forall x\in \Lambda.
\end{equation}

To verify \eqref{eq:011}, it is helpful to first bound the quantity $\mathbb{E}_u \left( ||u||^2 \langle \nabla f(x), u \rangle ^2 \right)$.
We follow \cite{Nesterov}.
Applying Cauchy-Schwarz and write
$||u||^2 \langle \nabla f(x), u \rangle ^2 \leq ||\nabla f(x)||^2 ||u||^4$.
Here we could apply $\mathbb{E}_u(\cdot)$ and use \eqref{eq:07} to bound $M_4(u)$, yielding
$$\mathbb{E}_u \left( ||u||^2 \langle \nabla f(x), u \rangle ^2 \right)\leq  (P+4)^2 ||\nabla f(x)||^2 .$$

However, the authors in \cite{Nesterov} obtain a tighter bound by minimizing the integral form of $\mathbb{E}_u \left( ||u||^2 \langle \nabla f(x), u \rangle ^2 \right)$. 
The proof is technical but mainly involves parameterizing and then minimizing the argument of the exponential function appearing in the integral form associated with the expectation over $u$.
\cite{Nesterov} show
$$\mathbb{E}_u \left( ||u||^2 \langle \nabla f(x), u \rangle ^2 \right)\leq (P+4) ||\nabla f(x)||^2 .$$

For the case of $g_\mu^\A$, we can show

\begin{equation} \label{eq:active-011}
\mathbb{E}_u \left( ||g_\mu(x)||^2 \right) \leq 2(j+4)||\nabla f(x)||^2+\frac{\mu^2}{2}L_1^2(j+6)^3 \quad \forall x\in \Lambda.
\end{equation}

Justifying \eqref{eq:active-011} requires the same substitutions as in justifying \eqref{eq:moment-act-3}.
Again, one must substitute $u$ with $\underline{u_{\A}}$ and $||\cdot||$ with $||\cdot||_{\underline{\A}}$ and rewrite the integrals corresponding to the expectation over all $\R^P$ as integrals over $\R^j$.
Then $P$'s in the arguments above may be safely replaced with $j$'s as in \eqref{eq:active-011}.

Next, the authors in \cite{CW} introduce two more pieces of notation integral to the STARS process. First, for an iterate $k$, let

\begin{equation} \label{eq:012}
s_{\mu_k}(x^{(k)};u^{(k)},\xi_{k-1},\xi_{k}):= \frac{\hat{f}(x^{(k)}+\mu_k u^{(k)}; \xi_{k-1}) - \hat{f}(x^{(k)}; \xi_{k})}{\mu_k}u^{(k)},
\end{equation}

\noindent which is the "stochastic gradient-free oracle" to the directional derivative of $f$ (with noise) in the direction $u$. The authors also define the error between the oracle (the forward-difference approximation) and the true directional derivative of $f$ in the direction $u$ as

\begin{equation} \label{eq:013}
\mathcal{E}(\mu)= \mathcal{E}(\mu; x, u, \xi_1,\xi_2) = ||s_\mu(x;u,\xi_1,\xi_2) - \langle \nabla f(x),u \rangle u ||^2.
\end{equation}

Note when $u=u_\A$ in either \eqref{eq:012} or \eqref{eq:013}, we write $s^\A_{\mu_k}$ and $\mathcal{E}^\A(\mu)$ to emphasize that these values are computed for $u_\A \in \A$.
If we are working with $\tilde{\A}$ to approximate $\A$, one must replace $\A$ with $\tilde{\A}$ and $j$ with $\tilde{j}$ to obtain the analogous definitions for the approximate $\tilde{\A}$ case.

\subsection{STARS Convergence with Estimated Hyperparameters} \label{ss:stars-conv}

In the following, we follow \cite{CW} closely.
However, we provide more detail and in fact correct some minor details from their unpublished manuscript while generalizing their results to the case in which hyperparameters are estimated. 
In the latter sections, we build on results in \cite{CW} and \cite{ConstantineK}, a paper which contains crucial theoretical results regarding the approximation of active subspaces.

Let the positive, finite values $\hat{L_1}$ and $\hat{\sigma}$ denote estimators to the true (also positive and finite values of $L_1$ and $\sigma$. Recall, in our setting, we will assume $0<L_1<\infty$ and $0<\sigma^2<\infty.$ We won't let $L_1=0$ since that would imply $f$ is constant; we won't let $\sigma^2=0$, since that would imply zero noise. Then there exist $K_1> 0$ and $K_2 
> 0$ so that

\begin{equation} \label{eq:014}
L_1^2 = K_1 \hat{L_1}^2 \quad \quad  \sigma^2 = K_2 \hat{\sigma}^2.
\end{equation}

Note that if a $K_i<1$, $i=1$ or $2$, then we have overestimated the corresponding value, $L_1$ or $\sigma^2$ and similarly when $K_i>1$, the corresponding value has been underestimated. Hence, as a particular $K_i \to 1$, the corresponding estimate to either $L_1$ or $\sigma^2$ approaches the truth. Finally, note that when the true values $L_1$ and $\sigma^2$ are unknown, $K_1$ and $K_2$ are also generally unknown.

Below, we recall the approximate smoothing size and step length to replace the STARS hyperparameters in the case that $L_1$ and $\sigma$ are unknown and estimated by values $\hat{L_1}$ and $\hat{\sigma}$.

\begin{equation} \label{eq:015}
\reallywidehat{\mu^*}:=\left( \frac{8\hat{\sigma}^2 P}{\hat{L_1}^2(P+6)^3}\right)^{1/4} \quad \quad \reallywidehat{h}:=(4\hat{L_1}(P+4))^{-1}
\end{equation}

Shortly, we will precisely show how the bound $\mathcal{E}(\hat{\mu^*})$ is just a modification of the bound on $\mathcal{E}(\mu^*)$ proven in \cite{CW}. 
With this point of view, one will see our choice of $\hat{\mu^*}$ is the best that we can do with uncertain estimators to $L_1$ and $\sigma$. We note that both in \cite{CW} and here in this work, $h$ (and also $\hat{h}$) is a fixed choice, not necessarily optimal nor sub-optimal; more investigation could be done into the step length, but $\hat{h}$ works fine in our numerical experiments.

It is helpful to quote \cite{CW} here: "Our goal is to find $\mu^*$ that minimizes an upper bound on" $\mathbb{E}_{u,\xi_1,\xi_2}(\mathcal{E}(\mu))$, where we have $\mathcal{E}$ defined in \eqref{eq:013} and the expectation is taken over the random vector $u$, as well as two draws of additive noise from the two function evaluations that occur in $s_{\mu}$. The noise draws are denoted as $\epsilon(\xi_1)$ and $\epsilon(\xi_2)$.

The major difference between our result and the STARS result is that we will use estimations to $L_1$ and $\sigma^2$ rather than their true values, which we have postulated lack of access to.

\textbf{STARS Theorem 4.3 (Modified):} \textit{We assume random vectors $u^{(k)}$ are drawn according to \eqref{eq:01}; $f\in \mathcal{C}^{1,1}(\Lambda)$ and $f$ is convex; and that the i.i.d. noise draws $\epsilon(\xi)$ are additive, zero mean, with bounded variance $\sigma^2$ for all $\xi$.
Let $u$ be drawn in the fashion described in \eqref{eq:01}. If a smoothing stepsize is chosen as $\mu=\hat{\mu^*}$ in \eqref{eq:015}, then for any iterate $x \in \Lambda$ and random vector $u$, noting $K_1>0$ and $K_2>0$ as in \eqref{eq:014},}

\begin{equation} \label{eq:016}
\mathbb{E}_{u,\xi_1,\xi_2} \left( \mathcal{E}\left(\hat{\mu^*}\right) \right) \leq \frac{K_1+K_2}{\sqrt{2K_1 K_2}}\sigma L_1\sqrt{P(P+6)^3}.
\end{equation}

\textit{Proof:} Let $u\in \Lambda$ be a random vector as in \eqref{eq:01}, $x \in \Lambda$ denote a general STARS iterate, and $\epsilon(\xi_1)$ and $\epsilon(\xi_2)$ denote the two i.i.d. draws of the additive noise in $\hat{f}$ which appear in $\mathcal{E}(\mu)$, \eqref{eq:013}. Plugging equation \eqref{eq:012} into equation \eqref{eq:013}, we obtain

\begin{equation} \label{eq:017}
\mathcal{E}(\mu) = \left| \left|  \frac{f(x+\mu u)+\epsilon(\xi_1)-\left(f(x)+\epsilon(\xi_2)\right)}{\mu}u  - \langle \nabla f(x),u \rangle u  \right| \right|^2.
\end{equation}

Rearranging,

\begin{equation} \label{eq:018}
\mathcal{E}(\mu) =\left| \left|  \left(  \frac{\left(\epsilon(\xi_1) - \epsilon(\xi_2)\right)+ \left(f(x+\mu u)-f(x) - \left\langle \nabla f(x),\mu u \right\rangle \right)}{\mu}  \right)u \right| \right|^2.
\end{equation}

We have

\begin{equation} \label{eq:019}
\mathcal{E}(\mu) \leq \frac{X^2}{\mu^2}||u||^2, \, \, \text{where} \, \, X:=\left(\epsilon(\xi_1) - \epsilon(\xi_2)\right)+ \left(f(x+\mu u)-f(x) - \left\langle \nabla f(x),\mu u \right\rangle \right).
\end{equation}

Expanding the form of $X$,

\begin{equation} \label{eq:020}
\begin{split}
X^2 &=\epsilon(\xi_1)^2-2\epsilon(\xi_1)\epsilon(\xi_2) + \epsilon(\xi_2)^2 + 2(\epsilon(\xi_1) - \epsilon(\xi_2))\left(f(x+\mu u)-f(x) - \left\langle \nabla f(x),\mu u \right\rangle \right) \\
&+ \left(f(x+\mu u)-f(x) - \left\langle \nabla f(x),\mu u \right\rangle \right)^2.
\end{split}
\end{equation}

We begin by examining the expectation of $X^2$ with respect to the two stochastic noise draws. Recall that $\mathbb{E}_\xi(\epsilon(\xi))=0$ and Var$(\epsilon(\xi))=\sigma^2>0$ for all draws $\xi,$ hence, we also have $\mathbb{E}_\xi(\epsilon(\xi)^2)=\sigma^2$ for all draws $\xi$.

\begin{equation} \label{eq:021}
\mathbb{E}_{\xi_1,\xi_2}(X^2)=2\sigma^2+\left(f(x+\mu u)-f(x) - \left\langle \nabla f(x),\mu u \right\rangle \right)^2.
\end{equation}

Noting that neither $u$ nor $\mu$ depends on noise draws $\xi$, we have 

\begin{equation} \label{eq:022}
\mathbb{E}_{u,\xi_1,\xi_2}(\mathcal{E}(\mu))=\mathbb{E}_u\left(\frac{\mathbb{E}_{\xi_1,\xi_2}(X^2)||u||^2}{\mu^2}\right)=\frac{\mathbb{E}_u\left(2\sigma^2||u||^2+\left(f(x+\mu u)-f(x) - \left\langle \nabla f(x),\mu u \right\rangle \right)^2||u||^2\right)}{\mu^2}.
\end{equation}

Now replacing $y$ with $y=x+\mu u$ and squaring both sides, we can re-write \eqref{eq:02} as

\begin{equation} \label{eq:023}
(f(x+\mu u)-f(x)-\langle \nabla f(x), \mu u \rangle)^2  \leq \frac{L_1^2}{4}\mu^4|| u||^4 .
\end{equation}

Now, combining \eqref{eq:022} and \eqref{eq:023}, we have 

\begin{equation} \label{eq:024}
\mathbb{E}_{u,\xi_1,\xi_2}(\mathcal{E}(\mu))\leq \frac{1}{\mu^2}\mathbb{E}_u\left(2\sigma^2||u||^2+\frac{L_1^2}{4}\mu^4|| u||^6\right).
\end{equation}

Using the bounds on the moments $M_p$ of $||u||^p$ given in \eqref{eq:07}, using $p=2$ and $p=6$, we have

\begin{equation} \label{eq:025}
\mathbb{E}_{u,\xi_1,\xi_2}(\mathcal{E}(\mu)) \leq \frac{2\sigma^2}{\mu^2}M_2+\frac{ L_1^2 \mu^2}{4}M_6 \leq \frac{2\sigma^2P}{\mu^2}+\frac{ L_1^2\mu^2(P+6)^3}{4}.
\end{equation}

We have

\begin{equation} \label{eq:026}
\mathbb{E}_{u,\xi_1,\xi_2}(\mathcal{E}(\mu)) \leq (2\sigma^2P)\frac{1}{\mu^2}+\left( \frac{L_1^2(P+6)^3}{4}\right)\mu^2.
\end{equation}

The authors in \cite{CW} observe that the right-hand side of the above inequality is uniformly convex for $\mu>0$, taking the form $t(\mu)=a\mu^{-2}+b\mu^2$ for positive constants $a$ and $b$; calculus shows that the minimizer of $t$ for $\mu>0$ is $\mu^*:=\left(a/b\right)^{1/4}$ with $t(\mu^*)=2\sqrt{ab}$.

Using $a=2\sigma^2P$ and $b=(L_1^2(P+6)^3)/4$, we recover $\mu^*=\left( \frac{8\sigma^2 P}{L_1^2(P+6)^3}\right)^{1/4}$, the optimal (in the sense of minimizing the upper bound on $\mathcal{E}$) smoothing step length proven in \cite{CW}. Our optimal choice of smoothing, given the information we have available, will require us to swap out $L_1$ and $\sigma^2$ in $\mu^*$ with their estimates, $\hat{L_1}$ and $\hat{\sigma}^2$, recovering $\hat{\mu}^*$ \eqref{eq:015}. This particular choice $\mu=\hat{\mu^*}$ can be plugged into \eqref{eq:026}, which gives us the bound

\begin{equation} \label{eq:027}
\mathbb{E}_{u,\xi_1,\xi_2}(\mathcal{E}(\hat{\mu^*}))
\leq \frac{K_1+K_2}{\sqrt{2K_1 K_2}}\sigma L_1\sqrt{P(P+6)^3}
,
\end{equation}

\noindent our main result.  $\blacksquare$


Note that we can recover the exact result of STARS Theorem 4.3 by taking $K_1=K_2=1$, the case in which our hyperparameters are estimated exactly.

We next derive an upper bound on $\mathbb{E}(||s_{\mu_k}||^2)$, where $\mathbb{E}$ will now denote the expectation over every noise draw and random vector used in STARS up to (and including) the $k$-th iterate; that is, the expectations are now taken with respect to $\xi_0,\ldots,\xi_{k}$ and $u^{(1)},\ldots,u^{k}$ unless stated otherwise.
Recall that $s_{\mu_k}$, the stochastic gradient-free oracle, was defined in \eqref{eq:012}.
We prove this result in a modified Lemma very similar to STARS Lemma 4.4 in \cite{CW}; similarly to the result above, we will only use the estimated values to $L_1$ and $\sigma$.

\textbf{STARS Lemma 4.4 (Modified):}
\textit{We assume random vectors $u^{(k)}$ are drawn according to \eqref{eq:01}; $f\in \mathcal{C}^{1,1}(\Lambda)$ and $f$ is convex; and that the i.i.d. noise draws $\epsilon(\xi)$ are additive, zero mean, with bounded variance $\sigma^2$ for all $\xi$. If we use $\mu_k = \hat{\mu^*}$ as in \eqref{eq:015} for all STARS iterates $k$, then noting $K_1>0$ and $K_2>0$ as in \eqref{eq:014}, STARS generates steps satisfying}

\begin{equation} \label{eq:028}
\mathbb{E}(||s_{\mu_k}||^2)\leq 2(P+4)||\nabla f(x^{(k)})||^2 + \frac{3K_1+K_2}{\sqrt{2K_1K_2}}  L_1\sigma\sqrt{P(P+6)^3}.
\end{equation}

\textit{Proof:} First, we set $\mu_k=\hat{\mu^*}$. For a STARS iterate $k$, let

\begin{equation} \label{eq:029}
g_0(x^{(k)}):=\langle \nabla f(x^{(k)}), u^{(k)} \rangle u^{(k)},
\end{equation}

\noindent which is the exact directional derivative of $f$ in the direction of $u$ at the point $x^{(k)} \in \Lambda$. We can use this notation to re-express \eqref{eq:027} as

\begin{equation} \label{eq:030}
\mathbb{E}\left(||s_{\mu_k}||^2 - 2\langle s_{\mu_k},g_0(x^{(k)}) \rangle + ||g_0(x^{(k)})||^2\right)\leq \frac{K_1+K_2}{\sqrt{2K_1 K_2}}\sigma L_1\sqrt{P(P+6)^3},
\end{equation}

\noindent where we have also expanded $\mathcal{E}$, defined in \eqref{eq:013}. Recalling that all draws of the noise are zero-mean, the expectation of the oracle $s_{\mu_k}$ (defined in \eqref{eq:012}) with respect to the appearing noise draws $\xi_{k-1}$ and $\xi_{k}$ is given by

\begin{equation} \label{eq:031}
\mathbb{E}_{\xi_{k-1},\xi_k}(s_{\mu_k})= \frac{f(x^{(k)}+\mu_k u^{(k)})-f(x^{(k)})}{\mu_k}u^{(k)} =g_\mu(x^{(k)}),
\end{equation}

\noindent which is the (noise-free) first-order approximation to the directional derivative of $f$ in the direction of $u$, defined in \eqref{eq:05}. The linearity of $\mathbb{E}$ allows us to rewrite \eqref{eq:030} as

\begin{equation} \label{eq:032}
\mathbb{E}\left(||s_{\mu_k}||^2\right) \leq \mathbb{E} \left( 2\langle s_{\mu_k},g_0(x^{(k)}) \rangle - ||g_0(x^{(k)})||^2\right)+ C_1,
\end{equation}

\noindent where $C_1:=\frac{K_1+K_2}{\sqrt{2K_1 K_2}}\sigma L_1\sqrt{P(P+6)^3}$. The only term involving noise draws on the right-hand side of \eqref{eq:032} is $s_{\mu_k}$; thus, passing through the expectation with respect to all noise draws $\xi_k$, we can use our result in \eqref{eq:031} to write

\begin{equation} \label{eq:033}
\begin{split}
\mathbb{E}\left(||s_{\mu_k}||^2\right)  & \leq \mathbb{E} \left( 2\langle s_{\mu_k},g_0(x^{(k)}) \rangle - ||g_0(x^{(k)})||^2\right)+ C_1  \\ & =\mathbb{E}_{u^{(k)}}\left( 2\langle g_{\mu}(x^{(k)}),g_0(x^{(k)}) \rangle - ||g_0(x^{(k)})||^2\right)+ C_1.
\end{split}
\end{equation}

Adding and subtracting by $||g_\mu (x^{(k)})||^2$ inside of the $\mathbb{E}_{u^{(k)}}$ (and then factoring) in \eqref{eq:033} gives

\begin{equation} \label{eq:034}
\mathbb{E}\left(||s_{\mu_k}||^2\right)  \leq \mathbb{E}_{u^{(k)}}\left( -||g_0(x^{(k)})-g_\mu(x^{(k)}) ||^2+ ||g_0(x^{(k)})||^2\right)+ C_1.
\end{equation}

Using the linearity of $\mathbb{E}_{u^{(k)}}$ and observing that $-||x||^2 \leq 0$ for all $x \in \Lambda$ gives

\begin{equation} \label{eq:035}
\mathbb{E}\left(||s_{\mu_k}||^2\right)  \leq \mathbb{E}_{u^{(k)}}\left( ||g_0(x^{(k)})||^2\right)+ C_1.
\end{equation}

Recalling the result in \eqref{eq:011}, we have arrived at

\begin{equation} \label{eq:036}
\mathbb{E}\left(||s_{\mu_k}||^2\right)  \leq 2(P+4)||\nabla f(x^{(k)}) ||^2 + \frac{\mu_k^2 L_1^2}{2}(P+6)^3 + C_1.
\end{equation}

Equivalently, recalling \eqref{eq:014} -- which equates $L_1^2$ to $\hat{L_1}^2$ scaled by a positive constant $K_1$ --  we also have

\begin{equation} \label{eq:037}
\mathbb{E}\left(||s_{\mu_k}||^2\right)  \leq 2(P+4)||\nabla f(x^{(k)}) ||^2 + \frac{K_1 \mu_k^2 \hat{L_1}^2}{2}(P+6)^3 + C_1.
\end{equation}

Recall we have set $\mu_k=\hat{\mu^*}$ (from \eqref{eq:015}) in \eqref{eq:036} for all iterations $k$. Plugging in this value, we obtain

\begin{equation} \label{eq:038}
\mathbb{E}\left(||s_{\mu_k}||^2\right)  \leq 2(P+4)||\nabla f(x^{(k)}) ||^2 + C_2,
\end{equation}

\noindent where $C_2:=\frac{3K_1+K_2}{\sqrt{2K_1K_2}}L_1\sigma\sqrt{P(P+6)^3}=\frac{3K_1+K_2}{\sqrt{2}}  \hat{L_1}\hat{\sigma}\sqrt{P(P+6)^3}$, our main result. $\blacksquare$

Note that in a fashion analogous to our modification of STARS Theorem 4.3, we can recover the exact result of STARS Lemma 4.4 by taking $K_1=K_2=1$.

We can now present the final result which shows that STARS converges with estimates replacing the exact values for $L_1$ and $\mu$. We need just a bit more notation, borrowed directly from \cite{CW}.
Let $x^* \in \Lambda$ denote a minimizer with the associated stochastic-free function evaluation $f^* := f(x^*)$. Also, define $\mathcal{Q}_k :=\{\xi_0, \ldots, \xi_k\}$and  $\mathcal{U}_k :=\{u^{(1)}, \ldots, u^{(k)}\}$, which are two sets containing all random variables that appear in STARS up through iteration $k$.
Let $\phi_0:= f(x^{(0)})$ and $\phi_k := \mathbb{E}_{\mathcal{Q}_{k-1},\mathcal{U}_{k-1}}(f(x^{(k)}))$, $k \geq 1.$ Define $M\in \N$ as the total number of STARS iterates performed.

\textbf{STARS Theorem 4.5 (Modified):}
\textit{Let Assumptions 3.1, 4.1, and 4.2 hold -- here those assumptions mean that random vectors $u^{(k)}$ are drawn according to \eqref{eq:01}; $f\in \mathcal{C}^{1,1}(\Lambda)$ and $f$ is convex; and that the i.i.d. noise draws $\epsilon(\xi)$ are additive, zero mean, with bounded variance $\sigma^2$ for all $\xi$. Let $\{x^{(k)}\}_{k\geq 0}$ denote a sequence of STARS iterates formed using a fixed step length $h_k=\hat{h}$ and fixed smoothing $\mu_k = \hat{\mu^*}$ (both given in \eqref{eq:015}) for all STARS iterates $k$. \textbf{Finally, we require $0<K_1<4$ and $K_2>0$, the values defined in \eqref{eq:014}.} Then for any total number of STARS iterations $M$,}

\begin{equation} \label{eq:039}
\sum_{k=0}^M\frac{\phi_k-f^*}{M+1} \leq   \frac{4L_1(P+4)||x^{(0)}-x^*||^2}{\sqrt{K_1}(2-\sqrt{K_1})(M+1)} + \frac{4\sigma(P+4)}{\sqrt{2K_2}(2-\sqrt{K_1})}C_5,
\end{equation}

\noindent where $C_5:=\sqrt{K_1}\cdot0.036+\frac{3K_1+K_2}{16}\cdot 1.034$.

\textit{Proof:} For a STARS iterate $k\geq 0$, let $r_k:=||x^{(k)}-x^*||$, the distance from a given STARS iterate to a true minimizer of $f$, denoted $x^* \in \Lambda$. Along the lines of \cite{CW}, we will bound $\mathbb{E}(r_{k+1}^2)-r_k$, `the expected change of $x$ after each iteration in our setting. We note that with this viewpoint, every iterate so far is known; that is, the sequence of vectors $\{x^i\}_{i=0}^k$ are known/fixed up until index $k$, and thus, the sequence $\{r_i\}_{i=0}^k$ of distances are also known/fixed. In particular, both sequences are non-stochastic, meaning that they are constant with respect to any expected values we may apply upon them, in the context of a given step.

First, observe that by using definitions, we may write

\begin{equation} \label{eq:040}
r_{k+1}^2 = ||x^{k+1}-x^*||^2=||x^{(k)} - \hat{h}s_{\mu_k}-x^*||^2.
\end{equation}

Rearranging and expanding,

\begin{equation} \label{eq:041}
r_{k+1}^2=||(x^{(k)}-x^*)-hs_{\mu_k}||^2= r_k^2 -2\hat{h} \langle s_{\mu_k}, x^{(k)} -x^* \rangle+\hat{h}^2||s_{\mu_k}||^2.
\end{equation}

Let $\mathbb{E}$ continue to denote the expectation over $\mathcal{Q}_{k}$ and $\mathcal{U}_{k}$, all of the random vectors and noise draws defining our first $k$ iterates. Recall that one of our current assumptions is that all of the $x^{(k)}$'s (and thus also all of the $r_k$'s) are given/fixed, as well as $\hat{h}$ and $x^*$. Hence, $\mathbb{E}(r_k^2)=r_k^2$, $\mathbb{E}(x^{(k)})=x^{(k)}$, $\mathbb{E}(\hat{h})=\hat{h}$, and $\mathbb{E}(x^*)=x^*$ as we already have these constant objects in hand. However, the next iterate, $x^{k+1}$, \textit{will} depend on the stochastic direction $u^{(k)}$, as well as the stochastic noise values $\xi_{k-1}$ and $\xi_k$ -- and since these stochastic objects literally define $r_{k+1}$ and $s_{\mu_k}$, the expectations will not drop from these terms.

Applying $\mathbb{E}$ to both sides of \eqref{eq:041},

\begin{equation} \label{eq:042}
\mathbb{E}(r_{k+1}^2) = r_k^2-2\hat{h} \langle \mathbb{E}\left(s_{\mu_k}\right),x^{(k)} -x^* \rangle+\hat{h}^2\mathbb{E}\left(||s_{\mu_k}||^2\right).
\end{equation}

We begin by noting the appearance of $\mathbb{E}\left(s_{\mu_k}\right)$, which we can characterize using a pair of previous results. First, we found in \eqref{eq:031} that $\mathbb{E}_{\xi_{k-1},\xi_k}(s_{\mu_k})=g_\mu(x^{(k)})$. Next, we recall \eqref{eq:010}, $\mathbb{E}_u \left( g_\mu (x) \right) = \nabla f_\mu (x) \quad \forall x \in \Lambda$. Putting these results together, we have $\mathbb{E}\left(s_{\mu_k}\right)=\nabla f_\mu (x^{(k)})$.
We can invoke the main result in STARS Lemma 4.4 (Modified) (summarized by \eqref{eq:038}) to bound $\mathbb{E}\left(||s_{\mu_k}||^2\right)$ and our new characterization of $\mathbb{E}\left(s_{\mu_k}\right)=\nabla f_\mu (x^{(k)})$ to  write

\begin{equation} \label{eq:043}
\mathbb{E}(r_{k+1}^2) \leq r_k^2-2\hat{h} \langle \nabla f_\mu (x^{(k)}), x^{(k)} -x^* \rangle+\hat{h}^2\left(2(P+4)||\nabla f(x^{(k)}) ||^2 + C_2\right).
\end{equation}

Now, to reach our next key result, we shall need to verify that $f_\mu$ is convex. Let $x,y \in \Lambda$ and $\mu>0$. By definition, $f_\mu (y) = \mathbb{E}_u(f(y+\mu u))$. Recalling that we have assumed the convexity of $f$ we invoke \eqref{eq:04}, writing

\begin{equation} \label{eq:044}
f_\mu (y) = \mathbb{E}_u(f(y+\mu u))\geq \mathbb{E}_u\left(f(x+\mu u)+\langle \nabla f(x+\mu u),y-x \rangle\right).
\end{equation}

Now using the properties of $\mathbb{E}_u$ and defintions,

\begin{equation} \label{eq:045}
f_\mu (y) \geq \mathbb{E}_u\left(f(x+\mu u)\right)+\langle \nabla \left(\mathbb{E}_u\left( f(x+\mu u)\right)\right),y-x \rangle=f_\mu (x)+\langle \nabla f_\mu (x),y-x\rangle,
\end{equation}

\noindent which holds for any $x,y \in \Lambda$ and $\mu >0$, proving that $f_\mu$ is convex (so long as $f$ is also convex). We shall require a restatement of \eqref{eq:045} for our purposes; we also have for any $x,y\in \Lambda$ and $\mu>0$:

\begin{equation} \label{eq:046}
\langle \nabla f_\mu (x),x-y\rangle \geq  f_\mu (x)- f_\mu (y).
\end{equation}

Hence, plugging $x=x^{(k)}$ and $y=x^*$ into \eqref{eq:046} and multiplying both sides of the inequality by $-2\hat{h}$, we obtain

\begin{equation} \label{eq:047}
-2\hat{h}\langle \nabla f_\mu (x^{(k)}),x^{(k)}-x^*\rangle \leq -2\hat{h}\left(  f_\mu (x^{(k)})- f_\mu (x^*)\right).
\end{equation}

Now \eqref{eq:08} implies $-2\hat{h} f_\mu (x^{(k)}) \leq -2\hat{h}f(x^{(k)})$. Hence,

\begin{equation} \label{eq:048}
-2\hat{h}\langle \nabla f_\mu (x^{(k)}),x^{(k)}-x^*\rangle \leq -2\hat{h}\left(  f(x^{(k)})- f_\mu (x^*)\right).
\end{equation}

Recalling \eqref{eq:03}, we can write $||\nabla f(x^{(k)})||^2 \leq 2L_1 (f(x^{(k)})-f(x^*))$. Using this result along with \eqref{eq:048}, we can update our bound in \eqref{eq:043}, writing

\begin{equation} \label{eq:049}
\mathbb{E}(r_{k+1}^2) \leq r_k^2-2\hat{h}\left(  f(x^{(k)})- f_\mu (x^*)\right)+\hat{h}^2\left(4L_1(P+4)(f(x^{(k)})-f(x^*)) + C_2\right).
\end{equation}

Now we add and subtract $-2\hat{h}f(x^*)$ on the RHS of \eqref{eq:049}, obtaining

\begin{equation} \label{eq:050}
\mathbb{E}(r_{k+1}^2) \leq r_k^2-2\hat{h}\left(  f(x^*)- f_\mu (x^*)\right)-2\hat{h}(f(x^{(k)})-f(x^*))+\hat{h}^2\left(4L_1(P+4)(f(x^{(k)})-f(x^*)) + C_2\right).
\end{equation}

Now \eqref{eq:09} implies $-\frac{\mu^2}{2}L_1P \leq f(x)-f_\mu(x) \implies -2\hat{h}(f(x)-f_\mu(x)) \leq \hat{h}\mu^2 L_1P$ for all $x$, in particular for $x=x^*$; hence,

\begin{equation} \label{eq:051}
\mathbb{E}(r_{k+1}^2) \leq r_k^2+ \hat{h}\mu^2 L_1P -2\hat{h}(f(x^{(k)})-f(x^*))+\hat{h}^2\left(4L_1(P+4)(f(x^{(k)})-f(x^*)) + C_2\right).
\end{equation}

Manipulating and rearranging a little,

\begin{equation} \label{eq:052}
\mathbb{E}(r_{k+1}^2) \leq r_k^2 -2\hat{h}(f(x^{(k)})-f(x^*))(1-2\hat{h}L_1(P+4)) + C_3,
\end{equation}

\noindent where $C_3:=\hat{h}\mu^2 L_1P +\hat{h}^2 C_2.$ Now recall that we have fixed $\hat{h}=(4\hat{L_1}(P+4))^{-1}$; plugging in this particular value, we obtain

\begin{equation} \label{eq:053}
\mathbb{E}(r_{k+1}^2) \leq r_k^2 -\frac{\sqrt{K_1}(2-\sqrt{K_1})}{4L_1(P+4)}(f(x^{(k)})-f(x^*)) + C_3,
\end{equation}

We rewrite $C_3$ (with our $\hat{h}$ and $\hat{\mu}^*$ plugged in) with a function $g(P)$, $P=\dim \Lambda$ given by:

\begin{equation} \label{eq:054}
C_3=\frac{\sqrt{K_1}}{\sqrt{K_2}}\cdot \frac{\sigma}{\sqrt{2}L_1}g(P),\quad \, g(P):=\left( \frac{\sqrt{K_1}}{P+4}\cdot \left(\frac{P}{P+6} \right)^{3/2} +\frac{3K_1+K_2}{16} \cdot \frac{\sqrt{P(P+6)^3}}{(P+4)^2} \right) .
\end{equation}

We may recover the function called $g_1$ from \cite{CW} by setting $K_1=K_2=1$, so that our estimates to $L_1$ and $\sigma^2$ are exact. We shall bound our more general $g$ by bounding each of its terms. Note that the asymptotic analysis of each term gives us hope to find such a bound, as the first term tends to zero as $P \to \infty$ and the second term tends to the constant $(3K_1+K_2)/16$. Using calculus and numerics along the lines of \cite{CW}, we find

\begin{equation} \label{eq:055}
g(P)\leq \sqrt{K_1}\cdot0.036+\frac{3K_1+K_2}{16}\cdot 1.034.
\end{equation}

With $K_1=K_2=1$, we have $g(P) \approx 0.2895 < 3/10$ again matching \cite{CW} (using the bound of $3/10$). We may now define a constant, $C_4$, which bounds $C_3$ over all dimensions $P$ in terms of $L_1$, $\sigma$, $K_1$, and $K_2$:

\begin{equation} \label{eq:056}
C_3\leq \frac{\sqrt{K_1}}{\sqrt{K_2}}\cdot \frac{\sigma}{\sqrt{2}L_1} \left(\sqrt{K_1}\cdot0.036+\frac{3K_1+K_2}{16}\cdot 1.034\right)=:C_4.
\end{equation}

Thus, we can update the bound in \eqref{eq:053} to write

\begin{equation} \label{eq:057}
\mathbb{E}(r_{k+1}^2) \leq r_k^2 -\frac{\sqrt{K_1}(2-\sqrt{K_1})}{4L_1(P+4)}(f(x^{(k)})-f(x^*)) + C_4.
\end{equation}

Applying the expectation over $\mathcal{U}_k$ and $\mathcal{P}_k$,

\begin{equation} \label{eq:058}
\mathbb{E}_{\mathcal{U}_k,\mathcal{P}_k}(r_{k+1}^2) \leq \mathbb{E}_{\mathcal{U}_{k-1},\mathcal{P}_{k-1}}(r_k^2) -\frac{\sqrt{K_1}(2-\sqrt{K_1})}{4L_1(P+4)}(\phi_k-f^*) + C_4.
\end{equation}

Rearranging, we have

\begin{equation} \label{eq:059}
\phi_k-f^* \leq   \frac{4L_1(P+4)}{\sqrt{K_1}(2-\sqrt{K_1})}\left(\mathbb{E}_{\mathcal{U}_{k-1},\mathcal{P}_{k-1}}(r_k^2) - \mathbb{E}_{\mathcal{U}_k,\mathcal{P}_k}(r_{k+1}^2)  + C_4\right).
\end{equation}

Summing over $k=0,\ldots,M$ and dividing by $M+1$, we obtain

\begin{equation} \label{eq:060}
\sum_{k=0}^M\frac{\phi_k-f^*}{M+1} \leq   \frac{4L_1(P+4)}{\sqrt{K_1}(2-\sqrt{K_1})(M+1)}\left(r_0-\mathbb{E}_{\mathcal{U}_k,\mathcal{P}_k}(r_{k+1}^2)\right)+\frac{4L_1(P+4)}{\sqrt{K_1}(2-\sqrt{K_1})}C_4.
\end{equation}

Dropping the strictly negative term on the RHS (involving $\mathbb{E}_{\mathcal{U}_k,\mathcal{P}_k}(r_{k+1}^2)$) and plugging in the definition of $r_0$ (which is a constant, and not stochastic) and noting that we require $0<K_1<4$ and $K_2>0$,
we have found

\begin{equation} \label{eq:061}
\sum_{k=0}^M\frac{\phi_k-f^*}{M+1} \leq   \frac{4L_1(P+4)||x^{(0)}-x^*||^2}{\sqrt{K_1}(2-\sqrt{K_1})(M+1)} + \frac{4\sigma(P+4)}{\sqrt{2K_2}(2-\sqrt{K_1})}C_5,
\end{equation}

\noindent where $C_5:=\sqrt{K_1}\cdot0.036+\frac{3K_1+K_2}{16}\cdot 1.034$, our main result. $\blacksquare$

Again, we recover the exact result of STARS Theorem 4.5 by taking $K_1=K_2=1$.

\textbf{Remark:} We now mimic the analysis in \cite{CW} to explain the implications of our modified Theorem 4.5. First, let $||x^{(0)}-x^*||^2 \leq R^2$. Define $x^\dagger:= \text{argmin}_x\{f(x):x \in \{x^{(0)}, \ldots, x^M\} \}$ and $\phi_\dagger:=\mathbb{E}_{\mathcal{U}_{k-1},\mathcal{P}_{k-1}}(f(x^\dagger))$. Then the value $\phi_\dagger-f^*$ must be less than or equal to the average improvement for any given run of STARS; that is, \eqref{eq:061}, along with our new definitions, implies that

\begin{equation} \label{eq:062}
\phi_\dagger-f^* \leq \sum_{k=0}^M\frac{\phi_k-f^*}{M+1} \leq   \frac{4L_1(P+4)}{\sqrt{K_1}(2-\sqrt{K_1})(M+1)}R^2 + \frac{4\sigma(P+4)}{\sqrt{2K_2}(2-\sqrt{K_1})}C_5.
\end{equation}

Along the lines of \cite{CW}, let us now assume that we wish to achieve a final accuracy of $\epsilon_{\text{tol}}>0$. Then we will need $\phi_\dagger-f^* \leq \epsilon_{\text{tol}}$. If we take 

\begin{equation} \label{eq:063}
 \frac{4\sigma(P+4)}{\sqrt{2K_2}(2-\sqrt{K_1})}C_5 \leq \frac{\epsilon_{\text{tol}}}{2},
\end{equation}

\noindent then we must require that the noise not exceed the following threshold:

\begin{equation} \label{eq:064}
\sigma \leq \frac{\sqrt{2K_2}(2-\sqrt{K_1})\epsilon_{\text{tol}}}{8(P+4)C_4},
\end{equation}

If we satisfy \eqref{eq:064}, then we can achieve $\epsilon_{\text{tol}}$ accuracy as long as

\begin{equation} \label{eq:065}
\frac{4L_1(P+4)R^2}{\sqrt{K_1}(2-\sqrt{K_1})(M+1)} \leq \frac{\epsilon_{\text{tol}}}{2} \iff M \geq \frac{8L_1(P+4)R^2}{\sqrt{K_1}(2-\sqrt{K_1})\epsilon_{\text{tol}}}-1.
\end{equation}

Hence, we achieve $\epsilon_{\text{tol}}$ accuracy as long as the noise is small enough, and $M$ is large enough, with details of those bounds given by \eqref{eq:064} and \eqref{eq:065}. Also, we achieve $\epsilon_{\text{tol}}$ accuracy in

\begin{equation} \label{eq:066}
M \sim \mathcal{O}\left( \frac{L_1 P R^2}{\sqrt{K_1}(2-\sqrt{K_1})\epsilon_{\text{tol}}} \right).
\end{equation}

This analysis also shows that given a particular variance in the noise, $\sigma$, the achievable accuracy can be no better (i.e., less) than the value given below:

\begin{equation} \label{eq:067}
\epsilon_{\text{tol}} \geq \frac{8(P+4)C_5}{\sqrt{2K_2}(2-\sqrt{K_1})} \sigma.
\end{equation}

As usual, we recover the results in \cite{CW} by setting $K_1=K_2=1$.

\subsection{ASTARS Convergence} \label{ss:astars-conv}

We now investigate the convergence of ASTARS. 
We will build upon the theoretical results of the last section, meaning \cite{CW} will be heavily invoked again in this section.
Given the exact $j$-dimensional AS $\A$ of $f$, we shall also need results generally regarding the distance between the minimum of $f$ and the minimum that ASTARS obtains.
We shall also need to discuss the corresponding minimizers.
Recall that we denote the minimizer of $f$ with $x^*$ and we have the stochastic-free minimum of $f$, $f^*=f(x^*)$, as before.
Since ASTARS steps in $\A$ only, given an initial iterate $x^{(0)}$, the minimizer ASTARS is able to attain will be of the form $x^*_\A := P_\A(x^*)+P_\I(x^{(0)})$.
We analogously define the stochastic-free $f^*_\A:=f(x^*_\A)$.
Since the initial iterate will not be changed in $\I$ during ASTARS, the components in $x^*_\A$ are fixed in $\I$ at the given initial iterate, given by the inactive coordinates of $P_\I (x^{(0)})$.
However, we do step towards the true $x^*$ in its coordinates corresponding to $\A$, which is why we also obtain $P_\A(x^*)$ in the definition of $x^*_\A$.
Notice that with our definitions, $x^*-x^*_\A = P_\I (x^{(0)}-x^*)$.
Again -- the difference between $x^*$ and $x^*_\A$ will be in the inactive subspace $\I$, given exactly by the projection of $x^{(0)}$ into $\I$, since ASTARS iterations do not perturb $\I$-coordinates.

ASTARS estimates directional derivatives of $f$ strictly for directions in $\A$ -- iterates are not perturbed in $\I$.
Consequently, the ASTARS gradient oracle can only provide gradient information in the $j$ active directions of $f$.
Hence, the gradients we approximate in ASTARS are denoted $\nabla_\A f(x) \in \A$, which is the gradient of $f$ in $\A$.
Gradients in $\I$ will also be needed for our proofs, and they are defined similarly with $\nabla_\I f(x) \in \I$.
Note the \textit{subspace} gradients $\nabla_\A f(x)$ and $\nabla_\I f(x)$ are still computed in $\Lambda$, but each will fall into their respective subspaces upon computation.
In particular, $\nabla_\A f(x)_i = 0$ for $i=j+1,\ldots, P$ and $\nabla_\I f(x)_i = 0$ for $i=1,\ldots, j$.

We first present a lemma which will be used in both the ASTARS and FAASTARS convergence analyses.
Recall that the vectors $r^{(k)}$ and $\tilde{r}^{(k)}$ have components which are drawn from a $N(0,1)$ distribution.
These vectors are used to form random coefficients in a linear combination in $\A$ (or $\tilde{\A}$) to perform ASTARS steps.
Here, we write $r^{(k)} \sim N(0,I_j)$, a multivariate normal distribution, where $0$ denotes the zero vector in $\Lambda^j$ and $I_j$ is the $j\times j$ identity matrix, so that the covariance is $1$ for every element in $r^{(k)}$ but all elements are independent, with zero covariance between elements.
Analogously, we have $\tilde{r}^{(k)}_p \sim N(0,I_{\tilde{j}})$.
In the first lemma, we show that the random directions for ASTARS steps $u^{(k)}$ and $\tilde{u}^{(k)}$ are also distributed normally with zero mean and unit covariance.

\textbf{ASTARS/FAASTARS Lemma 1:} \textit{Let $\tilde{\A}$ denote a $\tilde{j}$-dimensional AS of $\hat{f}$ and let $\A$ denote the true $j$-dimensional AS of $\hat{f}$. Recall $V_\A:=V_{1:P,1:j}$, where $V$ comes for the eigendecomposition (ED) of the exact sensitivity matrix $W$; as well, recall that $\tilde{V}_{\tilde{\A}}:=\tilde{V}_{1:P,1:\tilde{j}}$, where $\tilde{V}$ comes from the ED of the sensitivity matrix $\tilde{W}$, approximated from samples of $\hat{f}$. Let $r^{(k)}$ denote a random vector such that $r^{(k)} \sim N(0,I_{j})$; likewise, let $\tilde{r}^{(k)}$ denote a random vector such that $\tilde{r}^{(k)} \sim N(0,I_{\tilde{j}})$. Let $u^{(k)}:=V_{\A}r^{(k)}$ and $\tilde{u}^{(k)}:=\tilde{V}_{\tilde{\A}}\tilde{r}^{(k)}$. Then both $u^{(k)}$ and $\tilde{u}^{(k)}$ are normal random vectors; i.e., $u^{(k)} \sim N(0,I_j)$ and $\tilde{u}^{(k)}_p \sim N(0,I_{\tilde{j}})$. Also, $u^{(k)} \in \A$ and $\tilde{u}^{(k)} \in \tilde{\A}$ for all $k$.}

\textit{Proof:} We begin by considering the case in which we have the exact AS of $\hat{f}$, $\A$. We recall that since $W$ is a real $P\times P$ symmetric matrix, its ED is $W=VQV^\top$ where $V$ contains the $P$ eigenvectors of $W$ which are orthonormal in this case, due to the symmetry of $W$, meaning $V$ is a unitary matrix. (Note $Q$ contains the eigenvalues of $W$ along its diagonal in descending order.) Recall that $V_\A:=V_{1:P,1:j}$; hence, $V_\A$ is also a unitary matrix. 

By definition, $u^{(k)}= V_\A r^{(k)}$. Since every component of $r^{(k)}$ is a $N(0,1)$ random variable, $u^{(k)}$ is distributed as $u^{(k)} \sim N(0,(V_\A)^\top (V_\A))$. Now since $V_\A$ is unitary, we know $(V_\A)^\top (V_\A)=I_j$. 
Therefore $u^{(k)} \sim N(0,I_j)$, our desired result for the case of an exact AS.

In the case that we are dealing with an approximated $\tilde{j}$-dimensional AS $\tilde{\A}$, it is still the case that $\tilde{V}_{\tilde{\A}}$ is unitary by construction, and so we can follow the proof above analogously, only replacing $j$ with $\tilde{j}$, and state $\tilde{u}^{(k)}_p \sim N(0,I_{\tilde{j}})$ as well. 

To verify $u^{(k)} \in \A$ and $\tilde{u}^{(k)} \in \tilde{\A}$, recall $u^{(k)}=V_{\A}r^{(k)}$.
Then $u^{(k)}$ is a linear combination of the columns of $V_\A$ with coefficients given by $r^{(k)}$.
The columns of $V_A$ are the eigenvectors $v^i$, $i=1,\ldots,j$ of the associated sensitivity matrix $W$ meaning $u^{(k)}$ a linear combination of the first $j$ eigenvectors of $W$.
Since the span of those $j$ eigenvectors equals $\A$ by definition, $u^{(k)}\in \A$.
(Subspaces are closed under linear combinations of their elements.) 
The argument is analogous for $\tilde{u}^{(k)} \in \tilde{\A}$. $\blacksquare$

We now formulate a series of ASTARS results, where we assume $\A$ is correct and not estimated. 
We now show that using a fixed step size $h_\A$ in \eqref{eq:18}, the active smoothing parameter $\mu_\A^*$ in \eqref{eq:18} is optimal, in the sense that the error in the gradient oracle used in \ref{alg:ASTARS-algo} is minimized. 
This result is a direct corollary of ASTARS Proposition 1 above and STARS Theorem 4.3 (Modified) in the previous section, but with $K_1=K_2=1$, so that we have $L_1$ and $\sigma^2$ exactly to form the active hyperparameters.

\textbf{ASTARS Corollary 2:}
\textit{Let the vectors $u_\A^{(k)}$ denote those drawn using Algorithm \ref{alg:ASTARS-algo}; let $f\in \mathcal{C}^{1,1}(\Lambda)$ and assume $f$ is convex; and assume that the i.i.d. noise draws $\epsilon(\xi)$ are additive, zero mean, with bounded variance $\sigma^2$ for all $\xi$.
By fixing the step size $h_\A$ in \eqref{eq:18}, the active smoothing parameter $\mu_\A^*$ in \eqref{eq:18} minimizes the error between the gradient oracle in Algorithm \ref{alg:ASTARS-algo} and the true directional derivative of $f$ in the direction $u_\A^{(k)}$ in the $j$-dimensional AS $\A$. That is, $\mathcal{E}^\A(\mu)$ in \eqref{eq:013}  (with $u=u_\A^{(k)}$) is minimized by the choice $\mu=\mu^*_\A$. In particular, we have the bound} 

\begin{equation} \label{eq:hot-fix-1}
\mathbb{E}_{u_\A^{(k)},\xi_1,\xi_2} \left( \mathcal{E}^\A \left(\mu^*_\A\right) \right) \leq \sqrt{2}\sigma L_1\sqrt{j(j+6)^3}.
\end{equation}

\textit{Proof:} Replacing $\mathcal{E}$ with $\mathcal{E}^\A$ and taking the expectation over the noise and $u=u_\A$, the proof is identical to the proof of STARS Theorem 4.3 (Modified), until we formulate (as in \eqref{eq:024}) 

\begin{equation} \label{eq:hot-fix-2}
\mathbb{E}_{u_\A,\xi_1,\xi_2}(\mathcal{E}^\A(\mu))\leq \frac{1}{\mu^2}\mathbb{E}_u\left(2\sigma^2||u||^2+\frac{L_1^2}{4}\mu^4|| u||^6\right)
\end{equation}

\noindent and proceed to bound the right hand side.
Applying ASTARS Proposition 1, we have $(u_\A^{(k)})_p \sim N(0,1)$ for $p=1,\ldots,P$.
Taking $u=u_\A^{(k)}$ and noting $u_\A^{(k)} \in \A$ (and $\dim \A=j$), we apply  the bounds on the moments $M_p$ of $||u||^p$ given in \eqref{eq:moment-act-3}.
Using $p=2$ and $p=6$ we have

\begin{equation} \label{eq:hot-fix-3}
\mathbb{E}_{u_\A^{(k)},\xi_1,\xi_2}(\mathcal{E}^\A(\mu)) \leq (2\sigma^2 j)\frac{1}{\mu^2}+\left( \frac{L_1^2(j+6)^3}{4}\right)\mu^2.
\end{equation}

We again observe that the right-hand side of the above inequality is uniformly convex for $\mu>0$ with minimizer $\mu^*_\A:=\left( \frac{8\sigma^2 j}{L_1^2(j+6)^3}\right)^{1/4}$.
This particular choice $\mu=\hat{\mu^*_\A}$ can be plugged into \eqref{eq:hot-fix-3}, and we obtain the bound in \eqref{eq:hot-fix-1}, our main result.  $\blacksquare$

Next, define $\mathcal{P}_k :=\{\xi_k\}_{k=1}^{M}$ and  $\mathcal{U}_k^\A :=\{u^{(k)}_\A \}_{k=1}^{M}$, which are two sets containing all random variables that appear in ASTARS, iterations $k=1,\ldots,M$. 
Let $\phi_0:= f(x^{(0)})$ and $\phi_k^\A := \mathbb{E}_{\mathcal{Q}_{k-1},\mathcal{U}^\A_{k-1}}(f(x^{(k)}))$, $k \geq 1,$ where the $x^{(k)}$'s are now ASTARS iterates.
$\mathbb{E}$ will now denote the expectation over every noise draw and random vector used in STARS up to (and including) the $k$-th iterate; that is, the expectations are now taken with respect to $\xi_0,\ldots,\xi_{k}$ and $u^{(1)}_\A,\ldots,u^{(k)}_\A$ unless stated otherwise.

Now, given that the active smoothing parameter $\hat{\mu^*_\A}$ is optimal, in the sense of minimizing $\mathcal{E}^\A$, we present the following result, showing the convergence of ASTARS.
The following result is a direct corollary of STARS Lemma 4.4 (Modified), STARS Theorem 4.5 (Modified), ASTARS/FAASTARS Lemma 1, and FAASTARS Corollary 2.

\textbf{ASTARS Corollary 3:} 
\textit{Let random vectors $u^{(k)}_\A$ be drawn according to \ref{alg:ASTARS-algo}; $f\in \mathcal{C}^{1,1}(\Lambda)$ and $f$ is convex; and that the i.i.d. noise draws $\epsilon(\xi)$ are additive, zero mean, with bounded variance $\sigma^2$ for all $\xi$.
Let $\{x^{(k)}\}_{k\geq 0}$ denote a sequence of ASTARS iterates formed using a fixed active step length $h_\A$ and fixed active smoothing $\mu = \mu^*_\A$ (both given in \eqref{eq:18}) for all ASTARS iterates $k$. Then for any total number of ASTARS iterations $M$,}

\begin{equation} \label{eq:hot-fix-4}
\sum_{k=0}^M\frac{\phi_k^\A-f^*_\A}{M+1} \leq   \frac{4L_1(j+4)||P_\A(x^{(0)}-x^*)||^2}{(M+1)} + \frac{3\sqrt{2}\sigma(j+4)}{5}.
\end{equation}

\textit{Proof:} The proof is almost identical to the proofs of STARS Lemma 4.4 (Modified) and STARS Theorem 4.5 (Modified), but with $j$'s replacing the roles of $P$'s (since we take steps with $j$-dimensional $u^{(k)}_\A$ vectors and not $u^{(k)} \in \Lambda$), $\mathcal{E}^\A$ replacing $\mathcal{E}$, $\mathcal{U}^\A$ replacing $\mathcal{U}$, $\mu^*_\A$ replacing $\hat{\mu}^*$, $\phi_k^\A$ replacing $\phi_k$, $x^*_\A$ replacing $x^*$ (since we are converging to the minimum of $f_\A$, and $K_1=K_2=1$, since we are assuming exact active hyperparameters, formed with the true values for $L_1$ and $\sigma^2$.
Note $||\cdot||=||\cdot||_\Lambda$, a norm on $\Lambda$ throughout.
We outline the required changes one must make to the proofs of STARS Lemma 4.4 (Modified) and STARS Theorem 4.5 (Modified) to obtain our desired result.

We begin by obtaining a bound analogous to that of STARS Lemma 4.4 (Modified), but note that ASTARS steps are taken with random vectors $u^{(k)}_\A \in \A$, and $\dim \A = j$.
First, we replace $u^{(k)}$ in \eqref{eq:029} with $u^{(k)}_\A$, so $g_0(x^{(k)}):=\langle \nabla f(x^{(k)}),u^{(k)}_\A \rangle u^{(k)}_\A$.
Similarly, note that we set $u^{(k)}=u^{(k)}_\A$ in \eqref{eq:012}
Then, using ASTARS Corollary 2 -- instead of the STARS Theorem 4.3 (Modified) -- we have

\begin{equation} \label{eq:hot-fix-5}
\mathbb{E}\left(||s_{\mu^*_\A}^\A||^2 - 2\langle s_{\mu^*_\A}^\A,g_0(x^{(k)}) \rangle + ||g_0(x^{(k)})||^2\right)\leq \sqrt{2}\sigma L_1\sqrt{j(j+6)^3}.
\end{equation}

Using \eqref{eq:active-011}, we have

\begin{equation} \label{eq:hot-fix-6}
\mathbb{E}\left(||s_{\mu^*_\A}^\A||^2\right)  \leq 2(j+4)||\nabla f(x^{(k)}) ||^2 + \frac{(\mu^*_\A)^2 L_1^2}{2}(j+6)^3 + C_1,
\end{equation}

\noindent where we recall we have $C_1=\sqrt{2}\sigma L_1\sqrt{j(j+6)^3}$ here. Plugging in the value of $\mu^*_\A$, we obtain the bound

\begin{equation} \label{eq:hot-fix-7}
\mathbb{E}\left(||s_{\mu^*_\A}^\A||^2\right)  \leq 2(j+4)||\nabla f(x^{(k)}) ||^2 + C_2,
\end{equation}

\noindent where $C_2=2\sqrt{2}L_1\sigma \sqrt{j(j+6)^3}$. The bound in \eqref{eq:hot-fix-6} is the analogous result to STARS Lemma 4.4 (Modified) in the case of ASTARS performed in the known and exact $\A$ and with $K_1=K_2=1$ (exact hyperparameters).

We now proceed to proving the analogous result to STARS Theorem 4.5 (Modified) in our case.
We redefine $r_k:=||x^{(k)}-x^*_\A||$ for ASTARS iterates $x^{(k)}$.
The first three equations appearing in the proof of STARS Theorem 4.5 (Modified), \eqref{eq:040} through \eqref{eq:042}, are nearly identical for this proof -- one must replace $x^*$ with $x^*_\A$,  replace $\hat{h}$ with $h_{\A}$, and note that here we have $s_{\mu^{(k)}}=s_{\mu^*_\A}^\A$. Then, invoking \eqref{eq:hot-fix-7}, we rewrite \eqref{eq:043} in our case as

\begin{equation} \label{eq:hot-fix-8}
\mathbb{E}(r_{k+1}^2) \leq r_k^2-2h_\A \langle \nabla f_\mu^\A (x^{(k)}), x^{(k)} -x^*_\A \rangle+h_\A^2\left(2(j+4)||\nabla f(x^{(k)}) ||^2 + C_2\right),
\end{equation}

\noindent where we recall $C_2=2\sqrt{2}L_1\sigma \sqrt{j(j+6)^3}$ here. Now, again, equations \eqref{eq:044} through \eqref{eq:048} are identical for this proof as long as $x^*_\A$ replaces $x^*$ and $h_\A$ replaces $\hat{h}$ throughout, and similarly for equations \eqref{eq:049} through \eqref{eq:051} with the additional needed replacement of $j$ for $P$ and using $C_2$ as we have defined it here.
Taking $C_3=h_\A \mu^*_\A L_1 j + h_\A^2 C_2$, \eqref{eq:052} holds (with the usual replacments) and plugging $h_\A$ into the modified \eqref{eq:052} gives the following modification to \eqref{eq:053}:

\begin{equation} \label{eq:hot-fix-9}
\mathbb{E}(r_{k+1}^2) \leq r_k^2 -\frac{1}{4L_1(j+4)}(f(x^{(k)})-f(x^*_\A)) + C_3.
\end{equation}

Now plugging $h_\A$ into $C_3$, we obtain $C_3 \leq \frac{3\sigma}{10\sqrt{2}L_1}=:C_4$, so

\begin{equation} \label{eq:hot-fix-10}
\mathbb{E}(r_{k+1}^2) \leq r_k^2 -\frac{1}{4L_1(j+4)}(f(x^{(k)})-f(x^*_\A)) + C_4.
\end{equation}

We may now apply the expectation over $\mathcal{P}_k$ and $\mathcal{U}_k^\A$, rearrange, and sum over $k=0,\ldots, M$ as before. We obtain a modification to \eqref{eq:060}
with

\begin{equation} \label{eq:hot-fix-11}
\sum_{k=0}^M\frac{\phi_k^\A -f^*_\A}{M+1} \leq   \frac{4L_1(j+4)}{(M+1)}\left(r_0-\mathbb{E}_{\mathcal{U}_k^\A,\mathcal{P}_k}(r_{k+1}^2)\right)+4L_1(j+4)C_4.
\end{equation}

We drop the strictly negative term (again involving $\mathbb{E}_{\mathcal{P}_k,\mathcal{U}_k^\A}$) and plug in the definition of $r_0=||x^{(0)}-x^*_\A||^2$.
Then, recalling $x^*_\A=P_\A(x^*)+P_\I(x^{(0)})$, writing $x^{(0)}=P_\A(x^{(0)})+P_\I(x^{(0)})$, and noting that $P_\A$ is linear, we have $r_0=||P_\A(x^{(0)}-x^*)||^2$.
We obtain \eqref{eq:hot-fix-4}. $\blacksquare$

We have shown the convergence of ASTARS to the minimum $f^*_\A$ of $f_\A$ with correct hyperparameters and the correct and known $\A$.
Ultimately, to obtain complexity results for ASTARS, we desire a statement about the convergence of ASTARS to $f^*$, the minimum of $f$.
We pay a price for stepping only in active variables in ASTARS, which is that inactive variables are not minimized or even perturbed at all. 
Because ASTARS converges to $f_\A^*$, we will not minimize $f$ in its inactive variables $\I$, and $|f^*-f_\A^*|$ may be nonzero.
By modifying results in \cite{ConstantineK}, we show that this difference will usually be negligible, as it is bounded by the square root of the sum of the eigenvalues associated with $\I$ (which are usually small), scaled by a constant related to the distance from $x^*$ to $x^*_\A$ (which is also small in many cases).

\textbf{ASTARS Corollary 4:} \textit{Let $x^{(0)}$ denote any initial iterate for ASTARS. Let $x^*$ denote a true minimizer of $f$ with the corresponding stochastic-free function evaluation given by $f^*$. Let $\A$ continue to denote the true $j$-dimensional AS of $f$. Denote the ASTARS minimizer with $x^*_\A$ and corresponding noiseless function evaluation $f^*_\A$.
Assume $||x^*-x^*_\A||_\Lambda^2 < \infty$, where $||\cdot||_\Lambda$ denotes a norm.
Then we may bound the difference between $|f^*-f^*_\A|$ with}

\begin{equation} \label{eq:hot-fix-12}
|f^*-f^*_\A| \leq \sqrt{a_1 (q_{j+1}+ \cdots + q_P)},
\end{equation}

\noindent \textit{where $0 \leq a_1<\infty$ is an eigenvalue of the positive-semi definite matrix $(x^*-x^*_\A)(x^*-x^*_\A)^\top$ and the $q$'s are our notations for the exact eigenvalues associated with the eigendecomposition of the sensitivity matrix $W$. Also, $a_1 = ||P_\I(x^{(0)}-x^*)||_\Lambda^2.$}

\textit{Proof:} We bound the quantity $(f^*-f^*_\A)^2$.
First, we expand $f^*=f(x^*)$ around $x^*_\A$ using a special case of Taylor's theorem, sometimes called the \textit{Extended Mean Value Theorem}.
For $c \in [0,1]$ and $z := c x_\A^* + (1-c)x^*$ we have $f^*=f(x^*)=f(x_{\A}^*)+ \nabla f (z)^\top (x^* - x^*_\A)  $.
Note that since the components of $x_\A^*$ and $x^*$ must match for indices $i=1,\ldots,j$ (by definition), the point $z$ varies along $\I$ only; that is, $z$ is fixed in $\A$.
Thus, $\nabla f (z) = \nabla_\I f(z)$, the gradient taken in the inactive subspace only.
The expansion of $f^*$ around $x^*_\A$ allows us to write $(f^*-f^*_\A)^2=\left(\nabla_\I f (z)^\top (x^* - x^*_\A)\right)^2= \nabla_\I f (z)^\top A \nabla_\I f (z)$, where $A:=(x^* - x^*_\A)(x^* - x^*_\A)^\top$ is a $P \times P$ matrix.
Note that $A$ is a square, positive semi-definite, rank 1 matrix.
Hence, it has 1 eigenvalue that is positive or zero, which we denote with $a_1 \geq 0$; all other $P-1$ eigenvalues are $0$.
We find $a_1 = ||x^*-x^*_\A||^2_\Lambda$ and so by definition, $a_1=||P_\I(x^{(0)}-x^*)||_\Lambda^2$.
Observe $a_1<\infty$ because $||x^*-x^*_\A||_\Lambda < \infty$ and also $a_1=0 \iff x^*=x^*_\A$, in which case $f^*=f^*_\A$ so that $|f^*-f^*_\A|=0$.
We have $(f^*-f^*_\A)^2 \leq a_1 \nabla_\I f (z)^\top \nabla_\I f (z)$.
Applying the expectation over $\I$ to both sides (where the left-hand side is constant with respect to this expectation) we have $(f^*-f^*_\A) \leq a_1 \mathbb{E}_{\I}\left(\nabla_\I f (z)^\top \nabla_\I f (z)\right)$.
Citing \cite{ConstantineK} Lemma 2.2, we have $\mathbb{E}_{\I}\left(\nabla_\I f (z)^\top \nabla_\I f (z)\right)<q_{j+1}+\cdots + q_P$, where $q_i$, $i=j+1,\ldots,P$ are the last $P-j$ eigenvalues of the sensitivity matrix $W$.
Hence, $(f^*-f^*_\A)^2 \leq a_1 (q_{j+1}+\cdots + q_P)$.
Applying a square root to both sides, we obtain \eqref{eq:hot-fix-12}. $\blacksquare$

We now present a statement about the convergence of ASTARS to $f^*$ by combining the previous two results.
Recall $\phi_0:= f(x^{(0)})$ and $\phi_k^\A := \mathbb{E}_{\mathcal{Q}_{k-1},\mathcal{U}_{k-1}^\A}(f(x^{(k)}))$, $k \geq 1,$ where the $x^{(k)}$'s are ASTARS iterates.

\textbf{ASTARS Theorem 5:}
\textit{Let random vectors $u^{(k)}_\A$ be drawn according to \ref{alg:ASTARS-algo}; $f\in \mathcal{C}^{1,1}(\Lambda)$ and $f$ is convex; and that the i.i.d. noise draws $\epsilon(\xi)$ are additive, zero mean, with bounded variance $\sigma^2$ for all $\xi$.
Let $\{x^{(k)}\}_{k\geq 0}$ denote a sequence of ASTARS iterates formed using a fixed active step length $h_\A$ and fixed active smoothing $\mu = \mu^*_\A$ (both given in \eqref{eq:18}) for all ASTARS iterates $k$. For any total number of ASTARS iterations $M\geq 1$,}

\begin{equation} \label{eq:hot-fix-13}
\sum_{k=0}^M\frac{\phi_k^\A-f^*}{M+1} \leq \frac{4L_1(j+4)||P_\A(x^{(0)}-x^*)||^2 }{(M+1)} + \frac{3\sigma(j+4)}{5\sqrt{2}} + \sqrt{a_1 (q_{j+1}+ \cdots + q_P)}
\end{equation}

\textit{Proof:} By the triangle inequality, we have $|\phi_k-f^*| \leq |\phi_k - f_\A^*| + |f^*-f_\A^*|$.
Noting that $\phi_k - f_\A^*>0$ for all $k$, we can bound the left-hand side of \eqref{eq:hot-fix-11}, writing 

\begin{equation} \label{eq:hot-fix-14}
\sum_{k=0}^M\frac{\phi_k-f^*}{M+1} \leq \sum_{k=0}^M\frac{\phi_k-f^*_\A}{M+1} + |f^*-f_\A^*|.
\end{equation}

Now the first term on the right-hand side of \eqref{eq:hot-fix-14} is bounded by ASTARS Corollary 3 and the second term is bounded by ASTARS Corollary 4.
Plugging in those bounds, we obtain \eqref{eq:hot-fix-13}. $\blacksquare$

We use the results above to analyze the complexity of ASTARS in the following remark.

\textbf{Remark:} We now mimic the complexity analysis we performed in the preceding section for STARS with approximated hyperparameters for our case in this section, ASTARS with correct hyperparameters and $\A$.
Let $||\cdot||$ denote a norm in $\Lambda$ throughout.
Define $R_\A^2$ as a bound $||P_\A(x^{(0)}-x^*)||^2 \leq R_\A^2$.
We recall $a_1 = ||P_\I(x^{(0)}-x^*)||_\Lambda^2$ and define a bound $a_1 \leq R_\I^2$.

Now define $x^\dagger:= \text{argmin}_x\{f(x):x \in \{x^{(0)}, \ldots, x^{(M)}\} \}$ and $\phi_\dagger:=\mathbb{E}_{\mathcal{U}_{k-1},\mathcal{P}_{k-1}}(f(x^\dagger))$. Then the value $\phi_\dagger-f^*$ must be less than or equal to the average improvement for any given run of STARS; that is, \eqref{eq:hot-fix-13}, along with our new definitions, implies that

\begin{equation} \label{eq:hot-fix-15}
\phi_\dagger-f^* \leq \sum_{k=0}^M\frac{\phi_k-f^*}{M+1} \leq   \frac{4L_1(j+4)}{(M+1)}R_\A^2 + \frac{3\sigma(j+4)}{5\sqrt{2}} + R_\I \sqrt{(q_{j+1}+ \cdots + q_P)}.
\end{equation}

We assume that we wish to achieve a final accuracy of $\epsilon_{\text{tol}}>0$. Then we will need $\phi_\dagger-f^* \leq \epsilon_{\text{tol}}$. If we take 

\begin{equation} \label{eq:hot-fix-16}
\frac{3\sigma(j+4)}{5\sqrt{2}} \leq \frac{\epsilon_{\text{tol}}}{3},
\end{equation}

\noindent then we must require that the noise not exceed the following threshold:

\begin{equation} \label{eq:hot-fix-17}
\sigma \leq \frac{5\sqrt{2}\epsilon_{\text{tol}}}{9\sigma(j+4)}.
\end{equation}

\noindent If we satisfy \eqref{eq:hot-fix-17} and also have

\begin{equation} \label{eq:hot-fix-18}
R_\I \sqrt{ (q_{j+1}+ \cdots + q_P)}\leq \frac{\epsilon_{\text{tol}}}{3},
\end{equation}

\noindent then we can achieve $\epsilon_{\text{tol}}$ accuracy as long as 

\begin{equation} \label{eq:hot-fix-19}
\frac{4L_1(j+4)}{(M+1)}R_\A^2 \leq \frac{\epsilon_{\text{tol}}}{3} \iff M \geq \frac{12 L_1(j+4)R_\A^2}{\epsilon_{\text{tol}}}-1.
\end{equation}

Hence, we achieve $\epsilon_{\text{tol}}$ accuracy as long as: the noise is small enough; the eigenvalues of the inactive subspace and distance from $x^*$ to $x^*_\A$ is small enough; and $M$ is large enough.
Details of those required bounds given by \eqref{eq:hot-fix-17}, \eqref{eq:hot-fix-18}, and \eqref{eq:hot-fix-19}, respectively. With these assumptions, we achieve $\epsilon_{\text{tol}}$ accuracy in

\begin{equation} \label{eq:hot-fix-20}
M \sim \mathcal{O}\left( \frac{L_1 j R_\A^2}{\epsilon_{\text{tol}}} \right).
\end{equation}

This analysis also shows that given a particular variance in the noise, $\sigma$, as well as the term 
involving the eigenvalues of the inactive subspace and distance from $x^*$ to $x^*_\A$, the achievable accuracy can be no better (i.e., less) than the value given below:

\begin{equation} \label{eq:hot-fix-21}
\epsilon_{\text{tol}} \geq \max \left\{ \frac{9(j+4)}{5\sqrt{2}} \sigma, 3 R_\I \sqrt{(q_{j+1}+\cdots+q_p)}  \right\}.
\end{equation}


\subsection{FAASTARS Convergence} \label{ss:faastars-conv}

Now we focus on analyzing the convergence of FAASTARS.
Here, we must consider that FAASTARS uses approximate information both for hyperparameters and for the AS in its phases.
We have already analyzed the convergence of performing STARS with estimated hyperparameters in \ref{ss:stars-conv}, which corresponds to the first phase of FAASTARS.
Before we can state our main result about the convergence of FAASTARS, we will need results analogous to those in \ref{ss:astars-conv}, but with estimated hyperparameters \textit{and} and estimated AS $\tilde{A}$.
We first reintroduce the approximately-optimal ASTARS hyperparameters, need to perform the third and final phase of FAASTARS:

\begin{equation} \label{eq:068}
\hat{\mu}^*_{\tilde{\A}}:=\left( \frac{8\hat{\sigma}^2 \tilde{j}}{\hat{L_1}^2(\tilde{j}+6)^3}\right)^{1/4} \quad \quad \hat{h}_{\tilde{\A}}:=(4\hat{L_1}(\tilde{j}+4))^{-1}.
\end{equation}

We also must modify our definitions from the previous section to account for the approximated subspace.
Recall that the third phase of FAASTARS will begin with the initial iterate $x^{(M_\A)}$ (the last iterate from phase two).
Define $x^*_{\tilde{\A}} := P_{\tilde{\A}}(x^*)+P_{\tilde{\I}}(x^{(M_\A)})$ with its associated stochastic-free $f^*_\A:=f(x^*_\A)$.
Observe that with our definitions, $x^*-x^*_{\tilde{\A}} = P_{\tilde{\I}} (x^{(M_\A)})$.
We define $f|_{\tilde{\A}}(\lambda):=f(V_{\tilde{\A}} V_{\tilde{\A}}^\top \lambda)=f(P_{\tilde{\A}}(\lambda))$ and let $f_{\tilde{\A}}:=f|_{\tilde{\A}}$.
Note that $f_{\tilde{\A}}$ is convex since $f$ is convex and we can define $f_{\tilde{\I}}$ analogously.
Also, when we evaluate gradients for points $x \in \tilde{\A}$, we note we obtain the object $\nabla_{\tilde{\A}} f(x) \in \tilde{\A}$, and similarly for $\tilde{\I}$.

We begin by providing a modification to ASTARS Corollary 2 for the case of estimated hyperparameters and $\tilde{\A}$.
The proof is a blend of the proofs of STARS Theorem 4.3 (Modified) and ASTARS Corollary 2, with additional consideration for the now $\tilde{j}$-dimensional $\tilde{\A}$.

\textbf{FAASTARS Corollary 1 (Modified ASTARS Corollary 2):}
\textit{Let the vectors $u_{\tilde{\A}}^{(k)}$ denote those drawn using Algorithm \ref{alg:FAASTARS-algo-3}; let $f\in \mathcal{C}^{1,1}(\Lambda)$ and assume $f$ is convex; and assume that the i.i.d. noise draws $\epsilon(\xi)$ are additive, zero mean, with bounded variance $\sigma^2$ for all $\xi$.
We assume we have fixed estimates $\hat{L_1}$ and $\hat{\sigma}$ with $K_1>0$ and $K_2>0$ as in \eqref{eq:014}.
By fixing the step size as $\hat{h}_{\tilde{\A}}$ in \eqref{eq:068}, the approximately active smoothing parameter $\hat{\mu}_{\tilde{\A}}^*$ in \eqref{eq:068} minimizes the error between the gradient oracle in Algorithm \ref{alg:FAASTARS-algo-3} and the true directional derivative of $f$ in the direction $u_{\tilde{\A}}^{(k)}$ in the $\tilde{j}$-dimensional AS $\tilde{\A}$. That is, $\mathcal{E}(\mu)$ in \eqref{eq:013}  (with $u=u_{\tilde{\A}}^{(k)}$) is minimized by the choice $\mu=\hat{\mu}_{\tilde{\A}}^*$. In particular, we have the bound} 

\begin{equation} \label{eq:faastars-1}
\mathbb{E}_{u_{\tilde{\A}}^{(k)},\xi_1,\xi_2} \left( \mathcal{E}^{\tilde{\A}} \left(\mu^*_{\tilde{\A}} \right) \right) \leq \frac{K_1+K_2}{\sqrt{2K_1 K_2}}\sigma L_1\sqrt{\tilde{j}(\tilde{j}+6)^3}.
\end{equation}

\textit{Proof:} The proof is identical to the proof of STARS Theorem 4.3 (Modified) and ASTARS Corollary 2 (making the usual substitutions), until we formulate (as in \eqref{eq:024}) 

\begin{equation} \label{eq:faastars-2}
\mathbb{E}_{u,\xi_1,\xi_2}(\mathcal{E}^{\tilde{\A}}(\mu))\leq \frac{1}{\mu^2}\mathbb{E}_u\left(2\sigma^2||u||^2+\frac{L_1^2}{4}\mu^4|| u||^6\right)
\end{equation}

\noindent and proceed to bound the right hand side.
Applying ASTARS Proposition 1, we have $(u_{\tilde{\A}}^{(k)})_p \sim N(0,1)$ for $p=1,\ldots,P$.
Taking $u=u_{\tilde{\A}}^{(k)}$ and noting $u_{\tilde{\A}}^{(k)} \in \tilde{\A}$ (and $\dim \tilde{\A}=\tilde{j}$), we apply  the bounds on the moments $M_p$ of $||u||^p$ given in \eqref{eq:moment-act-3}.
Using $p=2$ and $p=6$ -- but replacing the AS dimension $j$ with $\tilde{j}$ -- we have

\begin{equation} \label{eq:faastars-3}
\mathbb{E}_{u_{\tilde{\A}}^{(k)},\xi_1,\xi_2}(\mathcal{E}(\mu)) \leq (2\sigma^2 \tilde{j})\frac{1}{\mu^2}+\left( \frac{L_1^2(\tilde{j}+6)^3}{4}\right)\mu^2.
\end{equation}

We again observe that the right-hand side of the above inequality is uniformly convex for $\mu>0$ with minimizer $\mu^*_{\tilde{\A}}:=\left( \frac{8\sigma^2 \tilde{j}}{L_1^2(\tilde{j}+6)^3}\right)^{1/4}$.
Analogously to the proof of STARS Theorem 4.3 (Modified), our optimal choice of smoothing, given the information we have available, will require us to swap out $L_1$ and $\sigma^2$ in $\mu^*_\A$ with their estimates, $\hat{L_1}$ and $\hat{\sigma}^2$, and to swap out $j$ for $\tilde{j}$, recovering $\hat{\mu}^*_{\tilde{\A}}$ \eqref{eq:068}. This particular choice $\mu=\hat{\mu}^*_{\tilde{\A}}$ can be plugged into \eqref{eq:faastars-3}, which gives us the bound in \eqref{eq:faastars-1}, our main result.  $\blacksquare$

Next, we redefine $\mathcal{Q}_k :=\{\xi_k\}_{k=1}^{M_\A}$ and  $\mathcal{U}_k :=\{u^{(k)}\}_{k=1}^{M_\A}$, which are two sets containing all random variables that appear in FAASTARS' regular STARS burn-in phase, iterations $k=1,\ldots,M_\A$. Likewise, we extend the definitions of each set so that $\mathcal{Q}_k =\{\xi_k\}_{k=M_\A+1}^{M}$ and  $\mathcal{U}_k^{\tilde{\A}} =\{u^{(k)}_{\tilde{\A}}\}_{k=M_\A+1}^{M}$ also contain all random variables that appear in FAASTARS' approximate ASTARS phase, iterations $k=M_\A+1,\ldots,M$. 

Let $\phi_0:= f(x^{(0)})$, $\phi_k := \mathbb{E}_{\mathcal{Q}_{k-1},\mathcal{U}_{k-1}}(f(x^{(k)}))$, $1 \leq k \leq  M_\A,$ and $\phi_k := \mathbb{E}_{\mathcal{Q}_{k-1},\mathcal{U}_{k-1}^{\tilde{\A}}}(f(x^{(k)}))$, $M_\A+ 1 \leq k \leq M$ where the $x^{(k)}$'s are now STARS iterates for $1 \leq k \leq  M_\A,$ and FAASTARS iterates for $M_\A+ 1 \leq k \leq M$.

\textbf{FAASTARS Corollary 2 (FAASTARS, approximate ASTARS phase result):}
\textit{For all FAASTARS iterates in phase 3, $k=M_\A+1,\ldots,M$, let the vectors $u_{\tilde{\A}}^{(k)}$ denote those drawn using Algorithm \ref{alg:FAASTARS-algo-3}.
Also, let $f\in \mathcal{C}^{1,1}(\Lambda)$ with $f$ convex, and let i.i.d. noise draws $\epsilon(\xi)$ be additive, zero mean, with bounded variance $\sigma^2$ for all appearing $\xi$.
We assume we have fixed estimates $\hat{L_1}$ and $\hat{\sigma}$ with $K_1>0$ and $K_2>0$ as in \eqref{eq:014}.
Let $\tilde{\A}$ denote the approximated $\tilde{j}$-dimensional AS of $f$.
Let $x^{(M_\A)}$ be fixed and given and let $\{x^{(k)}\}_{k=M_\A +1}^M$ denote a sequence of FAASTARS iterates formed using the approximate active hyperparameters in \eqref{eq:068}). 
\textbf{Finally, we require $0<K_1<4$ and $K_2>0$, the values defined in \eqref{eq:014}.} Then for any $M-M_\A$ total number of approximate ASTARS iterations within FAASTARS, $k=M_\A+1,\ldots , M$,}

\begin{equation} \label{eq:070}
\sum_{k=M_\A}^M\frac{\phi_k-f^*_{\tilde{\A}}}{M-M_\A} \leq   \frac{4L_1(\tilde{j}+4)||P_{\tilde{\A}}(x^{(M_\A)}-x^*)||^2}{\sqrt{K_1}(2-\sqrt{K_1})(M-M_\A+1)} + \frac{4\sigma(\tilde{j}+4)}{\sqrt{2K_2}(2-\sqrt{K_1})}C_5,
\end{equation}

\noindent where $C_5:=\sqrt{K_1}\cdot0.036+\frac{3K_1+K_2}{16}\cdot 1.034$.

\textit{Proof:} The proof is almost identical to the proof of ASTARS Corollary 3, but with $\tilde{j}$'s replacing the roles of $j$'s (since we take steps with $\tilde{j}$-dimensional $u^{(k)}_{\tilde{\A}}$ vectors and not $u^{(k)}_\A \in \A$), $\hat{\mu}^*_{\tilde{\A}}$ replacing $\mu^*_\A$, and $x^*_{\tilde{\A}}$ replacing $x^*_\A$ (since we are converging to the minimum of $f_{\tilde{\A}}$). 
Here, $K_1$ and $K_2$ are not necessarily equal to $1$ as before, since we are assuming inexact active hyperparameters, formed with the estimates $\hat{L_1}$ and $\hat{\sigma}^2$.
We outline the required changes one must make to the proofs of STARS Lemma 4.4 (Modified) and STARS Theorem 4.5 (Modified) to obtain our desired result.
These changes essentially amount to keeping the logic from STARS Lemma 4.4 (Modified) and STARS Theorem 4.5 (Modified) to account for $K_1$ and $K_2$ but to replacing $j$ with $\tilde{j}$.

We begin by obtaining a bound analogous to that of STARS Lemma 4.4 (Modified), but note that the approximate ASTARS steps are taken with random vectors $u^{(k)}_{\tilde{\A}} \in \tilde{\A}$, and $\dim \tilde{\A} = \tilde{j}$.
First, we replace $u^{(k)}$ in \eqref{eq:029} with $u^{(k)}_{\tilde{\A}}$, so $g_0(x^{(k)}):=\langle \nabla f(x^{(k)}),u^{(k)}_{\tilde{\A}} \rangle u^{(k)}_{\tilde{\A}}$.
Similarly, note that we set $u^{(k)}=u^{(k)}_{\tilde{\A}}$ in \eqref{eq:012}.
Also, let $s_{\hat{\mu}^*_{\tilde{\A}}}=s_{\hat{\mu}^*_{\tilde{\A}}}^{\tilde{\A}}$ for cleaner notation.
Then, using FAASTARS Corollary 1 -- instead of ASTARS Corollary 2 -- \eqref{eq:030} becomes

\begin{equation} \label{eq:faastars-4}
\mathbb{E}\left(||s_{\hat{\mu}^*_{\tilde{\A}}}||^2 - 2\langle s_{\hat{\mu}^*_{\tilde{\A}}},g_0(x^{(k)}) \rangle + ||g_0(x^{(k)})||^2\right)\leq \frac{K_1+K_2}{\sqrt{2K_1 K_2}}\sigma L_1\sqrt{\tilde{j}(\tilde{j}+6)^3}.
\end{equation}

Continuing with $u^{(k)}_{\tilde{\A}}$ replacing $u^{(k)}$, we proceed identically, noting $C_1=\frac{K_1+K_2}{\sqrt{2K_1 K_2}}\sigma L_1\sqrt{\tilde{j}(\tilde{j}+6)^3}$.
Modifying \eqref{eq:active-011} for $x^{(k)} \in \tilde{\A}$ with $\dim \tilde{\A} =\tilde{j}$, we have

\begin{equation} \label{eq:faastars-5}
\mathbb{E}\left(||s_{\hat{\mu}^*_{\tilde{\A}}}||^2\right)  \leq 2(\tilde{j}+4)||\nabla f(x^{(k)}) ||^2 + \frac{(\hat{\mu}^*_{\tilde{\A}})^2 L_1^2}{2}(\tilde{j}+6)^3 + C_1,
\end{equation}

Plugging in the value of $\hat{\mu}^*_{\tilde{\A}}$, we obtain the bound

\begin{equation} \label{eq:faastars-6}
\mathbb{E}\left(||s_{\hat{\mu}^*_{\tilde{\A}}}||^2\right)  \leq 2(\tilde{j}+4)||\nabla f(x^{(k)}) ||^2 + C_2,
\end{equation}

\noindent where $C_2:=\frac{3K_1+K_2}{\sqrt{2K_1K_2}}L_1\sigma\sqrt{\tilde{j}(\tilde{j}+6)^3}=\frac{3K_1+K_2}{\sqrt{2}}  \hat{L_1}\hat{\sigma}\sqrt{\tilde{j}(\tilde{j}+6)^3}$. The bound in \eqref{eq:faastars-6} is the analogous result to STARS Lemma 4.4 (Modified) in the case of ASTARS performed in the estimated $\tilde{\A}$ and with inexact hyperparameters.

We now proceed to proving the analogous result to STARS Theorem 4.5 (Modified) in our case.
We redefine $r_k:=||x^{(k)}-x^*_{\tilde{\A}}||$ for approximate ASTARS iterates $x^{(k)}$.
The first three equations appearing in the proof of STARS Theorem 4.5 (Modified), \eqref{eq:040} through \eqref{eq:042}, are nearly identical for this proof -- one must replace $x^*$ with $x^*_{\tilde{\A}}$,  replace $\hat{h}$ with $\hat{h}_{\tilde{\A}}$, and note that here we have $s_{\mu^{(k)}}=s^{\tilde{\A}}_{\hat{\mu}^*_{\tilde{\A}}}$. Then, invoking \eqref{eq:faastars-6}, we rewrite \eqref{eq:043} in our case as

\begin{equation} \label{eq:faastars-7}
\mathbb{E}(r_{k+1}^2) \leq r_k^2-2 \hat{h}_{\tilde{\A}} \langle \nabla f_\mu (x^{(k)}), x^{(k)} -x^*_{\tilde{\A}} \rangle+\hat{h}_{\tilde{\A}}^2\left(2(\tilde{j}+4)||\nabla f(x^{(k)}) ||^2 + C_2\right),
\end{equation}

\noindent where we recall $C_2=\frac{3K_1+K_2}{\sqrt{2}}  \hat{L_1}\hat{\sigma}\sqrt{\tilde{j}(\tilde{j}+6)^3}$ here. Now, again, equations \eqref{eq:044} through \eqref{eq:048} are identical for this proof as long as $x^*_{\tilde{\A}}$ replaces $x^*$ and $\hat{h}_{\tilde{\A}}$ replaces $\hat{h}$ throughout, and similarly for equations \eqref{eq:049} through \eqref{eq:051} with the additional needed replacement of $\tilde{j}$ for $P$ and using $C_2$ as we have defined it here.
Taking $C_3=\hat{h}_{\tilde{\A}} \hat{\mu}^*_{\tilde{\A}} L_1 \tilde{j} + \hat{h}_{\tilde{\A}}^2 C_2$, \eqref{eq:052} holds (with the usual replacments) and plugging $\hat{h}_{\tilde{\A}}$ in \eqref{eq:faastars-7}:

\begin{equation} \label{eq:faastars-8}
\mathbb{E}(r_{k+1}^2) \leq r_k^2 -\frac{1}{4L_1(\tilde{j}+4)}(f(x^{(k)})-f(x^*_{\tilde{\A}})) + C_3.
\end{equation}

Now plugging $\hat{h}_{\tilde{\A}}$ into $C_3$, we obtain $C_3\leq \frac{\sqrt{K_1}}{\sqrt{K_2}}\cdot \frac{\sigma}{\sqrt{2}L_1} \left(\sqrt{K_1}\cdot0.036+\frac{3K_1+K_2}{16}\cdot 1.034\right)=:C_4$, so

\begin{equation} \label{eq:faastars-9}
\mathbb{E}(r_{k+1}^2) \leq r_k^2 -\frac{1}{4L_1(\tilde{j}+4)}(f(x^{(k)})-f(x^*_{\tilde{\A}})) + C_4.
\end{equation}

We may now apply the expectation over $\mathcal{P}_k$ and $\mathcal{U}_k^{\tilde{\A}}$, $k=M_\A,\ldots, M$, rearrange, and sum over $k=M_\A,\ldots, M$, similarly to before. We obtain a modification to \eqref{eq:060}
with

\begin{equation} \label{eq:faastars-10}
\sum_{k=M_\A}^M\frac{\phi_k-f^*_{\tilde{\A}}}{M-M_\A+1} \leq   \frac{4L_1(\tilde{j}+4)}{(M-M_\A+1)}\left(r_{M_\A}-\mathbb{E}_{\mathcal{U}_k,\mathcal{P}_k}(r_{k+1}^2)\right)+4L_1(\tilde{j}+4)C_4.
\end{equation}

We drop the strictly negative term (again involving $\mathbb{E}_{\mathcal{P}_k,\mathcal{U}_k^{\tilde{\A}}}$) and plug in the definition of $r_{M_\A} =||x^{(M_\A)}-x^*_{\tilde{\A}}||^2$.
Then, recalling $x^*_{\tilde{\A}}=P_{\tilde{\A}}(x^*)+P_{\tilde{\I}}(x^{(M_\A)})$, writing $x^{(M_\A)}=P_{\tilde{\A}}(x^{(M_\A)})+P_{\tilde{\I}}(x^{(M_\A)})$, and noting that $P_{\tilde{\A}}$ is linear, we have $r_{M_\A}=||P_{\tilde{\A}}(x^{(M_\A)}-x^*)||^2$.
Plugging in the value of $C_4$, we obtain \eqref{eq:070}. $\blacksquare$

FAASTARS Corollary 2 shows that its iterates during its third phase (approximate ASTARS) converge to $f^*_{\tilde{\A}}$; however, we need a result for the convergence to $f^*$ like we did in the last section.
In particular, we need a result analogous to ASTARS Corollary 4, generally regarding the distance between evaluations of $f$ and $f_{\tilde{\A}}$ (instead of $f_\A$) at their respective minimizers, $x^*$ and $x^*_{\tilde{\A}}$.
Recall, by $f_{\tilde{\A}}$, we mean $f_{\A}(\lambda):=f(V_{\tilde{\A}} V_{\tilde{\A}}^\top \lambda)$, where the application of $ V_{\tilde{\A}} V_{\tilde{\A}}^\top$ is a linear transformation of $\lambda$ into $\tilde{\A}$, as in \cite{ConstantineK}, which we again rely on heavily for the following result.

\textbf{FAASTARS Corollary 3:} \textit{Let $x^{(M_\A)}$ denote any initial iterate for phase 3 of FAASTARS. 
Let $\A$ continue to denote the true $j$-dimensional AS of $f$ and let $\tilde{\A}$ denote the approximate $\tilde{j}$-dimensional AS of $f$ and the true minimizer of $f_{\tilde{\A}}$ with $x^*_{\tilde{\A}}$ and corresponding noiseless function evaluation $f^*_{\tilde{\A}}$. 
Assume $||x^*-x^*_{\tilde{\A}}||_\Lambda < \infty$, where $||\cdot||_\Lambda$ denotes a norm.
Assume $\tilde{j}=j$ and for $\delta>0$, $||V-\tilde{V}||_2<\delta$, where $||\cdot ||_2$  here denotes the matrix 2-norm induced by the 2-norm in $\Lambda$.
Also, assume the sign of $(\tilde{V}_{\tilde{\I}})_i $, the $i$-th column of $\tilde{V}_{\tilde{\I}}$ to be chosen so that $||(\tilde{V}_{\tilde{\I}})_i - (V_\I)_i||_2 $ is minimized for $i=j+1,\ldots,P$, where $||\cdot||_2$ denotes the vector 2-norm.
Then the difference between $|f^*-f^*_{\tilde{\A}}|$ is bounded by}

\begin{equation} \label{eq:072}
|f^*-f^*_{\tilde{\A}}| \leq \sqrt{\tilde{a}_1}\left(\delta \sqrt{q_1 +\cdots + q_{j}} + \sqrt{q_{j+1}+ \cdots + q_{P}} \right),
\end{equation}

\noindent \textit{where $0 \leq \tilde{a}_1 < \infty$ is an eigenvalue of the positive-semi definite matrix $(x^*-x^*_{\tilde{\A}})(x^*-x^*_{\tilde{\A}})^\top$ and the $q$'s are our notations for the exact eigenvalues associated with the eigendecomposition of the sensitivity matrix $W$. Also, $\tilde{a}_1 = ||P_{\tilde{\I}}(x^{(M_\A)}-x^*)||_\Lambda^2. $}

\textit{Proof:} We bound the quantity $(f^*-f^*_{\tilde{\A}})^2$.
First, we expand $f*=f(x^*)$ around $x^*_{\tilde{\A}}$ by applying \textit{Extended Mean Value Theorem} analogously to ASTARS Corollary 4.
For $c \in [0,1]$ and $\tilde{z} := c x_{\tilde{\A}}^* + (1-c)x^*$ we have $f^*=f(x^*)=f(x_{\tilde{\A}}^*)+ \nabla f (\tilde{z})^\top (x^* - x^*_{\tilde{\A}})  $.
Note that since the components of $x_{\tilde{\A}}^*$ and rotated $V_{\tilde{\A}} V_{\tilde{\A}}^\top x^*$ must match for indices $i=1,\ldots,j$, the point $\tilde{z}$ varies along $\tilde{\I}$ only; that is, $\tilde{z}$ is fixed in $\tilde{\A}$.
Thus, $\nabla f (\tilde{z}) = \nabla_{\tilde{\I}} f(\tilde{z})$, the gradient taken in the approximate inactive subspace only.

The expansion of $f^*$ around $x^*_{\tilde{\A}}$ writes $(f^*-f^*_{\tilde{\A}})^2=\left(\nabla_{\tilde{\I}} f (\tilde{z})^\top (x^* - x^*_{\tilde{\A}})\right)^2= \nabla_{\tilde{\I}} f (\tilde{z})^\top \tilde{A} \nabla_{\tilde{\I}} f (\tilde{z})$, where $\tilde{A}:=(x^* - x^*_{\tilde{\A}})(x^* - x^*_{\tilde{\A}})^\top$ is a $P \times P$ matrix.
Note that $\tilde{A}$ is a square, positive semi-definite, rank 1 matrix.
Hence, it has 1 eigenvalue that is positive or zero, which we denote with  $\tilde{a}_1 \geq 0$; all other $P-1$ eigenvalues are $0$.
We find $\tilde{a}_1 = ||x^* - x^*_{\tilde{\A}} ||^2_\Lambda $ so by definition, $\tilde{a}_1 =  ||P_{\tilde{\I}}(x^{(M_\A)}-x^*)||_\Lambda^2.$
Observe $\tilde{a}_1<\infty$ because $||x^*-x^*_{\tilde{\A}}||_\Lambda < \infty$ and also $\tilde{a}_1=0 \iff x^* = x^*_{\tilde{\A}}$, which case $f^*=f^*_{\tilde{\A}}$ so that $|f^*-f^*_{\tilde{\A}}|=0$.
Then $(f^*-f^*_{\tilde{\A}})^2 \leq \tilde{a}_1 \nabla_{\tilde{\I}} f (\tilde{z})^\top \nabla_{\tilde{\I}} f (\tilde{z})$.

Recall $\tilde{j}=j$ and for $\delta>0$, $||V-\tilde{V}||_2<\delta$.
Then we have \textit{comparable partitions}, in the sense that the submatrices $V_\A$ and $\tilde{V}_{\tilde{\A}}$ (of $V$ and $\tilde{V}$ respectively) are both $j$-dimensional, and likewise the submatrices $V_\I$ and $V_{\tilde{\I}}$ are both $P-j$-dimensional.
We shall need the results from Lemma 3.4 in \cite{ConstantineK} which state $|| V_{\I}^\top \tilde{V}_{\tilde{\I}} ||_2 \leq 1$ and $|| V_{\A}^\top \tilde{V}_{\tilde{\I}} ||_2 \leq \delta$, where $||\cdot||_2$ denotes the matrix 2-norm.
Note that Lemma 3.4 requires that the sign of $(\tilde{V}_{\tilde{\I}})_i $, the $i$-th column of $\tilde{V}_{\tilde{\I}}$ to be chosen so that $||(\tilde{V}_{\tilde{\I}})_i - (V_\I)_i||_2 $ is minimized for $i=j+1,\ldots,P$, where $||\cdot||_2$ denotes the vector 2-norm.
Now the chain rule provides $\nabla_{\tilde{\I}} f=V_{\I}^\top \tilde{V}_{\tilde{\I}} \nabla_{\I} f + V_{\A}^\top \tilde{V}_{\tilde{\I}} \nabla_{\A} f$ (\cite{ConstantineK}, pp. A1510).

Applying the expectation over both $\I$ and $\A$ to both sides (where the left-hand side is constant with respect to this expectation) we have $(f^*-f^*_{\tilde{\A}}) \leq \tilde{a}_1 \mathbb{E}_{\A,\I} \left(\nabla_{\tilde{\I}} f (\tilde{z})^\top \nabla_{\tilde{\I}} f (\tilde{z})\right)$.
Using the chain rule and $|| V_{\I}^\top \tilde{V}_{\tilde{\I}} ||_2 \leq 1$ and $|| V_{\A}^\top \tilde{V}_{\tilde{\I}} ||_2 \leq \delta$ gives

\begin{equation} \label{eq:faastars-11}
(f^*-f^*_{\tilde{\A}})^2 \leq \tilde{a}_1\left( \mathbb{E}_{\I}\left(\nabla_\I f(\tilde{z})^\top \nabla_\I f(\tilde{z})\right) +2 \delta \, \mathbb{E}_{\A,\I}\left( \nabla_\I f(\tilde{z})^\top \nabla_\A f(\tilde{z}) \right) + \delta^2 \mathbb{E}_{\A} \left( \nabla_\A f(\tilde{z})^\top \nabla_\A f(\tilde{z})  \right)\right),
\end{equation}

\noindent where linearity of $\mathbb{E}$ allows for breaking the expectation of the appearing sum into a sum of expectations taken over only $\A$, only $\I$, or both $\A$ and $\I$, depending on which types of gradients appear.
The Cauchy-Schwarz inequality may be applied to the second term on the right-hand side of \eqref{eq:faastars-11} to write 

\begin{equation} \label{eq:faastars-12}
2 \delta \mathbb{E}_{\A,\I}\left( \nabla_\I f(\tilde{z})^\top \nabla_\A f(\tilde{z}) \right) \leq 2 \delta \, \mathbb{E}_{\I}\left( \nabla_\I f(\tilde{z})^\top \nabla_\I f(\tilde{z}) \right) \mathbb{E}_{\A} \left( \nabla_\A f(\tilde{z})^\top \nabla_\I f(\tilde{z})) \right),
\end{equation}

\noindent where we can split up the two terms inside $\mathbb{E}_{\A,\I}(\cdot)$ into two separate expectations taken over $\I$ and $\A$ individually, due to their independence. 
Substituting this bound into \eqref{eq:hot-fix-11} and then factoring the resulting terms will yield

\begin{equation} \label{eq:faastars-13}
(f^*-f^*_{\tilde{\A}})^2 \leq \tilde{a}_1 \left( \mathbb{E}_{\I}\left( \nabla_\I f(\tilde{z})^\top \nabla_\I f(\tilde{z}) \right)^{1/2} + \delta \, \mathbb{E}_{\A}\left( \nabla_\A f(\tilde{z})^\top \nabla_\A f(\tilde{z}) \right)^{1/2}    \right)^2,
\end{equation}

Citing \cite{ConstantineK} Lemma 2.2, we have $\mathbb{E}_{\A}\left(\nabla_\A f (\tilde{z})^\top \nabla_\A f (\tilde{z})\right)<q_{1}+\cdots + q_j$, where $q_i$, $i=1,\ldots,j$ are the first $j$ eigenvalues of the sensitivity matrix $W$ and $\mathbb{E}_{\I}\left(\nabla_\I f (\tilde{z})^\top \nabla_\I f (\tilde{z})\right)<q_{j+1}+\cdots + q_P$, where $q_i$, $i=j+1,\ldots,P$ are the last $P-j$ eigenvalues of the sensitivity matrix $W$.
Hence, $(f^*-f^*_\A)^2 \leq \tilde{a}_1 \left( \delta \, \sqrt{q_1 + \cdots + q_j} + \sqrt{q_{j+1}+\cdots + q_P} \right)^2$.
Applying a square root to both sides, we obtain \eqref{eq:072}. $\blacksquare$

We now present a statement about the convergence of FAASTARS to $f^*$ by combining the previous two results, along with STARS Theorem 4.5 (Modified).

Recall for $k=1,\ldots,M_\A$ iterates correspond to STARS with estimated hyperparameters and for $k=M_\A+1, \ldots, M$ the iterates correspond to ASTARS with estimated hyperparameters and estimated $\A$.

\textbf{FAASTARS Theorem 4:}
\textit{For $k=M_\A, \ldots, M$, let the assumptions from FAASTARS Corollaries 3 and 4 hold. For any total number of FAASTARS iterations $M\geq 1$,}

\begin{equation} \label{eq:074}
\begin{split}
\sum_{k=M_\A}^M\frac{\phi_k-f^*}{M-M_\A+1}
\leq   
& \frac{4L_1(j+4)||P_{\tilde{\A}}(x^{(M_\A)}-x^*)||^2}{\sqrt{K_1}(2-\sqrt{K_1})(M-M_\A+1)} + \frac{4\sigma(j+4)}{\sqrt{2K_2}(2-\sqrt{K_1})}C_5 \\
&+ \sqrt{\tilde{a}_1}\left(\delta \sqrt{q_1 +\cdots + q_{j}} + \sqrt{q_{j+1}+ \cdots + q_{P}} \right),
\end{split}
\end{equation}

\noindent where $C_5:=\sqrt{K_1}\cdot0.036+\frac{3K_1+K_2}{16}\cdot 1.034$.

\textit{Proof:} Using the triangle inequality, we have $|\phi_k-f^{*}|\leq |\phi_k-f^*_{\tilde{\A}}|+ |f^*-f^*_{\tilde{\A}}|.$ Thus, noting that the quantity $\phi_k-f^*_{\tilde{\A}}$ is always nonnegative, we use the triangle inequality to rewrite the summation on the left-hand side of \eqref{eq:074} as

\begin{equation} \label{eq:075}
\sum_{k=M_\A}^M\frac{\phi_k-f^*}{M-M_\A+1} 
\leq 
\sum_{k=M_\A}^M\frac{\phi_k-f^*_{\tilde{\A}}}{M-M_\A+1}+ |f^*-f^*_{\tilde{\A}}|.
\end{equation}

Now the second sum on the right-hand side of \eqref{eq:076} is in a more useful form for us, since the final phase of iterates, $k=M_\A,\ldots, M$ are (approximate) ASTARS iterates converging to $f^*_{\tilde{\A}}$ (not $f^*$); thus, we recall that we already proved a bound for this sum in FAASTARS Corollary 2.
Note that since we assume in FAASTARS Corollary 3 that $\tilde{j}=j$, we replace $\tilde{j}$ with $j$ in the result from FAASTARS Corollary 2.
Otherwise, we just invoke FAASTARS Corollary 3, bounding the last term on the right-hand side of \eqref{eq:075} to obtain \eqref{eq:074}. $\blacksquare$

We use the results above to analyze the complexity of FAASTARS in the following remark.

\textbf{Remark:} We again mimic the complexity analyses we performed in the preceding sections for FAASTARS.
Let $||\cdot||$ denote a norm in $\Lambda$ throughout.
Define $R_{\tilde{\A}}^2$ as a bound $||P_{\tilde{\A}}(x^{(M_\A)}-x^*)||^2 \leq R_{\tilde{\A}}^2$.
We recall $\tilde{a}_1 = ||P_{\tilde{\I}}(x^{(M_\A)}-x^*)||_\Lambda^2$ and define a bound $\tilde{a}_1 \leq R_{\tilde{\I}}^2$.

Define $x^\dagger:= \text{argmin}_x\{f(x):x \in \{x^{(0)}, \ldots, x^M\} \}$ and $\phi_\dagger:=\mathbb{E}_{\mathcal{U}_{k-1},\mathcal{P}_{k-1}}(f(x^\dagger))$ as before. Then the value $\phi_\dagger-f^*$ must be less than or equal to the average improvement for any given run of STARS; that is, \eqref{eq:074}, along with our new definitions, implies that

\begin{equation} \label{eq:076}
\begin{split}
\phi_\dagger-f^* 
\leq   
& \frac{4L_1(j+4)}{\sqrt{K_1}(2-\sqrt{K_1})(M-M_\A+1)}R_{\tilde{\A}}^2 + \frac{4\sigma(j+4)}{\sqrt{2K_2}(2-\sqrt{K_1})}C_5 \\
&+ R_{\tilde{\I}} \left(\delta \sqrt{q_1 +\ldots + q_{j}} + \sqrt{q_{j+1}+ \cdots + q_{P}} \right),
\end{split}
\end{equation}

Now in this analysis, we are only taking enough approximate STARS steps to learn a surrogate for $f$ to obtain $\tilde{\A}$
Thus, we have $M_\A \sim \mathcal{O}(L(P))$, where $L$ is a function of $P$ that depends on the surrogate method used to learn $\A$.
For instance, if we use a linear surrogate or RBF's  -- which begin by fitting a linear surrogate, and then later a quadratic surrogate once enough ASTARS steps are taken -- then $L(P)=P$.
If we use quadratic surrogates, then $L(P)=P^2$.
(Typically we use quadratic surrogates for higher-quality active subspaces.)

The terms involving Phase 3 of FAASTARS (approximate ASTARS phase) are the three terms on the right-hand side of \eqref{eq:076}.
We assume that we wish to achieve a final accuracy of $\epsilon_{\text{tol}}>0$. Then we will need $\phi_\dagger-f^* \leq \epsilon_{\text{tol}}$. If we take 

\begin{equation} \label{eq:end-complexity-1}
\frac{4\sigma(j+4)}{\sqrt{2K_2}(2-\sqrt{K_1})}C_5 \leq \frac{\epsilon_{\text{tol}}}{3},
\end{equation}

\noindent then we must require that the noise not exceed the following threshold:

\begin{equation} \label{eq:end-complexity-2}
\sigma \leq \frac{\sqrt{2K_2}(2-\sqrt{K_1})\epsilon_{\text{tol}}}{12\sigma(j+4)}.
\end{equation}

\noindent If we satisfy \eqref{eq:end-complexity-2} and also have

\begin{equation} \label{eq:end-complexity-3}
R_{\tilde{\I}}  \left(\delta \sqrt{q_1 +\ldots + q_{j}} + \sqrt{q_{j+1}+ \cdots + q_{P}} \right) \leq \frac{\epsilon_{\text{tol}}}{3},
\end{equation}

\noindent then we can achieve $\epsilon_{\text{tol}}$ accuracy as long as 

\begin{equation} \label{eq:end-complexity-4}
\frac{4L_1(j+4)}{\sqrt{K_1}(2-\sqrt{K_1})(M-M_\A+1)}R_{\tilde{\A}}^2 \leq \frac{\epsilon_{\text{tol}}}{3} \iff M \geq M_{\A} + \frac{12 L_1(j+4)R_{\tilde{\A}}^2}{\epsilon_{\text{tol}} \sqrt{K_1}(2-\sqrt{K_1})}-1.
\end{equation}

Hence, we achieve $\epsilon_{\text{tol}}$ accuracy as long as: the noise is small enough; the eigenvalues of the inactive subspace and distance from $x^*$ to $x^*_{\tilde{\A}}$ is small enough; and $M$ is large enough.
Details of those required bounds given by \eqref{eq:end-complexity-2}, \eqref{eq:end-complexity-3}, and \eqref{eq:end-complexity-4}, respectively. With these assumptions, we achieve $\epsilon_{\text{tol}}$ accuracy in

\begin{equation} \label{eq:end-complexity-5}
M \sim \mathcal{O}\left( \max \left\{ L(P), \, \frac{L_1 j R_{\tilde{\A}}^2}{\epsilon_{\text{tol}}\sqrt{K_1}(2-\sqrt{K_1})} \right\} \right).
\end{equation}

This analysis also shows that given a particular variance in the noise, $\sigma$, as well as the term 
involving the eigenvalues of the inactive subspace and distance from $x^*$ to $x^*_{\tilde{\A}}$, the achievable accuracy can be no better (i.e., less) than the value given below:

\begin{equation} \label{eq:end-complexity-6}
\epsilon_{\text{tol}} \geq \max \left\{ \frac{12(j+4)C_5}{\sqrt{K_1}(2-\sqrt{K_1})} \sigma, \, 3 R_{\tilde{\I}}  \left(\delta \sqrt{q_1 +\cdots + q_{j}} + \sqrt{q_{j+1}+ \cdots + q_{P}} \, \right)  \right\},
\end{equation}

\noindent showing that the achievable accuracy will either be limited (mainly) by hyperparameter approximations or (mainly) by the error in $\tilde{\A}$.

We find similar complexity results to that of ASTARS, but also pay a price for approximations to $\tilde{\A}$, especially when the approximation is poor, usually do to both an insufficient number of samples and samples that do not explore $\Lambda$ sufficiently.
(Recall that since $\tilde{\A}$ is formed using information from a surrogate trained on STARS samples in the FAASTARS routine, it is really poor \textit{surrogates} that hurt convergence.)